\input amstex
\documentstyle{amsppt}
\nologo
%\NoRunningHeads
\TagsOnRight
\input pictex
%\magnification 1200
\parskip=\medskipamount
\NoBlackBoxes

\def\ve{\varepsilon}

\hsize 6.0 true in
\vsize=8.5 true in
\def\nint{\mathop{\diagup\kern-13.0pt\int}}
\def\textnoint{\matpop
{\raise1.25pt\hbox{${\ssize\diagup}$}\kern-9.5pt\int}}

\lineskip=2pt\baselineskip=18pt\lineskiplimit=0pt

\topmatter
\title
Decoupling for perturbed cones and mean square of $|\zeta (\frac 12+it)|$
\endtitle
\author {J. Bourgain and N.~Watt}
\endauthor
\address {Institute for Advanced Study, School of Mathematics, Princeton, NJ 08540}
\endaddress
\address
{Dunfermline}
\endaddress
\thanks {First author was partially supported by NSF grant DMS 1301619}
\endthanks

\abstract
An improved estimate is obtained for the mean square of the modulus of the zeta function on the critical line.
It is based on the decoupling techniques in harmonic analysis developed in 
[B-D].
\endabstract

\endtopmatter

\document

\beginsection
{1. Introduction}

The aim of this paper is to establish improved bounds for the mean square of $|\zeta(\frac 12+it)|$ on short intervals.

More precisely, we will prove

\noindent
{\bf Theorem 3.}
{\sl Let the function $I:[0, \infty)\times (0,\infty)\to\Bbb R$ be given by
$$
I(t, U) =\frac 1{2U} \int_{t-U}^{t+U} \Big|\zeta \Big(\frac 12+i\tau\Big)\Big|^2 d\tau.
$$
Then, for all $\ve>0$, one has
$$
I\Big(t, t^{\frac {1273}{4053}+\ve}\Big)= O(\log t)\ \text{ as $t\to\infty$}.
$$ }
\noindent
and

\noindent
{\bf Theorem 4.}
{\sl Define $E(T)$ by
$$
E(T) =\int^T_0\Big|\zeta\Big(\frac 12+it\Big)\Big|^2 dt -\Big(\log \Big(\frac T{2\pi}\Big)+2\gamma-1\Big)T \quad (T\geq 1)
$$
($\gamma$ = Euler-Mascheroni constant).
Then for all $\ve>0$
$$
E(T) =O(T^{\frac {1515}{4816}+\ve}) \ \text { as } \ T\to\infty.
$$}

Theorem 4 improves upon the estimate $E(T)= O(T^{\frac {131}{416}}(\log T)^{\frac{32587}{8320}})$ obtained in \cite{W10}, noting
that
$$
\frac {131}{416} = 0, 314903... \ \text { while \ } \frac {1515}{4816}= 0,314576\ldots
$$
(see the Remarks in Section 13 for a more detailed dimension).

The basic approach is the Bombieri-Iwaniec method and its further developments, in particular the contributions of M.~Huxley and
the second author.

Recall that there are two main parts to the Bombieri-Iwaniec approach, referred to as the First and Second Spacing Problem.
Roughly speaking, the content of the first spacing problem are certain moment inequalities while the second spacing problem is a
distributional issue.
These two components are then combined by an application of the large sieve.
See \cite{H96}.
A novelty in this work is a different treatment of the first spacing problem using recent developments around the `decoupling
principle' in harmonic analysis (see \cite {B-D}).
An earlier application of these results on bounding $\zeta\big(\frac 12+it\big)$, i.e. towards the Lindel\"of hypothesis, appears
in the first author's paper \cite {B}.
                                                                                                                  
While the original Bombieri-Iwaniec method deals with one - variable exponential sums, the present context involves exponential
sums with two variables (see Theorem 2 in Section 6).
As a consequence, the mean value estimates in the first spacing problem involves points on surfaces rather than curves.
These surfaces turn out to be perturbed cones.
Our improvement in treating the first spacing problem relies on exploiting the additional curvature in the radial direction and
moment inequalities for $q> 4$.

Next some more details.

Let $k\sim K, \ell\sim L, K>L$ and
$$
\omega (k, \ell)= \frac{(k+\ell)^{3/2} -(k-\ell)^{3/2}}{3} = k^{\frac
12}\ell+ck^{-\frac 32}\ell^3+\cdots \ (c\not= 0).
\eqno{(1.1)}
$$
Let $\eta>0$ be a small parameter.
Motivated by the first spacing problem, we are interested in bounding
moments
$$
\Big\Vert\sum_{k\sim K, \ell \sim L} a_{k\ell} \, e(\ell{x_1}+k\ell
x_2+\omega (k, \ell)x_3)\Big\Vert_{L^q_{\#} (|x_1|< 1, |x_2|<1, |x_3|<\frac
1{\eta L\sqrt
K})}\eqno{(1.2)}
$$
where $|a_{k, \ell}|\leq 1$ and $L^q_{\#}$ refers to the averaged
$L^q$-norm.
The case $q=4$ corresponds to the `classical' treatment (cf. \cite {H96})
and we aim at $q>4$ (not necessarily an integer) in order to achieve
better
estimates when exploiting the large sieve (which basically amounts to an
application of H\"older's inequality).
It is important to point out that this improvement uses essentially the
perturbative terms $ck^{-3/2}\ell^3+\cdots $ in (1.1) and that we are
unable to establish a similar result for the case $\omega(k, \ell)= k^{3/2} \ell$ of the
unperturbed cone.

Setting $s=\frac kK, t=\frac \ell L$, we consider the surface
$$
\tilde C:\cases
z=t\\ y = st\\
x= s^{\frac 12} t -c \ve^{2} s^{-\frac 32} t^2+\cdots\endcases
\eqno{(1.3)}
$$
with $\ve=\frac LK$, $s, t\sim 1$.

Thus $\tilde C$ parametrizes a perturbed cone and our approach to (1.2)
consists in invoking as initial step the so called `decoupling theory'
from
\cite{B-D}.
The role of this step is  to achieve a variable restriction, after which
we again exploit further arithmetical considerations as in earlier treatments.
But because of restrictions of the variables $k$ and $\ell$ to suitably
small intervals this arithmetical component becomes more
straightforward.
\medskip

In the next four sections, we prove the main analytical inequalities (see Propositions
10 and 10$'$)
needed for our enhanced treatment of the first spacing problem.
We will use without much explanation several techniques
and results from modern harmonic analysis in the presence of curvature.
The decoupling result for curved (hyper) surfaces (applied here to
surfaces in $\Bbb
R^3$) will be an essential ingredient.

\bigskip
\beginsection
{ 2. Decoupling inequalities for the cone}

Consider the truncated cone $C$ and let $C_{\frac 1N}$ be a $\frac
1N$-neighborhood of $C$

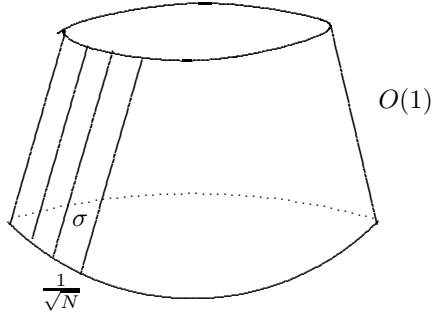
\font\thinlinefont=cmr5
$$
\hbox{\beginpicture
\setcoordinatesystem units <.80000cm,.80000cm>
\linethickness=1pt
{\circulararc 90.690 degrees from  1.175 20.955 center at  4.239 24.014
}%
%
% Fig CIRCULAR ARC object
%
\linethickness= 0.500pt
\setplotsymbol ({\thinlinefont .})
\setdots < 0.0953cm>
{\circulararc 31.373 degrees from  7.239 21.018 center at  4.239 10.335
}%
%
% Fig POLYLINE object
%
\linethickness= 0.500pt
\setplotsymbol ({\thinlinefont .})
\setsolid
{\plot  2.127 24.130  1.206 20.923 /
}%
%
% Fig POLYLINE object
%
\linethickness= 0.500pt
\setplotsymbol ({\thinlinefont .})
{\plot  6.541 24.194  7.303 20.923 /
}%
%
% Fig POLYLINE object
%
\linethickness= 0.500pt
\setplotsymbol ({\thinlinefont .})
{\plot  3.397 23.654  2.381 20.098 /
}%
%
% Fig POLYLINE object
%
\linethickness= 0.500pt
\setplotsymbol ({\thinlinefont .})
{\plot  2.910 23.793  1.926 20.396 /
}%
\linethickness= 0.500pt
\setplotsymbol ({\thinlinefont .})
{\plot  2.496 23.889  1.575 20.682 /
}%
%
% Fig POLYLINE object
%
\linethickness= 0.500pt
\setplotsymbol ({\thinlinefont .})
{\plot  2.000 24.098  2.004 24.100 /
\plot  2.004 24.100  2.013 24.107 /
\plot  2.013 24.107  2.030 24.117 /
\plot  2.030 24.117  2.051 24.132 /
\plot  2.051 24.132  2.081 24.149 /
\plot  2.081 24.149  2.115 24.168 /
\plot  2.115 24.168  2.153 24.189 /
\plot  2.153 24.189  2.191 24.208 /
\plot  2.191 24.208  2.231 24.227 /
\plot  2.231 24.227  2.273 24.244 /
\plot  2.273 24.244  2.318 24.263 /
\plot  2.318 24.263  2.366 24.280 /
\plot  2.366 24.280  2.419 24.297 /
\plot  2.419 24.297  2.477 24.314 /
\plot  2.477 24.314  2.540 24.331 /
\plot  2.540 24.331  2.580 24.342 /
\plot  2.580 24.342  2.623 24.352 /
\plot  2.623 24.352  2.667 24.363 /
\plot  2.667 24.363  2.714 24.376 /
\plot  2.714 24.376  2.762 24.386 /
\plot  2.762 24.386  2.815 24.399 /
\plot  2.815 24.399  2.870 24.412 /
\plot  2.870 24.412  2.929 24.424 /
\plot  2.929 24.424  2.991 24.437 /
\plot  2.991 24.437  3.054 24.450 /
\plot  3.054 24.450  3.120 24.460 /
\plot  3.120 24.460  3.188 24.473 /
\plot  3.188 24.473  3.258 24.486 /
\plot  3.258 24.486  3.330 24.498 /
\plot  3.330 24.498  3.401 24.509 /
\plot  3.401 24.509  3.476 24.519 /
\plot  3.476 24.519  3.550 24.530 /
\plot  3.550 24.530  3.624 24.541 /
\plot  3.624 24.541  3.698 24.549 /
\plot  3.698 24.549  3.772 24.558 /
\plot  3.772 24.558  3.846 24.566 /
\plot  3.846 24.566  3.918 24.572 /
\plot  3.918 24.572  3.990 24.579 /
\plot  3.990 24.579  4.062 24.585 /
\plot  4.062 24.585  4.132 24.589 /
\plot  4.132 24.589  4.202 24.591 /
\plot  4.202 24.591  4.271 24.594 /
\plot  4.271 24.594  4.339 24.596 /
\putrule from  4.339 24.596 to  4.407 24.596
\putrule from  4.407 24.596 to  4.477 24.596
\putrule from  4.477 24.596 to  4.547 24.596
\plot  4.547 24.596  4.616 24.594 /
\plot  4.616 24.594  4.686 24.589 /
\plot  4.686 24.589  4.758 24.585 /
\plot  4.758 24.585  4.830 24.581 /
\plot  4.830 24.581  4.902 24.577 /
\plot  4.902 24.577  4.976 24.570 /
\plot  4.976 24.570  5.048 24.564 /
\plot  5.048 24.564  5.122 24.555 /
\plot  5.122 24.555  5.194 24.547 /
\plot  5.194 24.547  5.266 24.539 /
\plot  5.266 24.539  5.336 24.530 /
\plot  5.336 24.530  5.406 24.519 /
\plot  5.406 24.519  5.474 24.511 /
\plot  5.474 24.511  5.541 24.500 /
\plot  5.541 24.500  5.605 24.490 /
\plot  5.605 24.490  5.666 24.479 /
\plot  5.666 24.479  5.726 24.469 /
\plot  5.726 24.469  5.783 24.460 /
\plot  5.783 24.460  5.836 24.450 /
\plot  5.836 24.450  5.886 24.439 /
\plot  5.886 24.439  5.935 24.431 /
\plot  5.935 24.431  5.980 24.420 /
\plot  5.980 24.420  6.024 24.412 /
\plot  6.024 24.412  6.064 24.403 /
\plot  6.064 24.403  6.102 24.395 /
\plot  6.102 24.395  6.164 24.380 /
\plot  6.164 24.380  6.219 24.367 /
\plot  6.219 24.367  6.270 24.354 /
\plot  6.270 24.354  6.316 24.342 /
\plot  6.316 24.342  6.358 24.329 /
\plot  6.358 24.329  6.394 24.316 /
\plot  6.358 24.329  6.394 24.316 /
\plot  6.394 24.316  6.428 24.304 /
\plot  6.428 24.304  6.456 24.293 /
\plot  6.456 24.293  6.477 24.280 /
\plot  6.477 24.280  6.496 24.270 /
\plot  6.496 24.270  6.511 24.259 /
\plot  6.511 24.259  6.519 24.249 /
\plot  6.519 24.249  6.526 24.238 /
\plot  6.526 24.238  6.530 24.229 /
\putrule from  6.530 24.229 to  6.530 24.219
\putrule from  6.530 24.219 to  6.530 24.208
\plot  6.530 24.208  6.528 24.200 /
\plot  6.528 24.200  6.524 24.189 /
\plot  6.524 24.189  6.515 24.177 /
\plot  6.515 24.177  6.505 24.164 /
\plot  6.505 24.164  6.492 24.151 /
\plot  6.492 24.151  6.473 24.136 /
\plot  6.473 24.136  6.452 24.122 /
\plot  6.452 24.122  6.424 24.105 /
\plot  6.424 24.105  6.394 24.088 /
\plot  6.394 24.088  6.358 24.071 /
\plot  6.358 24.071  6.318 24.052 /
\plot  6.318 24.052  6.272 24.035 /
\plot  6.272 24.035  6.223 24.016 /
\plot  6.223 24.016  6.170 23.997 /
\plot  6.170 23.997  6.111 23.975 /
\plot  6.111 23.975  6.049 23.954 /
\plot  6.049 23.954  6.009 23.944 /
\plot  6.009 23.944  5.969 23.931 /
\plot  5.969 23.931  5.925 23.918 /
\plot  5.925 23.918  5.878 23.906 /
\plot  5.878 23.906  5.829 23.891 /
\plot  5.829 23.891  5.779 23.878 /
\plot  5.779 23.878  5.723 23.863 /
\plot  5.723 23.863  5.666 23.848 /
\plot  5.666 23.848  5.607 23.834 /
\plot  5.607 23.834  5.544 23.821 /
\plot  5.544 23.821  5.480 23.806 /
\plot  5.480 23.806  5.412 23.791 /
\plot  5.412 23.791  5.345 23.777 /
\plot  5.345 23.777  5.275 23.764 /
\plot  5.275 23.764  5.203 23.749 /
\plot  5.203 23.749  5.131 23.736 /
\plot  5.131 23.736  5.059 23.724 /
\plot  5.059 23.724  4.985 23.713 /
\plot  4.985 23.713  4.913 23.700 /
\plot  4.985 23.713  4.913 23.700 /
\plot  4.913 23.700  4.839 23.692 /
\plot  4.839 23.692  4.767 23.681 /
\plot  4.767 23.681  4.695 23.673 /
\plot  4.695 23.673  4.623 23.666 /
\plot  4.623 23.666  4.553 23.660 /
\plot  4.553 23.660  4.483 23.654 /
\plot  4.483 23.654  4.413 23.650 /
\plot  4.413 23.650  4.343 23.645 /
\plot  4.343 23.645  4.276 23.643 /
\plot  4.276 23.643  4.208 23.641 /
\putrule from  4.208 23.641 to  4.138 23.641
\putrule from  4.138 23.641 to  4.068 23.641
\plot  4.068 23.641  3.998 23.643 /
\plot  3.998 23.643  3.926 23.645 /
\plot  3.926 23.645  3.854 23.647 /
\plot  3.854 23.647  3.782 23.652 /
\plot  3.782 23.652  3.708 23.658 /
\plot  3.708 23.658  3.634 23.662 /
\plot  3.634 23.662  3.560 23.671 /
\plot  3.560 23.671  3.486 23.677 /
\plot  3.486 23.677  3.412 23.685 /
\plot  3.412 23.685  3.340 23.694 /
\plot  3.340 23.694  3.268 23.705 /
\plot  3.268 23.705  3.198 23.715 /
\plot  3.198 23.715  3.128 23.726 /
\plot  3.128 23.726  3.061 23.736 /
\plot  3.061 23.736  2.997 23.747 /
\plot  2.997 23.747  2.934 23.760 /
\plot  2.934 23.760  2.874 23.772 /
\plot  2.874 23.772  2.817 23.783 /
\plot  2.817 23.783  2.762 23.796 /
\plot  2.762 23.796  2.711 23.806 /
\plot  2.711 23.806  2.663 23.819 /
\plot  2.663 23.819  2.618 23.832 /
\plot  2.618 23.832  2.576 23.842 /
\plot  2.576 23.842  2.536 23.855 /
\plot  2.536 23.855  2.498 23.865 /
\plot  2.498 23.865  2.436 23.884 /
\plot  2.436 23.884  2.383 23.906 /
\plot  2.383 23.906  2.335 23.925 /
\plot  2.335 23.925  2.290 23.946 /
\plot  2.290 23.946  2.250 23.965 /
\plot  2.250 23.965  2.216 23.986 /
\plot  2.216 23.986  2.187 24.005 /
\plot  2.187 24.005  2.161 24.024 /
\plot  2.161 24.024  2.140 24.043 /
\plot  2.140 24.043  2.123 24.060 /
\plot  2.123 24.060  2.110 24.075 /
\plot  2.110 24.075  2.102 24.090 /
\plot  2.102 24.090  2.095 24.102 /
\plot  2.095 24.102  2.089 24.113 /
\plot  2.089 24.113  2.087 24.122 /
\plot  2.087 24.122  2.085 24.130 /
\putrule from  2.085 24.130 to  2.085 24.132
\plot  2.085 24.132  2.083 24.136 /
\putrule from  2.083 24.136 to  2.083 24.138
\putrule from  2.083 24.138 to  2.083 24.141
\putrule from  2.083 24.141 to  2.083 24.143
\putrule from  2.083 24.143 to  2.083 24.145
\putrule from  2.083 24.145 to  2.081 24.145
\putrule from  2.081 24.145 to  2.081 24.147
\putrule from  2.081 24.147 to  2.081 24.149
\putrule from  2.081 24.149 to  2.083 24.149
\putrule from  2.083 24.149 to  2.085 24.149
\putrule from  2.085 24.149 to  2.085 24.147
\plot  2.085 24.147  2.087 24.145 /
\putrule from  2.087 24.145 to  2.087 24.143
\putrule from  2.087 24.143 to  2.089 24.143
\putrule from  2.089 24.143 to  2.089 24.141
\plot  2.089 24.141  2.091 24.138 /
\putrule from  2.091 24.138 to  2.091 24.136
\plot  2.091 24.136  2.093 24.134 /
\putrule from  2.093 24.134 to  2.093 24.132
\putrule from  2.093 24.132 to  2.095 24.132
\putrule from  2.095 24.132 to  2.095 24.130
}%
%
% Fig TEXT object
%
\put{$\sigma$ } [lB] at  2.235 20.892
%
% Fig TEXT object
%
\put{$\frac 1{\sqrt N}$} [lB] at  1.700 19.730
%
% Fig TEXT object
%
\put{$O(1)$} [lB] at  7.325 22.860
\linethickness=0pt
\putrectangle corners at  1.158 24.621 and  7.557 19.181
\endpicture}
$$

Here $N$ is some scale (will be taken to be $\sim \frac
1{\ve^2}$ in (1.3) later on).

Partition $C_{\frac 1N}$ into plates $\{\sigma\}$ as indicated above.

Thus each $\sigma$ has dimensions $\sim \frac 1{\sqrt N}\times \frac 2N\times 1$.

\smallskip

Denoting $B_N\subset \Bbb
R^3$ a ball (not necessarily centered at
$0$) of size $N$, the following
decoupling inequality is proven in \cite {B-D}.

\noindent
{\bf Proposition 1.}
{\sl
Assume {\rm supp}$\,\hat f\subset C_{\frac 1N}$ and denote $f_\sigma=(\hat
f|_\sigma)^\vee$ the Fourier restriction of $f$ to $\sigma$.
Then
$$
\Vert f\Vert_{L^6{(B_N)}}\ll
N^{o+} \Big(\sum_\sigma \Vert
f_\sigma\Vert^2
_{L^6(B_N)}\Big)^{\frac 12}\eqno{(2.1)}
$$
where $\ll N^{0+}$ means $\leq C_\ve
N^\ve$ for any
$\ve>0$.}

\smallskip

The next inequality we need is a consequence of the multilinear theory from \cite {BCT}.

Given functions $f_1, f_2, f_3$, supp$\, \hat f_i\subset C_{\frac 1N}$, we say that $f_1, f_2, f_3$
are separated provided supp$\,\hat
f_i$ are contained in slabs that are angularly $O(1)$-separated.

\noindent
{\bf Proposition 2. } {\sl (trilinear inequality)

Let $f_1, f_2, f_3$, supp$\, \hat f_i\subset C_{\frac 1N}$ be separated. Then
$$
\Vert \,
|f_1f_2f_3|^{\frac 13}\Vert
_{L^3_{\#} (B_N)} \ll
N^{o+} \prod^3_{i=1} \Vert
f_i\Vert_{L^2_{\#}(B_N)}^{1/3}.\eqno{(2.2)}
$$}

More generally, assume supp$\,\hat f_i\subset C_{\frac 1N}\cap R_\delta$, where $R_\delta$ is an
angular sector of size $\delta$ and that      moreover
supp$\, \hat f_i$ are $O(\delta)$-separated $(\delta> \frac 1{\sqrt N})$

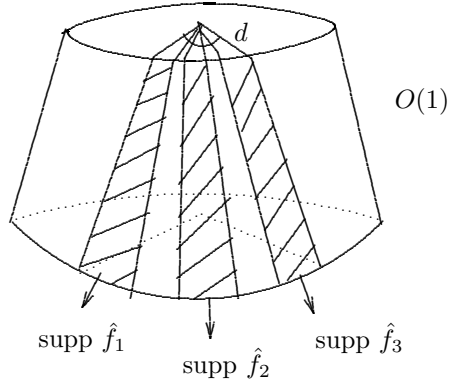
\font\thinlinefont=cmr5
$$
\hbox{\beginpicture
\setcoordinatesystem units <.80000cm,.80000cm>
% Fig CIRCULAR ARC object
%
\linethickness= 0.500pt
\setplotsymbol ({\thinlinefont .})
{\circulararc 90.690 degrees from  1.175 20.955 center at  4.239 24.014
}%
%
% Fig CIRCULAR ARC object
%
\linethickness= 0.500pt
\setplotsymbol ({\thinlinefont .})
\setdots < 0.0953cm>
{\circulararc 31.373 degrees from  7.239 21.018 center at  4.239 10.335
}%
%
% Fig CIRCULAR ARC object
%
\linethickness= 0.500pt
\setplotsymbol ({\thinlinefont .})
\setsolid
{\circulararc 120.510 degrees from  4.015 24.083 center at  4.326 24.240
}%
%
% Fig POLYLINE object
%
\linethickness= 0.500pt
\setplotsymbol ({\thinlinefont .})
{\plot  6.541 24.194  7.303 20.923 /
}%
%
% Fig POLYLINE object
%
\linethickness= 0.500pt
\setplotsymbol ({\thinlinefont .})
{\plot  2.055 24.172  1.135 20.966 /
}%
%
% Fig POLYLINE object
%
% Fig POLYLINE object
%
\linethickness= 0.500pt
\setplotsymbol ({\thinlinefont .})
{\plot  5.842 21.558  5.366 20.828 /
}%
%
% Fig POLYLINE object
%
\linethickness= 0.500pt
\setplotsymbol ({\thinlinefont .})
{\plot  6.032 21.018  5.524 20.320 /
}%
%
% Fig POLYLINE object
%
\linethickness= 0.500pt
\setplotsymbol ({\thinlinefont .})
{\plot  6.191 20.574  5.747 20.003 /
}%
%
% Fig POLYLINE object
%
\linethickness= 0.500pt
\setplotsymbol ({\thinlinefont .})
{\plot  5.302 23.178  4.985 22.479 /
}%
%
% Fig POLYLINE object
%
\linethickness= 0.500pt
\setplotsymbol ({\thinlinefont .})
{\plot  5.175 23.622  4.826 22.892 /
}%
%
% Fig POLYLINE object
%
\linethickness= 0.500pt
\setplotsymbol ({\thinlinefont .})
{\plot  4.255 24.289  5.143 23.685 /
\plot  5.143 23.685  6.287 20.225 /
}%
%
% Fig POLYLINE object
%
\linethickness= 0.500pt
\linethickness= 0.500pt
\setplotsymbol ({\thinlinefont .})
{\plot  4.381 23.559  4.000 22.924 /
}%
%
% Fig POLYLINE object
%
\linethickness= 0.500pt
\setplotsymbol ({\thinlinefont .})
{\plot  4.318 24.225  4.032 23.685 /
\plot  4.032 23.685  3.937 19.685 /
}%
%
% Fig POLYLINE object
%
\linethickness= 0.500pt
\setplotsymbol ({\thinlinefont .})
{\plot  4.508 23.050  3.969 22.225 /
}%
%
% Fig POLYLINE object
%
\linethickness= 0.500pt
\setplotsymbol ({\thinlinefont .})
{\plot  4.572 22.415  3.937 21.590 /
}%
%
% Fig POLYLINE object
%
\linethickness= 0.500pt
\setplotsymbol ({\thinlinefont .})
{\plot  4.667 21.876  3.937 21.050 /
}%
%
% Fig POLYLINE object
%
\linethickness= 0.500pt
\setplotsymbol ({\thinlinefont .})
{\plot  4.794 21.368  3.905 20.510 /
}%
%
% Fig POLYLINE object
%
\linethickness= 0.500pt
\setplotsymbol ({\thinlinefont .})
{\plot  4.826 20.701  3.905 20.034 /
}%
%
% Fig POLYLINE object
%
\linethickness= 0.500pt
\setplotsymbol ({\thinlinefont .})
{\plot  4.889 20.161  4.191 19.685 /
}%
%
% Fig POLYLINE object
%
\linethickness= 0.500pt
\setplotsymbol ({\thinlinefont .})
{\plot  3.238 20.161  2.826 19.939 /
}%
%
% Fig POLYLINE object
%
\linethickness= 0.500pt
\setplotsymbol ({\thinlinefont .})
{\plot  4.299 24.136  4.585 23.787 /
\plot  4.585 23.787  5.601 19.945 /
}%
%
% Fig POLYLINE object
%
\linethickness= 0.500pt
\setplotsymbol ({\thinlinefont .})
{\plot  4.259 24.293  4.386 23.658 /
\plot  4.386 23.658  4.925 19.689 /
}%
%
% Fig POLYLINE object
%
\linethickness= 0.500pt
\setplotsymbol ({\thinlinefont .})
{\plot  3.620 22.240  2.872 21.844 /
}%
% Fig POLYLINE object
%
\linethickness= 0.500pt
\setplotsymbol ({\thinlinefont .})
{\plot  3.539 21.685  2.699 21.304 /
}%
%
% Fig POLYLINE object
%
\linethickness= 0.500pt
\setplotsymbol ({\thinlinefont .})
{\plot  3.444 21.114  2.491 20.718 /
}%
%
% Fig POLYLINE object
%
\linethickness= 0.500pt
\setplotsymbol ({\thinlinefont .})
{\plot  3.317 20.559  2.318 20.225 /
}%
%
% Fig POLYLINE object
%
\linethickness= 0.500pt
\setplotsymbol ({\thinlinefont .})
{\plot  3.715 22.701  3.048 22.384 /
}%
%
% Fig POLYLINE object
%
\linethickness= 0.500pt
\setplotsymbol ({\thinlinefont .})
{\plot  4.314 24.253  3.520 23.713 /
\plot  3.520 23.713  2.282 20.189 /
}%
%
% Fig POLYLINE object
%
\linethickness= 0.500pt
\setplotsymbol ({\thinlinefont .})
{\plot  4.246 24.174  3.833 23.603 /
\plot  3.833 23.603  3.135 19.698 /
}%
% Fig POLYLINE object
%
\linethickness= 0.500pt
\setplotsymbol ({\thinlinefont .})
{\plot  5.508 22.574  5.080 21.971 /
}%
%
% Fig POLYLINE object
%
\linethickness= 0.500pt
\setplotsymbol ({\thinlinefont .})
{\plot  5.730 22.003  5.207 21.353 /
}%
%
% Fig POLYLINE object
%
\linethickness= 0.500pt
\setplotsymbol ({\thinlinefont .})
{\plot  5.952 20.083  6.174 19.448 /
%
% arrow head
%
\plot  6.030 19.667  6.174 19.448  6.150 19.709 /
}%
%
% Fig POLYLINE object
%
\linethickness= 0.500pt
\setplotsymbol ({\thinlinefont .})
{\plot  2.635 20.115  2.335 19.543 /
%
% arrow head
%
\plot  2.397 19.798  2.335 19.543  2.509 19.738 /
}%
%
% Fig POLYLINE object
% Fig POLYLINE object
%
\linethickness= 0.500pt
\setplotsymbol ({\thinlinefont .})
{\plot  4.445 19.702  4.462 19.035 /
%
% arrow head
%
\plot  4.392 19.287  4.462 19.035  4.519 19.291 /
}%
%
% Fig POLYLINE object
%
\linethickness= 0.500pt
\setplotsymbol ({\thinlinefont .})
\setdots < 0.0953cm>
{\plot  2.303 20.115  4.318 21.099 /
\plot  4.318 21.099  6.287 20.242 /
}%
%
% Fig POLYLINE object
%
\linethickness= 0.500pt
\setplotsymbol ({\thinlinefont .})
\setsolid
{\plot  2.000 24.098  2.004 24.100 /
\plot  2.004 24.100  2.013 24.107 /
\plot  2.013 24.107  2.030 24.117 /
\plot  2.030 24.117  2.051 24.132 /
\plot  2.051 24.132  2.081 24.149 /
\plot  2.081 24.149  2.115 24.168 /
\plot  2.115 24.168  2.153 24.189 /
\plot  2.153 24.189  2.191 24.208 /
\plot  2.191 24.208  2.231 24.227 /
\plot  2.231 24.227  2.273 24.244 /
\plot  2.273 24.244  2.318 24.263 /
\plot  2.318 24.263  2.366 24.280 /
\plot  2.366 24.280  2.419 24.297 /
\plot  2.419 24.297  2.477 24.314 /
\plot  2.477 24.314  2.540 24.331 /
\plot  2.540 24.331  2.580 24.342 /
\plot  2.580 24.342  2.623 24.352 /
\plot  2.623 24.352  2.667 24.363 /
\plot  2.667 24.363  2.714 24.376 /
\plot  2.714 24.376  2.762 24.386 /
\plot  2.762 24.386  2.815 24.399 /
\plot  2.815 24.399  2.870 24.412 /
\plot  2.870 24.412  2.929 24.424 /
\plot  2.929 24.424  2.991 24.437 /
\plot  2.991 24.437  3.054 24.450 /
\plot  3.054 24.450  3.120 24.460 /
\plot  3.120 24.460  3.188 24.473 /
\plot  3.188 24.473  3.258 24.486 /
\plot  3.258 24.486  3.330 24.498 /
\plot  3.330 24.498  3.401 24.509 /
\plot  3.401 24.509  3.476 24.519 /
\plot  3.476 24.519  3.550 24.530 /
\plot  3.550 24.530  3.624 24.541 /
\plot  3.624 24.541  3.698 24.549 /
\plot  3.698 24.549  3.772 24.558 /
\plot  3.772 24.558  3.846 24.566 /
\plot  3.846 24.566  3.918 24.572 /
\plot  3.918 24.572  3.990 24.579 /
\plot  3.990 24.579  4.062 24.585 /
\plot  4.062 24.585  4.132 24.589 /
\plot  4.132 24.589  4.202 24.591 /
\plot  4.202 24.591  4.271 24.594 /
\plot  4.271 24.594  4.339 24.596 /
\putrule from  4.339 24.596 to  4.407 24.596
\putrule from  4.407 24.596 to  4.477 24.596
\putrule from  4.477 24.596 to  4.547 24.596
\plot  4.547 24.596  4.616 24.594 /
\plot  4.616 24.594  4.686 24.589 /
\plot  4.686 24.589  4.758 24.585 /
\plot  4.758 24.585  4.830 24.581 /
\plot  4.830 24.581  4.902 24.577 /
\plot  4.902 24.577  4.976 24.570 /
\plot  4.976 24.570  5.048 24.564 /
\plot  5.048 24.564  5.122 24.555 /
\plot  5.122 24.555  5.194 24.547 /
\plot  5.194 24.547  5.266 24.539 /
\plot  5.266 24.539  5.336 24.530 /
\plot  5.336 24.530  5.406 24.519 /
\plot  5.406 24.519  5.474 24.511 /
\plot  5.474 24.511  5.541 24.500 /
\plot  5.541 24.500  5.605 24.490 /
\plot  5.605 24.490  5.666 24.479 /
\plot  5.666 24.479  5.726 24.469 /
\plot  5.726 24.469  5.783 24.460 /
\plot  5.783 24.460  5.836 24.450 /
\plot  5.836 24.450  5.886 24.439 /
\plot  5.886 24.439  5.935 24.431 /
\plot  5.935 24.431  5.980 24.420 /
\plot  5.980 24.420  6.024 24.412 /
\plot  6.024 24.412  6.064 24.403 /
\plot  6.064 24.403  6.102 24.395 /
\plot  6.102 24.395  6.164 24.380 /
\plot  6.164 24.380  6.219 24.367 /
\plot  6.219 24.367  6.270 24.354 /
\plot  6.270 24.354  6.316 24.342 /
\plot  6.316 24.342  6.358 24.329 /
\plot  6.358 24.329  6.394 24.316 /
\plot  6.394 24.316  6.428 24.304 /
\plot  6.428 24.304  6.456 24.293 /
\plot  6.456 24.293  6.477 24.280 /
\plot  6.477 24.280  6.496 24.270 /
\plot  6.496 24.270  6.511 24.259 /
\plot  6.511 24.259  6.519 24.249 /
\plot  6.519 24.249  6.526 24.238 /
\plot  6.526 24.238  6.530 24.229 /
\putrule from  6.530 24.229 to  6.530 24.219
\putrule from  6.530 24.219 to  6.530 24.208
\plot  6.530 24.208  6.528 24.200 /
\plot  6.528 24.200  6.524 24.189 /
\plot  6.524 24.189  6.515 24.177 /
\plot  6.515 24.177  6.505 24.164 /
\plot  6.505 24.164  6.492 24.151 /
\plot  6.492 24.151  6.473 24.136 /
\plot  6.473 24.136  6.452 24.122 /
\plot  6.452 24.122  6.424 24.105 /
\plot  6.424 24.105  6.394 24.088 /
\plot  6.394 24.088  6.358 24.071 /
\plot  6.358 24.071  6.318 24.052 /
\plot  6.318 24.052  6.272 24.035 /
\plot  6.272 24.035  6.223 24.016 /
\plot  6.223 24.016  6.170 23.997 /
\plot  6.170 23.997  6.111 23.975 /
\plot  6.111 23.975  6.049 23.954 /
\plot  6.049 23.954  6.009 23.944 /
\plot  6.009 23.944  5.969 23.931 /
\plot  5.969 23.931  5.925 23.918 /
\plot  5.925 23.918  5.878 23.906 /
\plot  5.878 23.906  5.829 23.891 /
\plot  5.829 23.891  5.779 23.878 /
\plot  5.779 23.878  5.723 23.863 /
\plot  5.723 23.863  5.666 23.848 /
\plot  5.666 23.848  5.607 23.834 /
\plot  5.607 23.834  5.544 23.821 /
\plot  5.544 23.821  5.480 23.806 /
\plot  5.480 23.806  5.412 23.791 /
\plot  5.412 23.791  5.345 23.777 /
\plot  5.345 23.777  5.275 23.764 /
\plot  5.275 23.764  5.203 23.749 /
\plot  5.203 23.749  5.131 23.736 /
\plot  5.131 23.736  5.059 23.724 /
\plot  5.059 23.724  4.985 23.713 /
\plot  4.985 23.713  4.913 23.700 /
\plot  4.913 23.700  4.839 23.692 /
\plot  4.839 23.692  4.767 23.681 /
\plot  4.767 23.681  4.695 23.673 /
\plot  4.695 23.673  4.623 23.666 /
\plot  4.623 23.666  4.553 23.660 /
\plot  4.553 23.660  4.483 23.654 /
\plot  4.483 23.654  4.413 23.650 /
\plot  4.413 23.650  4.343 23.645 /
\plot  4.343 23.645  4.276 23.643 /
\plot  4.276 23.643  4.208 23.641 /
\putrule from  4.208 23.641 to  4.138 23.641
\putrule from  4.138 23.641 to  4.068 23.641
\plot  4.068 23.641  3.998 23.643 /
\plot  3.998 23.643  3.926 23.645 /
\plot  3.926 23.645  3.854 23.647 /
\plot  3.854 23.647  3.782 23.652 /
\plot  3.782 23.652  3.708 23.658 /
\plot  3.708 23.658  3.634 23.662 /
\plot  3.634 23.662  3.560 23.671 /
\plot  3.560 23.671  3.486 23.677 /
\plot  3.486 23.677  3.412 23.685 /
\plot  3.412 23.685  3.340 23.694 /
\plot  3.340 23.694  3.268 23.705 /
\plot  3.268 23.705  3.198 23.715 /
\plot  3.198 23.715  3.128 23.726 /
\plot  3.128 23.726  3.061 23.736 /
\plot  3.061 23.736  2.997 23.747 /
\plot  2.997 23.747  2.934 23.760 /
\plot  2.934 23.760  2.874 23.772 /
\plot  2.874 23.772  2.817 23.783 /
\plot  2.817 23.783  2.762 23.796 /
\plot  2.762 23.796  2.711 23.806 /
\plot  2.711 23.806  2.663 23.819 /
\plot  2.663 23.819  2.618 23.832 /
\plot  2.618 23.832  2.576 23.842 /
\plot  2.576 23.842  2.536 23.855 /
\plot  2.536 23.855  2.498 23.865 /
\plot  2.498 23.865  2.436 23.884 /
\plot  2.436 23.884  2.383 23.906 /
\plot  2.383 23.906  2.335 23.925 /
\plot  2.335 23.925  2.290 23.946 /
\plot  2.290 23.946  2.250 23.965 /
\plot  2.250 23.965  2.216 23.986 /
\plot  2.216 23.986  2.187 24.005 /
\plot  2.187 24.005  2.161 24.024 /
\plot  2.161 24.024  2.140 24.043 /
\plot  2.140 24.043  2.123 24.060 /
\plot  2.123 24.060  2.110 24.075 /
\plot  2.110 24.075  2.102 24.090 /
\plot  2.102 24.090  2.095 24.102 /
\plot  2.095 24.102  2.089 24.113 /
\plot  2.089 24.113  2.087 24.122 /
\plot  2.087 24.122  2.085 24.130 /
\putrule from  2.085 24.130 to  2.085 24.132
\plot  2.085 24.132  2.083 24.136 /
\putrule from  2.083 24.136 to  2.083 24.138
\putrule from  2.083 24.138 to  2.083 24.141
\putrule from  2.083 24.141 to  2.083 24.143
\putrule from  2.083 24.143 to  2.083 24.145
\putrule from  2.083 24.145 to  2.081 24.145
\putrule from  2.081 24.145 to  2.081 24.147
\putrule from  2.081 24.147 to  2.081 24.149
\putrule from  2.081 24.149 to  2.083 24.149
\putrule from  2.083 24.149 to  2.085 24.149
\putrule from  2.085 24.149 to  2.085 24.147
\plot  2.085 24.147  2.087 24.145 /
\putrule from  2.087 24.145 to  2.087 24.143
\putrule from  2.087 24.143 to  2.089 24.143
\putrule from  2.089 24.143 to  2.089 24.141
\plot  2.089 24.141  2.091 24.138 /
\putrule from  2.091 24.138 to  2.091 24.136
\plot  2.091 24.136  2.093 24.134 /
\putrule from  2.093 24.134 to  2.093 24.132
\putrule from  2.093 24.132 to  2.095 24.132
\putrule from  2.095 24.132 to  2.095 24.130
}%
%
%
% Fig POLYLINE object
%
\linethickness= 0.500pt
\setplotsymbol ({\thinlinefont .})
{\plot  3.747 23.114  3.238 22.909 /
}%
%
% Fig POLYLINE object
%
\linethickness= 0.500pt
\setplotsymbol ({\thinlinefont .})
{\plot  3.793 23.527  3.349 23.305 /
}%
%
% Fig TEXT object
%
\put{$O(1)$} [lB] at  7.525 22.860
%
% Fig TEXT object
%
%
% Fig TEXT object
%
\put{$d$} [lB] at  4.858 23.971
%
% Fig TEXT object
%
\put{$\text{supp}\ \hat f_3$} [lB] at  6.208 18.860
%
% Fig TEXT object
%
\put{$\text{supp}\ \hat f_1$} [lB] at  1.600 18.876
%
% Fig TEXT object
%
\put{$\text{supp}\ \hat f_2$} [lB] at  4.000 18.447
\linethickness=0pt
\putrectangle corners at  1.109 24.621 and  7.557 18.415
\endpicture}
$$

\smallskip

We obtain then
\smallskip

\noindent
{\bf Proposition 2$'$.}
Under the above assumptions on $f_1, f_2, f_3$, one has the inequality
$$
\Vert\, |f_1 f_2
f_3|^{\frac
13}\Vert_{L^3_{\#} (B_N)}\ll
\delta^{-\frac 12}
N^{0+} \prod^3_{i=1} \Vert
f_i\Vert^{1/3}
_{L^2_{\#(B_N)}}.\eqno{(2.3)}
$$
\smallskip

We will not explain the deduction of Proposition 2$'$ from Proposition 2  in detail, but just point
out that it is based on the
rescaling map
$$
L^\sigma: C\to C: (\xi_1, \xi_2, \xi_3)\mapsto \Big(\frac {\xi_1-\xi_3}{2\sigma} + \frac
{\xi_1+\xi_3}2, \frac {\xi_2}{
\sqrt \sigma}, \frac {\xi_3-\xi_1}{2\sigma}+\frac {\xi_1+\xi_3}{2}\Big).\eqno{(2.4)}
$$
mapping $C\cap R_\delta$ to $C\cap R_{\frac \delta{\sqrt\sigma}}$.

Note that obviously $\Vert
f_i\Vert_{L^2_{\#}(B_N)}\sim \big(\sum_\sigma\Vert f_{i, \sigma}
\Vert^2_{L^2_{\#}(B_N)}\big)^{\frac 12}$
by orthogonality.
\smallskip

Next, we perform an interpolation between (2.1) and (2.3), setting
$$
\cases
\frac 14 =\frac {1-\theta}2+\frac \theta 6\\
\frac 1q =\frac {1-\theta}3+\frac \theta 6\endcases
$$
with $\theta=\frac 34, q=\frac {25}4$.

We obtain

\noindent
{\bf Proposition 3.}
{\sl Let $f_1, f_2, f_3$ be as in Proposition 2$'$. Then
$$
\Vert\, |
f_1f_2f_3|^{\frac
13}\Vert_{L^q_{\#} (B_N)}\ll
\delta^{-\frac 12(1-\theta)}
N^{o+} \prod^3_{i=1} \Big(\sum_\sigma\Vert
f_{i,
\sigma}\Vert^2_{L^4_{\#}(B_N)}\Big)^{\frac 16}.
$$}

Note that this interpolation is not trivial and requires the `balanced wave packet decompositions'
introduced in \cite {B-D}.

Our final step is to derive from the multi-linear inequalities (2.5) a linear inequality.
This is a relatively easy multiscale argument going back to \cite{B-G} and which will not be
repeated here.
The upshot of this argument is that we recover $\Vert
f\Vert_{L^q_{\#}(B_N)} $ from the contributions
$$
\eqalignno
{&\delta^{-\frac 18}
N^{o+}\Big(\sum_\alpha \Big(\sum_{\sigma\subset\alpha} \Vert
f_\sigma\Vert^2_{L^4_{\#}(B_N)}\Big)^{\frac
q2}\Big)^{\frac 1q}\leq\cr
&N^{o+} \delta^{-\frac 18} (\delta\sqrt
N)^{\frac 14} \Big(\sum_\alpha\Big(\sum_{\sigma\subset\alpha} \Vert
f_\sigma\Vert^4_{L^4_{\#}(B_N)}\Big)^{\frac
q4}\Big)^{\frac 1q}\leq\cr
&N^{\frac 18+} \Big(\sum_\sigma\Vert
f_\sigma\Vert_{L^4_{\#}(B_N)}\Big)^{\frac 14}
}
$$
where $\{\alpha\}$ refers to a partition of $C_{\frac 1N}$ in $\delta$-slabs with $\delta$ taking
dyadic values between $\frac
1{\sqrt N}$ and 1.

\smallskip
The final statement in this section is following

\noindent
{\bf Proposition 4.}
{\sl If supp\,${\hat f\subset C_{\frac 1N}}$, then for $q=\frac {25}4$
$$
\Vert f\Vert_{L^q_{\#}(B_N)} \ll
N^{\frac 18+} \Big(\sum\Vert
f_\sigma\Vert^4_{L^4_{{\#}(B_N)}}\Big)^{\frac 14}.\eqno{(2.6)}
$$}

\bigskip

\beginsection
3.  Decoupling for the perturbed cone

We are  now returning to the perturbed cone $\tilde C$ defined by (1.3).

\smallskip
Note that at scale $N\sim \frac
1{\ve^2}$, $\tilde C$ may be identified with the cone
$$
C:\cases z=t\\ y=st\\
x=s^{\frac 12}t\endcases
$$
and hence at this scale (2.6) remains applicable to $\tilde C$ as well.

Assuming $\eta
<\ve^2$, it follows that if supp\,$\hat f\subset \tilde C_{\eta}$, then
$$
\Vert f\Vert _{L^4_{\#}(B_{\frac 1\eta})
L^q_{\#}(B_{x, \frac
1{\ve^2}})}\ll \
\Big(\frac
1\ve\Big)^{\frac 14+} \Big(\sum_\sigma
\Vert
f_\sigma\Vert^4_{L^4_{\#}(B_{\frac
1\eta})}\Big)^{\frac 14}.\eqno{(3.1)}
$$
Here the left side stands for
$$
\Big(\frac 1{|B_{\frac 1\eta}|} \int_{B_{\frac 1\eta}}\Vert
f\Vert^4_{L^q_{\#}(B_{x, \frac 1
{\ve^2}})}
dx\Big)^{\frac 14}
$$
and (3.1) is deduced from (2.6) just by partitioning $B_{\frac 1\eta}$ in balls of size $N=\frac
1{\ve^2}$.
The slabs $\sigma$ have angular width $\ve$ .

Exploiting the perturbative term
$\ve^2 s^{-\frac 32}
t^2+O(\ve^4)$ in (1.3), we will perform a further decoupling of $f_\sigma$.

Let us first rewrite (1.3) as
$$
x=y^{\frac 12}z^{\frac 12}
-c\ve^2 \, \frac
{z^{7/2}}{y^{3/2}}+ O
(\ve^4)\eqno{(3.2)}
$$
and making the substitution
$$
\cases
z=z_1+y_1\\ y=z_1-y_1\sim 1
\endcases
$$
$$
x^2=z^2_1 -y_1^2
-2c\ve^2 \ \frac
{(z_1+y_1)^4}{z_1-y_1}+O(\ve^4).\eqno{(3.3)}
$$
Fixing $\sigma$, perform a rotation in $(x_1, y_1)$-plane
$$
\cases
x=(\cos \theta) x_2 -(\sin \theta) y_2\\ y_1= (\sin \theta) x_2+(\cos \theta) y_2
\endcases
$$
to put $\sigma$ in position $|y_2|<\ve$.
Writing $y_2=\ve y_3$, (3.3) becomes then
$$
\eqalignno
{z_1^2 &=x_2^2
+\ve^2 y^2_3+
2c\ve^2 \, \frac {(z_1 +(\sin\theta) x_2+\ve(\cos\theta)
y_3)^4}{z_1-(\sin\theta) x_2 -\ve(\cos
\theta) y_3}
+O(\ve^4)\cr
&=x_2^2 +\ve^2y^2_3 +
2c\ve^2 \ \frac {(z_1+(\sin \theta)
x_2)^4}{z_1-(\sin\theta) x_2}+
O(\ve^3).}
$$
Hence
$$
z_1=x_2+\frac
12\ve^2 \, \frac
{y^2_3}{x_2}
+c\ve^2 \, \frac
{(1+\sin\theta)^4}{1-\sin\theta}
x_2^2 +O(\ve^3).\eqno{(3.4)}
$$
Set $z_1-x_2 =
\ve^2 z_2$ to obtain
$$
z_2 =\frac 12 \,
\frac{y^2_3}{x_2} +C \, \frac
{(1+\sin\theta)^4}{1-\sin\theta}
x_2^2+O(\ve).\eqno{(3.5)}
$$
Note that the Hessian equals
$$
\left|\matrix \frac
{y_3^2}{x_2^3} + 2c \, \frac
{(1+\sin\theta)^4}{1-\sin \theta}&
-\frac{y_3}{x_2^2}\\
-\frac {y_3}{x_2^2} & \frac 1{x_2}\endmatrix\right| = 2c\, \frac
{(1+\sin\theta)^4}{1-\sin \theta} \,  \frac 1{x_2}.
\eqno{(3.6)}
$$
Since $z_1-y_1\sim 1, z_1 -x_2 \sin \theta \sim 1$ and, by (3.4), $(1-\sin\theta)x_2\sim 1$.
Thus there is a further decoupling in $(x_2, y_3)$ at scale $\frac{\sqrt\eta}\ve$, hence in $(x_2,
y_2)$ at scale $\big(\frac {\sqrt         \eta}\ve,
\sqrt\eta)\Rightarrow \sqrt\eta$-angular, $\frac{\sqrt \eta}\ve$-radial in $(x_1 y_1)$-space.
Since $t=z$, there is a decomposition in $t$ at scale $\frac{\sqrt\eta}{\ve}$.
Since $tg  \theta =\frac {y_1}x =\frac {z-y}{2x} =\frac
{1-s}{2(s^{\frac 12}
-c\ve^2 s^{-\frac 32} t)}$, the angular decomposition at scale
$\sqrt \eta$ corresponds to a decomposition in $s$ at scale $\sqrt\eta$.

Returning to (3.1), the preceding leads to the further decoupling
$$
\Vert f\Vert_{L^4_{\#}(B_{\frac
1\eta})L^q_{\#}(B_{x, \frac
1{\ve^2}})}\ll \Big(\frac
1\ve\Big)^{\frac 14+} \Big(\frac
\ve{\sqrt\eta}\Big)^{\frac 12}
\Big(\sum_\tau \Vert f_\tau
\Vert^4_{L^4_{\#}(B_{\frac 1\eta})}
\Big)^{\frac 14}\eqno{(3.7)}
$$
with $\{\tau\}$ a partition of $\tilde C$ at scale $\sqrt\eta$ in $s$ and scale $\frac {\sqrt
\eta}{\ve}$ in $t$.
Hence, with $q=\frac {24}5$ and assuming $\ve>\sqrt\eta$, we get

\noindent
{\bf Proposition 5.}
{\sl The following decoupling inequality holds for $\tilde C$.
Let supp\,$\hat f\subset\tilde C_\eta$. Then
$$
\Vert f\Vert_{L^4_{\#}(B_{\frac 1\eta})
L^q_{\#} (B_{x, \frac
1{\ve^2}})} \ll
\ve^{\frac
14-}\eta^{-\frac 14}
\Big(\sum_\tau \Vert
f_\tau\Vert^4_{L^4_{\#}(B_{\frac 1 \eta})}
\Big)^{\frac 14}\eqno{(3.8)}
$$
with $\{\tau \}$ a partition in $(\sqrt\eta, \frac{\sqrt\eta}\ve)$ rectangles in $(s, t)$.
}

Recall that $s=\frac kK, t=\frac\ell L, \ve=\frac LK, \omega(k,
\ell)=k^{\frac 12}\ell
+ck^{-\frac
32}\ell^3+\ldots
$ and we assumed $\eta < \frac
{L^2}{K^2}$.
In future applications, $\eta$ will moreover satisfy
$$
\frac 1{KL}<\eta< \frac 1K< \frac 1L.\eqno{(3.9)}
$$
In particular,
$K<L^3$.

Obviously the ball $B_{\frac 1\eta}$ in (3.8) may be replaced by any larger domain of the form \hfill\break
$[|x|<X_1]\times [|x_2|<X_2] \times [|x_3|< X_3]$ where $X_1, X_2, X_3\geq \frac 1\eta$.

\smallskip

Returning to (1.2), it follows that
$$
\eqalignno
{&\Big\Vert\sum_{k\sim K, \ell \sim L} a_{k\ell} \, e\Big(\frac \ell Lx_1 +\frac {k\ell}{KL} x_2+\frac {\omega(k, \ell)}{\sqrt KL}
x_3\Big)\Big\Vert_{L^4_{\#}[x_1|<\frac 1\eta, |x_2|< KL, |x_3|<\frac 1\eta] L^q_{\#}(B_{x, \frac {K^2}{L^2}})}\cr
&\ll \Big(\frac LK\Big)^{\frac 14-} \eta^{-\frac 14} \Big(\sum_{\alpha, \beta}\Big\Vert\sum_{k\in I_\alpha, \ell\in J_\beta}
a_{k\ell} \, e\Big(\frac \ell L x_1+\frac {k\ell}{KL} x_2+\frac {\omega(k, \ell)}{\sqrt kL}x_3\Big)\Big\Vert
^4_{L^4_{\#} [|x_1|<\frac 1\eta, |x_2|<KL, |x_3|
<\frac 1\eta]}\Big)^{\frac 14}}
$$
with $\{I_\alpha\}$ a partition of $[k\sim K]$ in $\sqrt\eta K$-intervals and $\{J_\beta\}$ a partition of $[\ell\sim L]$ in $\frac
{\sqrt   \eta}\ve
L=\sqrt\eta K$-intervals.

\smallskip

Note that the function $\sum_{k\sim K, \ell \sim L} a_{k\ell} \, e\big(\ell x_1+k\ell x_2+ \omega(k, \ell)x_3\big) $ in (1.2) is
1-periodic  in $x_1, x_2$.
Thus the previous inequality may be reformulated as
\medskip

\noindent
{\bf Proposition 6.}
$$
\eqalignno
{&\left\Vert \sum_{k\sim K, \ell \sim L} a_{k\ell} \, e(\ell x_1+k\ell x_2+\omega(k, \ell)x_3)\right\Vert
{\Sb {L^4_ {\#} [|x_1|< K, |x_2|< 1, |x_3|<\frac 1{\sqrt KL\eta}]}\\
{L^q_{\#} ([|x_1|<\frac {K^2}{L^3}, |x_2|<\frac K{L^3}, |x_3|<\frac {K^{3/2}}{L^3}]+x)}\endSb} \cr
&\ll\left(\frac LK\right)^{\frac 14-} \eta^{-\frac 14} \left(\sum_{\alpha, \beta} \left\Vert\sum_{k\in I_\alpha, \ell
\in J_\beta} e\left(\ell x_1+k\ell x_2+\omega(k, \ell)x_3\right)
\right\Vert^4_{L^4_{\#}[|x_1|<1, |x_2|<1, |x_3|< \frac 1{\sqrt K L\eta}]}\right)^{\frac 14}.
}$$
\smallskip

Clearly the expression
$$
\sum_{\alpha, \beta}\left\Vert\sum_{k\in I_\alpha, \ell \in J_\beta} e(\ell x_1+k\ell x_2+ \omega(k,
\ell)x_3)\right\Vert^4_{L^4_{\#}[|x_1|<1,
|x_2|<1,|x_3|<\frac 1{\sqrt K 2\eta}]}
$$
amounts to the number of integral solutions of the system
$$
\qquad\qquad \
\cases
\ell_1+\ell_2=\ell_3+\ell_4 \qquad &{(3.10)}\cr
\ell_1k_1+\ell_2k_2= \ell_3k_3+\ell_4k_4\qquad &{(3.11)}\cr
\omega(k_1, \ell_1)+\cdots -\omega(k_4, \ell_4)=O(\eta \sqrt K L)\qquad\qquad\qquad\qquad\qquad\qquad\qquad &{(3.12)}\cr
\endcases
$$
with
$$
\alignat2
&k\sim K, \ell \sim L &\cr
&\text{diam} \, (k_1, k_2, k_3, k_4)\leq \sqrt\eta K \qquad  \tag{3.13}\cr
&\text{diam}\, (\ell_1, \ell_2, \ell_3, \ell_4)< \frac {\sqrt \eta} \ve L = \sqrt\eta K\qquad\qquad\qquad\qquad\qquad
\qquad \tag{3.14}
\endalignat
$$
\medskip

\noindent
{\bf Proposition 7.}
{\sl The number of solutions of (3.10)-(3.14) is bounded by
$$
\eta^2 K^{5+} +\eta K^3L.\eqno{(3.15)}
$$}
\smallskip

\noindent
{\bf Proof.}
We discard (3.12) which in fact is easily seen to be redundant.
In what follows, we ignore the effect of divisor functions, which introduce an extra factor $K^{O+}$.

Set $\Delta k_i=k_i-k_4 (i=1, 2, 3)$. Thus $|\Delta k_i|\lesssim \sqrt\eta K$.

Since
$$
\ell_1\Delta k_1+\ell_2\Delta k_2=\ell_3 \Delta k_3 \eqno{(3.16)}
$$
$$
(\ell_1-\ell_3) \Delta k_1+(\ell_2-\ell_3)\Delta k_2 =\ell_3 (\Delta k_3-\Delta k_1 -\Delta k_2).\eqno{(3.17)}
$$
Assume $\Delta k_1+ \Delta k_2 \not= \Delta k_3$.
Choose $\ell_1-\ell_3, \ell_2-\ell_3, \Delta k_1, \Delta k_2$ ($\eta K^2  \frac \eta {\ve 2} L^2$-possibilities)
$\overset {(3.17)}\to \Rightarrow \ell_3, \Delta k_3$.

 Since there are $K$ possibilities for $k_4$, this gives
$$
\frac {\eta^2}{\ve^2} K^3 L^2 =\eta^2 K^5.
$$
If $\Delta k_1+\Delta k_2 =\Delta k_3, (\ell_1-\ell_3)\Delta k_1+ (\ell_2-\ell_3) \Delta k_2=0$.

If $(\ell_1-\ell_3) \Delta k_1\not= 0$, choose $\ell_1, \ell_1-\ell_3, \Delta k_1$ \, ($\frac {\sqrt\eta}\ve L^2 \sqrt\eta
K$-possibilities)
$$
\Rightarrow \ell_2, \Delta k_2, \Delta k_3
$$
which gives the contribution $\eta LK^3$.

Case $(\ell_1-\ell_3)\Delta k_1= (\ell_2-\ell_3) \Delta k_2=0$.
$$
\alignat2
&\bullet \ell_1=\ell_2=\ell_3=\ell_4 \Rightarrow \quad\qquad &L\eta K^2.K=\eta LK^3&\cr
&\bullet \ell_1=\ell_3, \ell_2=\ell_4, \Delta k_2=0\Rightarrow \quad \qquad
&L\frac {\sqrt \eta}\ve L\sqrt\eta KK=\eta K^3L \text { by
\, (3.11)}&\cr
&\bullet \Delta k_1 =\Delta k_2=\Delta k_3 =0\Rightarrow  k_1 =k_2=k_3 = k_4\Rightarrow
 & K.L. \Big(\frac {L\sqrt\eta}\ve\Big)^2 = \eta K^3L&.
\endalignat
$$
This proves Proposition 7.
$\quad\blacksquare$
\medskip

The final statement of this section becomes then

\medskip
\noindent
{\bf Proposition 8.}
{\sl Let $|a_{k, \ell}|\leq 1$ (arbitrary) and $q=\frac {24}5$.
Then
$$
\Big\Vert\sum_{k\sim K, \ell \sim L}  a_{k\ell} \, e(\ell x_1+ k\ell x_2+\omega (k, \ell)x_3)\Big\Vert{\Sb L^4_{\#}
[|x_1|<K, |x_2| < 1, |x_3|< \frac 1{\sqrt KL\eta}]\\
L^q_{\#} ([|x_1|<\frac {K^2}{L^3}, |x_2|<\frac K{L^3}, |x_3|<\frac {K^{3/2}}{L^3}]+x)\endSb}
$$
$$
\ll K^{\frac 12+}\, L^{\frac 12} \Big(1+\eta\frac {K^2}{L}\Big)^{\frac 14}
\eqno{(3.18)}
$$
under assumption (3.9) on $\eta$.}
\bigskip

\beginsection
{4.  The basic moment inequalities}

The main results from this section are Propositions 10 and 10$'$.

In order to get an estimate on (1.2) for some $q_\nu >4$, we interpolate (3.18) with a bound on
$$
\Big\Vert\sum_{k\sim K, \ell \sim L} a_{k\ell} \ e(\ell x_1+k\ell x_2+\omega (k, \ell)x_3)\Big\Vert\Sb
L^{2\nu}_{\#}[|x_1|<K, |x_2|<1, |x_3|< \frac 1{\sqrt KL\eta}]\\
L^2_{\#}([|x_1|<\frac {K^2}{L^3}, |x_2|<\frac K{L^3}, |x_3|<\frac {K^{3/2}}{L^3}]+x)\endSb\eqno{(4.1)}
$$
with $\nu\in\Bbb Z, \nu\geq 3$.

Note that if $\Delta\ell <\frac {L^3}{K^2}, \Delta(k\ell)< \frac {L^3}K$, then $\Delta k<\frac {L^2}K$.

Denoting $F(x_1, x_2, x_3)= \sum_{k\sim K, \ell \sim L} a_{k\ell} \, e\big(\ell x_1+k\ell x_2+\omega(k,
\ell)x_3\big)$, it follows that (4.1) is bounded by
$$
\Big\Vert\Big(\sum|F_\tau|^2\Big)^{\frac 12}\Big\Vert_{L^{2\nu}_{\#} [|x_1|<1, |x_2|<1, |x_3|<\frac 1{\sqrt K
L\eta}]}\eqno{(4.2)}
$$
with $\{\tau\}$ a partition of $[k\sim K]\times [\ell\sim L]$ in intervals $I_\alpha\times J_\beta$ with $|I_\alpha|< \frac
{L^2}K, |J_\beta|< \frac {L^3}{K^2}$.

We used here again periodicity of $F$ in $x_1$.
If necessary, we refine the partition further as to restrict $\ell$ to intervals of size 1. Thus
$$
(4.2) <\Big(1+\frac {L^3}{K^2}\Big)^{\frac 12-\frac 1{2\nu}} \Big\Vert\Big(\sum|F_{\tau'}|^2\Big)^{\frac
12}\Big\Vert_{L^{2\nu}_{\#} [|x_1|<1, |x_2|<1, |x_3|<\frac 1{\sqrt K L\eta}]}\eqno{(4.3)}
$$
with $\{\tau'\}$ a partition in intervals $I_\alpha \times J_{\beta}$, $|I_\alpha| =\frac {L^2}K, |J_\beta'|<1$.

Evaluation of
$$
\Big\Vert\Big(\sum|F_{\tau'}|^2\Big)^{\frac 12}\Big\Vert^{2\nu}_{L^{2\nu}_{\#} [|x_1|<1, |x_2|< 1, |x_3|< \frac 1{\sqrt K
L\eta}]}
$$
amounts to the number of integral solutions of
$$
\, \ \qquad\qquad\cases
(k_1-k_2) \ell_1+(k_3-k_4)\ell_3+\cdots+ (k_{2\nu-1}-k_{2\nu}) \ell_{2\nu-1} =0\qquad&{(4.4)}\cr
\omega(k_1, \ell_1)-\omega(k_2, \ell_1) +\cdots+\omega(k_{2\nu-1}, \ell_{2\nu-1}) -\omega (k_{2\nu}, \ell_{2\nu-1})=O (\sqrt
KL\eta)\qquad \quad \
 &{(4.5)}
\endcases
$$
in $k_1, \ldots, k_{2\nu}\sim K, \ell_1, \ell_3, \ldots, \ell_{2\nu-1} \sim L$ and with
$$
|k_1-k_2|, \ldots, |k_{2\nu-1} - k_{2\nu}|<\frac {L^2}K.\eqno{(4.6)}
$$
Equivalently, consider the system
$$
\ \, \quad\qquad
\cases
u_1\ell_1+u_2\ell_2+\cdots+ u_\nu \ell_\nu =0\qquad \ & {(4.7)}\cr
\omega(k_1+u_1, \ell_1)-\omega (k_1, \ell_1)+\cdots+\omega (k_\nu+u_\nu, \ell_\nu)-\omega(k_\nu, \ell_\nu)< O(\sqrt
KL\eta) \qquad\qquad &{(4.8)}
\endcases
$$
with $k_i \sim K, \ell_i\sim L$ and $u_i =O(\frac {L^2}{K})$.
\smallskip

Assume $|u_1|\geq |u_2|, \ldots, |u_\nu|$.
The contribution of $u_1=0$ is $L^\nu K^\nu$.
Next, consider the contribution of $|u_1|\sim U\not= 0$,
$U\lesssim L^2K^{-1}>1$.
Fix $k_2, \ldots, k_\nu, \ell_2, \ldots, \ell_\nu, u_2, \ldots, u_\nu$.
From (4.7) we retrieve $\ell_1, u_1$.
Considering (4.8) as an equation in $k_1$, we obtain the bound
$$
\align
&K^{\nu-1} L^{\nu-1} U^{\nu-1} \Big( 1+\frac {\eta K^2}{U}\Big)\leq\cr
&K^{\nu-1} L^{\nu-1} \Big(\frac {L^2}K\Big)^{\nu-1}+\eta K^{\nu+1} L^{\nu-1} \Big(\frac {L^2}K\Big)^{\nu-2}\leq\tag {4.9}\cr
&L^{3\nu-3}+\eta K^3 L^{3\nu-5}.
\endalign
$$
Thus the number of solutions of (8.8)-(8.10) is at most
$$
K^{\nu+\ve}L^\nu \Big(1+\frac {L^{2\nu-3}}{K^\nu}+(\eta KL)\frac {L^{2\nu-6}}{K^{\nu-2}}\Big).\eqno{(4.10)}
$$
Hence
\medskip

\noindent
{\bf Proposition 9.}
{\sl $$
(4.1) \ll \Big(1+\frac {L^3}{K^2}\Big)^{\frac 12-\frac 1{2\nu}} \Big( 1+\frac {L^{2\nu-3}}{K^\nu}+ (\eta KL)\frac
{L^{2\nu-6}}{K^{\nu-2}} \Big)^{\frac 1{2\nu}} K^{\frac 12+}L^{\frac 12}.\eqno{(4.11)}
$$}
\medskip
It remains to interpolate between (3.18) and (4.11).

\smallskip

We obtain the following
\medskip

\noindent
{\bf Proposition 10.}
{\sl For $\nu\geq 3$, take $q_\nu =\frac {13\nu -12}{3\nu -\frac 52}> 4$.

We have, assuming $\eta <\frac {L^2}{K^2} $ and $\frac
1{KL} <\eta <\frac 1K<\frac 1L$
$$
\aligned
&\Big\Vert\sum_{k\sim K, \ell \sim L} a_{k\ell} \, e(\ell x_1+k\ell x_3+\omega (k, \ell)x_3)\Big\Vert_{L^{q_\nu} [|x_1|<1,
|x_2|<1, |x_3|<\frac 1{\sqrt KL\eta}]}\cr
&\ll \Big(1+\eta \, \frac {K^2}L\Big)^{\frac {3(\nu-1)}{13\nu-12}} \Big(1+\frac {L^3}{K^2}\Big)^{\frac {\nu-1}{2(13\nu-12)}}
\Big(1+\frac {L^{2\nu-3}}{K^\nu} +(\eta KL) \frac {L^{2\nu-6}}{K^{\nu-2}}\Big)
^{\frac 1{2(13\nu-12)}} K^{\frac 12+}L^{\frac 12}.
\endaligned
\eqno{(4.12)}
$$}
\medskip

Note that if $\eta>\frac {L^2}{K^2}$, we may ignore the $\ve$-terms in (1.3) i.e. we are in the pure conical situation.
We get the inequality
$$
\Vert f\Vert_{L_{\#}^{\frac {24}5}(B_{\frac 1\eta})} \ll \Big(\frac 1\eta\Big)^{\frac 18+\ve} \Big(\sum\Vert
f_\sigma\Vert^4_{L^4_{\#}(B_{\frac 1\eta})} \Big)^{\frac 14}\eqno{(4.13)}
$$
instead of (2.6), with $\{\sigma\}$ a partition in $\sqrt\eta$-plates.

Hence instead of Proposition 6, we obtain, with $F$ defined as above,

$$
\aligned
&\Vert F\Vert_{L^4_{\#}[|x_1|<K, |x_2|< 1, |x_3|<\frac 1{\sqrt KL\eta}]}
L^{\frac {24}5}_{\#} \Big(x+\Big[|x_1|<\frac 1{L\eta}, |x_2| <\frac 1{KL\eta}, |x_3|< \frac 1{\sqrt KL\eta}\Big]\Big)\cr
&\ll \Big(\frac 1\eta\Big)^{\frac 18+} \Big(\sum_\alpha \Vert F_{I_\alpha}\Vert^4_{L^4_{\#}[|x_1|<1, |x_2|<1, |x_3|<\frac
1{\sqrt KL\eta}]}\Big)^{\frac 14}.
\endaligned
\eqno{(4.14)}
$$
with $\{I_\alpha\}$ a partition of $[k\sim K]$ in $\sqrt\eta K$-intervals.
\smallskip

The expression $\sum_\alpha \Vert F_{I_\alpha}\Vert^4_{L^4_{\#}}$ amounts to the number of integral solutions of
(3.10)-(3.12) under the only restriction (3.13) and from the analysis in Proposition 7, we get the bound
$$
\Big(\frac 1\eta\Big)^{\frac 18} K^\ve \big(L^2(\sqrt\eta K)^2 K +\sqrt\eta K^2 L^2 +KL^3\Big)^{\frac 14}.\eqno{(4.15)}
$$
Hence, since $\eta \gtrsim L^2/K^2$
$$
\aligned
(4.14)&\ll K^\ve (\sqrt\eta K^3 L^2+K^2L^2+\eta^{-\frac 12}KL^3)^{\frac 14}\cr
&\ll K^{\frac 12+} L^{\frac 12}\Big(1+\eta \frac {K^2}L\Big)^{\frac 14}
\endaligned
\eqno{(4.16)}
$$
which is the same as the r.h.s. of (3.18).

\smallskip

In the l.h.s. of (3.18), $\frac {K^2}{L^2}$ is replaced by $\frac 1\eta$ and, using only the $L^2$-norm, there is by (3.9)
the trivial bound $K^{\frac 12}$ $L^{\frac 12}$ on (4.1).

It follows that Proposition 10 remains valid without the assumption $\eta <\frac {L^2}{K^2}$.
Hence
\smallskip

\noindent
{\bf Theorem 1.}
{\sl For $\nu\geq 3$, $q_\nu =\frac {13\nu-12}{3\nu-\frac 52}$ and $\eta$ satisfying
$ \frac 1{KL}<\eta <\frac 1K<\frac 1L$, inequality (4.12) holds.
}
\bigskip

Note that for $\eta =\frac 1{KL}$, the first factor in (4.12) becomes $1+\frac K{L^2}$.
For $k>L^2$, we establish an alternative bound.

Assume
$$
K\geq L^2.\eqno{(4.17)}
$$

We may then replace Proposition 10 by

\noindent
{\bf Proposition 10'.}
{\sl
$$
\aligned
&\Big\Vert\sum_{k\sim K, \ell\sim L} a_{k\ell} \ e\big(\ell x_1+ k_1\ell x_2+ \omega (k, \ell)x_3\big)
\Big\Vert _{L^9_{\#} [|x_1|< 1, |x_2|<1, |x_3|<\frac 1{\sqrt K L\eta}]}\cr
&\ll (1+\eta KL)^{\frac 1{48}} \Big(1+\eta\frac {K^2}L\Big)^{\frac 5{24}} K^{\frac 12}+L^{\frac 12}
\endaligned
\eqno{(4.18)}
$$
where $q=\frac {48}{11}$.}

\medskip

\noindent
{\bf Sketch of the Argument}

Instead of (4.1), we will bound
$$
\Big\Vert\sum_{k\sim K, \ell\sim L}\cdots \Big\Vert _{L^{2\nu}_{\#} l^3_{\#}}\eqno{(4.19)}
$$
appealing to inequalities (2.1), (2.3) derived from the multi-linear theory in order to bound $L^3_{\#}$.

(Note that in (4.1) the $L^2_{\#}$ was bounded by a simple orthogonality argument which does not exploit the geometric structure).

In the separated case, Proposition 2 provides a bound on the inner $L^3_{\#}$-norm in (4.19) by
$$
\Big\Vert\Big(\sum_\tau|F_\tau|^2 \Big)^{\frac 12}\Big\Vert_{L^3_{\#}}\eqno{(4.20)}
$$
with $\{\tau\}$ a partition of $[k\sim K]\times [\ell\sim L]$ in intervals $I_\alpha\times J_\beta$ with $|I_\alpha|=L, |J_\beta|=1$
(we use
here the fact that $L^2\leq K$).

According to Proposition 2$'$, the $\delta$-separated case $(\frac LK <\delta<1)$ involves an additional factor $\delta^{-\frac 12}$.

Next, we bound
$$
\Big\Vert\Big(\sum_\tau|F_\tau|^2\Big)^{\frac 12} \Big\Vert_{L^{2\nu}_{\#}}\eqno{(4.21)}
$$
similarly to (4.2).
Thus (4.21)$^{2\nu}$ amounts to the number of solutions of (4.4), (4.5) where now
$$
|k_1 -k_2|, \ldots, |k_{2\nu-1} - k_{2\nu}|< L\eqno{(4.22)}
$$
instead of (4.6).
Thus a similar calculation as leading to (4.10) gives the bound (using the same notation)
$$
\aligned
&K^\nu L^\nu +K^{\nu-1} L^{\nu-1} U^{\nu-1} +\eta K^{\nu+1} L^{\nu-1} U^{\nu-2}\leq\cr
&K^\nu L^\nu+ K^{\nu-1} L^{2\nu-2}+\eta K^{\nu+1} L^{2\nu-3}.
\endaligned
\eqno{(4.23)}
$$
The natural choice is $\nu=4$, leading to
$$
K^{\frac 12+} L^{\frac 12} (1+\eta KL)^{\frac 18}\eqno{(4.24)}
$$
as bound for the transverse contribution to (4.19).

This contribution of the `$\delta$-separated' case is bounded by
$$
\delta^{-\frac 12} K^{\frac 12+} L^{\frac 12}(1+\eta KL)^{\frac 18}.\eqno{(4.25)}
$$
As before, we are interpolating (4.19) with $\Vert \ \Vert_{L^4_{\#} L^{\frac {24} 5}_{\#}}$
bounded by (3.18), i.e.
$$
K^{\frac 12+} L^{\frac 12}\Big(1+\eta\frac {K^2}L\Big)^{\frac 14}.\eqno{(4.26)}
$$
Note that for the $\delta$-separated contribution in Proposition 4 analyzed below Proposition 3, there is an extra factor
$\delta^{\frac 18}$
that was dropped.
This factor needs to be added to the r.h.s. of (3.18) so that instead of (4.26), one gets in fact
$$
\delta^{\frac 18} K^{\frac 12+} L^{\frac 12} \Big( 1+\eta\frac {K^2}L\Big)^{\frac 14}.\eqno{(4.27)}
$$
Interpolating (4.25), (4.27) with
$$
\frac 1q =\frac {1-\theta}8 +\frac\theta 4=\frac {1-\theta} 3 +\theta .\frac 5{24}
$$
$$
\theta =\frac 56, q=\frac {48}{11}
$$
leads to (4.18).

Above Propositions 10 and 10$'$ form the basis of our treatment of the first spacing problem.
\bigskip

\beginsection
{5.  A variant of the double large sieve}

Let $X, Y\subset\Bbb R^d$ be bounded sets, with
$$
\aligned
&|x_i|< U_i \text { for } x\in X\cr
&|y_i|<V_i \text { for } y\in Y.
\endaligned
\eqno{(5.1)}
$$
Estimate
$$
\sum\Sb x\in X\\ y\in Y\endSb e(x.y) =(5.2)
$$
Let $\nu$ be the discrete measure on $\Bbb R^d$ defined by
$$
\nu =\sum_{x\in X} \delta_x.
$$
By (5.1),
$$
\aligned
(5.2)&= \int_{U_1 \times \cdots \times U_d}\Big[\sum_{y\in Y} e(x.y)\big] \nu(dx)\sim\cr
|(10.2)|&\lesssim \int_ {U_1 \times\cdots \times U_d} \Big|\sum_{y\in Y} e(x.y)\Big| \ |\Bbb E[\nu] (x)|dx
\endaligned
\eqno{(5.3)}
$$
with $\Bbb E[\nu]$ the conditional expectation of $\nu$ at scale $\min(U_1, \frac 1{V_1})\times\cdots \times\min (U_d, \frac
1{V_d})$.
\smallskip

Take $2<q<\infty$.
By H\"older's inequality $(\frac 1p+\frac 1q=1)$
$$
(5.3) \leq \Big[\int_{U_1\times\cdots \times U_d} \Big|\sum_{y\in Y} \, e(x. y)  \Big|^q dx\Big]^{\frac 1q} \ \Vert \Bbb E[\nu]\Vert_p
$$
and by interpolation
$$
\Vert\Bbb E[\nu]\Vert_p \leq |X|^{1-\frac 2q} \Vert\Bbb E[\nu]\Vert_2^{\frac 1q}.
$$
Clearly
$$
\Vert\Bbb E[\nu]\Vert^2_2 = \Big(\frac 1{U_1}+V_1\Big)\cdots \Big(\frac 1{U_d}+V_d\Big). \ (5.4)
$$
with
$$
(5.4) =|\{(x, x')\in X\times X; |x_i-x_i'|< \frac 1{V_i} (1\leq i\leq d)\}|.
$$
Hence we arrive at
$$
\aligned
|(5.2)|\lesssim & (1+U_1V_1)^{\frac 1q}\cdots (1+U_dV_d)^{\frac 1q}. |X|^{1-\frac 2q}\cr
&|\{(x.x')\in X\times X; |x_i-x_i'|<\frac 1{V_i} (1\leq i\leq d)\}|^{\frac 1q}\cr
& \Big[ \nint_{U_1\times\cdots \times U_d}\Big|\sum_{y\in Y} e(x. y)\Big|^q dx\Big]^{\frac 1q}.
\endaligned
\eqno{(5.5)}
$$

\def\cal{\Cal}

\bigskip
\beginsection
{6. An application: bounds for exponential sums with a difference}

Proposition~10 and Proposition~10$'$ supply new information concerning the `First Spacing Problem' of the Bombieri-Iwaniec method for the estimation of exponential sums (see [H96] or [G\&K91] for descriptions of 
the Bombieri-Iwaniec method, and [H96, Part~III] and [W10, Section~3] for relevant previous 
results on the First Spacing Problem). With the aid of the variant of the Bombieri-Iwaniec `Double 
Large Sieve' developed in Section 5 we are able to exploit this new information, 
and so achieve a small but significant advance in the application of the Bombieri-Iwaniec method to a certain class of exponential sums that is of some significance in the analytic theory of numbers. Our results in this direction are contained in the following theorem, the proof of which forms the subject of both the remainder of the present section and the whole of the next six sections.  

\noindent
{\bf Theorem~2.} {\sl
Let $\varepsilon >0$ and $C_2,C_3,\ldots ,C_6\geq 2$ be real constants. Let $\nu\geq 6$ be an integer constant, and let 
$$q_{\nu}={2(13\nu -12)\over 6\nu -5}\;.$$
Let $F(x)$ be a real function that is five times continuously differentiable for
${1\over 3}\leq x\leq 3$,  and let $g(x)$, $G(x)$ be bounded functions of
bounded variation on ${1\over 2}\leq x\leq 1$.
Let $M$ and $T$ be large positive parameters, let $H\geq 1$, and let
$$
S=\sum_{H/2<h\leq H}g\left({h\over H}\right) \sum_{M/2<m\leq M}G\left({m\over M}\right)
{\text e}\left( TF\left({m+h\over M}\right) -TF\left({m-h\over M}\right)\right)\,.
$$ 
Suppose moreover that, on the interval $\left[{1\over 3}\,, 3\right]$, the derivatives $F^{(2)}(x),\ldots , F^{(5)}(x)$ satisfy:  
$$
\left| F^{(r)}(x)\right|\leq C_r\qquad\quad\hbox{($r=2,3,4,5$),}\eqno{(6.1)}
$$
$$
\left| F^{(r)}(x)\right|\geq C_r^{-1}\qquad\quad\hbox{($r=2,3,4$),}\eqno{(6.2)}
$$
and $$\left| F^{(2)}(x)F^{(4)}(x)-3F^{(3)}(x)^2\right|\geq C_5^{-1}\, .\eqno{(6.3)}
$$ 
Then one has the following, in which $B_5$ and $B_4$ are small positive
constants constructed from $C_2,\ldots , C_6$. \hfill\break 
\smallskip\noindent(A)\quad If $H$, $M$ and $T$ satisfy the three conditions 
$$
H\geq M^{-9} T^4 (\log T)^{171\over 140}\qquad\hbox{if}\quad M\leq T^{7\over 16}
(\log T)^{57\over 448},\eqno{(6.4)}
$$
$$
H\geq M^{11} T^{-6} (\log T)^{171\over 140}\qquad\hbox{if}\quad M\geq T^{9\over 16}
(\log T)^{-{57\over 448}},\eqno{(6.5)}
$$ 
$$
H\leq B_5 M T^{-{(149\nu -400)\over 16(29\nu -75)}}(\log T)^{969\nu\over 2240(29\nu -75)},\eqno{(6.6)}
$$
then either 
$$
H\ll M T^{-{149\over 464}} (\log T)^{969\over 64960}\eqno{(6.7)}
$$ 
and 
$$
S\ll T^{\varepsilon} H\min\left\{ \left({M\over H}\right)^{\!\!{277\over 600}} 
T^{397\over 2400} + \left( {H\over M}\right)^{\!\!{19\over 50}} T^{1133\over 2600} \;,\,   
\left( {H\over M}\right)^{\!\!{1\over 25}} 
T^{131\over 400} \right\} ,\eqno{(6.8)}
$$ 
or else 
$$
S\ll H\left({H\over M}\right)^{\!\left( {22\over 25}\right) q_{\nu}^{-1}-{9\over 50}} 
T^{\left( {33\over 100}\right) q_{\nu}^{-1}+{49\over 200}+\varepsilon}\;.\eqno{(6.9)}
$$ 
\medskip\noindent 
(B)\ If $H$, $M$ and $T$ satisfy the two conditions 
$$ 
M\leq C_6 T^{1\over 2},\eqno{(6.10)}
$$
$$
H\!\leq B_4 \min\!\left\{ 
\left( {M^{155\nu -480} 
(\log T)^{\left( {969\over 140}\right) \nu}\over T^{46\nu -160}}\right)^{\!\!{1\over 189\nu -480}} \;,\,   
\left( {M^3\over T}\right)^{\!\!{5\nu -12\over 13\nu -36}}
\right\}\!,\eqno{(6.11)}
$$
\noindent then either 
$$
H\ll\min\left\{ 
M^{155\over 189} T^{-{46\over 189}} (\log T)^{323\over 8820} \;,\; 
M^{15\over 13} T^{-{5\over 13}}\right\}\eqno{(6.12)}
$$
and 
$$
\eqalign{ S &\ll T^{\varepsilon}\min\Bigl\{ T^{13\over 160} M^{125\over 192} H^{179\over 320} +
T^{32\over 153} M^{19\over 51} H^{283\over 612}  
+T^{151\over 520} M^{-{11\over 208}} H^{1473\over 1040} +
T^{113\over 221} M^{-{118\over 221}} H^{276\over 221} \;,\cr 
 &\qquad\qquad\quad\   T^{17\over 80} M^{7\over 32} H^{171\over 160}  
+T^{79\over 204} M^{-{11\over 68}} H^{191\over 204} \Bigr\}\;, }\eqno{(6.13)}
$$
or else 
$$
S\ll  T^{\left( {11\over 20}\right) q_{\nu}^{-1}+{3\over 40}+\varepsilon} 
M^{{9\over 16}-\left( {11\over 8}\right) q_{\nu}^{-1}} 
H^{\left( {33\over 40}\right) q_{\nu}^{-1} +{69\over 80}} 
+ T^{\left( {11\over 51}\right) q_{\nu}^{-1} +{1\over 3}+\varepsilon} 
M^{-\left( {11\over 17}\right) q_{\nu}^{-1}} 
H^{\left( {55\over 51}\right) q_{\nu}^{-1} +{2\over 3}}\;.\eqno{(6.14)}
$$}

\bigskip 

\noindent{\bf Remarks.}\ 

\smallskip 

\noindent{\bf (i)} This theorem is not quite all that one can prove. We have omitted to include in it  
a `Part~(C)' that might be obtained by using `Case~(C)' of Huxley's results (Lemma~10.1 and Lemma~10.2, below) concerning the Second Spacing Problem of the Bombieri-Iwaniec method. Moreover, in (6.6) and (6.11), 
we have chosen to impose an upper bound on $H$ that is slightly stronger than our method requires: 
it would otherwise have been necessary to include in the upper bounds for $S$ certain extra terms 
associated with the perturbing effect of the first three factors of the bound given in (4.12), above. 
Our insistence on the conditions (6.6) and (6.11) may be considered harmless: for it is not 
one of the factors limiting what we are able to achieve through our applications, in 
Section~13, of Theorem~2. Our work in Section~13 is similarly unaffected by the omission of 
a `Part~(C)' from Theorem~2 (it being only cases with $M\ll T^{1/2}$ that are relevant for the applications 
considered in that section). 

\smallskip 

\noindent{\bf (ii)} A preference for simplicity has also led us to simplify the hypotheses on $F(x)$ by 
strengthening them beyond what is strictly necessary. We mention here that 
the condition that (6.2) hold for $r=4$ can be omitted when 
$M>T^{147/328} (\log T)^{-2907/45920}$, and that one can omit the condition (6.3) 
when $M<T^{181/328} (\log T)^{2907/45920}$. These restrictions on the enforcement 
of (6.2) (for $r=4$) and (6.3) are analogous to what occurs in [H03, (1.11) and (1.12)], and 
they have the same origin (in the works [H04] and [H05] of Huxley). Despite what has 
just been noted, we shall work with Theorem~2 as it is stated: this creates one slight 
complication in our proof of Lemma~13.1, where we are obliged to make certain that (6.3) holds. 

\smallskip 

\noindent{\bf (iii)} Although the bound on $S$ in (6.9) becomes stronger as $\nu$ increases, the extent to 
which this can be exploited is limited due to the fact that the upper bound on $H$ in (6.6) also 
strengthens as $\nu$ is increased. This condition (6.6) arises from the assumption that we 
make in (11.13), below. When $\nu\geq 8$, the condition (6.6) requires that we have $H<M T^{-99/314} (\log T)^{969/43960}$, and so, since $99/314=0.315286...\ > 131/416=0.314903...\ $,  
we are prevented from using the corresponding case of (6.9) 
to improve upon the bound $E(T)\ll T^{131/416} (\log T)^{32587/8320}$ obtained in [W10]. 
When one has instead $\nu =7$ the condition (6.6) becomes 
$H\leq B_5 M T^{-643/2048} (\log T)^{969/40960}$, with 
$643/2048=0.313964...\ <131/416$, and so does not prevent us from using the corresponding case of (6.9) 
to improve upon the above mentioned bound for $E(T)$. 
\par 
The assumption (11.13) is made for convenience (it simplifies many calculations). 
It could be replaced by the weaker assumption that $H\ll N^{2/3} R^{1/3}$. 
This would have the effect of replacing the condition (6.6) with the 
condition $H\ll M T^{-247/792} (\log T)^{323/12320}\,$ (which, if one takes the implicit constant to be $B_5$, 
is precisely the case $\nu =6$ of (6.6)). Since $247/792=0.311\dot 8\dot 6$, this 
last upper bound on $H$ is weak in comparison to those mentioned above. 
This relaxation of our assumption (11.13) would, at the same time, lead to the bound (6.9) being weakened 
through the appearance of an extra term arising from the factor $(1+K^{-\nu} L^{2\nu -3} + (\eta K L) K^{2-\nu} L^{2\nu -6})^{1/(26\nu -24)}$  that occurs in (4.12). That is, we would have 
$$S\ll \left(\left({H\over M}\right)^{\!\left( {22\over 25}\right) q_{\nu}^{-1}-{9\over 50}} 
T^{\left( {33\over 100}\right) q_{\nu}^{-1}+{49\over 200}} 
+\left({H\over M}\right)^{\!\!{547\over 350}-\left( {44\over 7}\right) q_{\nu}^{-1}} 
T^{{2341\over 2800} - \left( {59\over 28}\right) q_{\nu}^{-1}}\right) T^{\varepsilon} H$$ 
in place of (6.9). Then, in order that we could obtain the estimate $S\ll M/\log T\,$ (as we do in the result (13.8) of Lemma~13.1, below), we would need the parameter $H$ to satisfy 
$$H\ll M T^{-\max\{ a(q_{\nu}) , b(q_{\nu})\}}\;,$$ 
with  
$$
a(q)={49q+66\over 4(41q+44)}\quad\,{\text and}\quad\,b(q)={2341q-5900\over 8(897q-2200)}\;.$$
Note that $q_{\nu}\,$ (defined in Theorem~2) is an increasing function of $\nu$ on $[6,\infty)$, 
and that $a(q)$ and $b(q)$ are, respectively, decreasing and increasing on $[q_{6},\infty)$. 
A calculation shows that $a(q_7)=1273/4053 =0.3140883...\ > b(q_7)=0.3140809...\ $, whereas 
$a(q_8)=0.31406...\ <b(q_8)=6323/20128=0.31413...\ $. Therefore, even if we had not 
made the assumption (11.13), we would still have had to have 
$c>1273/4053$ in the hypothesis (13.1) of Lemma~13.1: note, in particular, that the size of the extra term that 
would appear in (6.9) could not be reduced by some adjustment of the parameter $N$ (for (10.14) would continue to give the optimal choice of $N$ to use with `Case~(A)' of Lemma~10.1).

\smallskip

\noindent{\bf (iv)} Our proof of Theorem~2 splits naturally into two cases, which are (roughly speaking) that in which $H\ll M T^{-1/3}$, and that in which $H\gg M T^{-1/3}$. 
It is only in the latter case that Proposition~10 and Proposition~10$'$ 
yield new information concerning the first spacing problem of the Bombieri-Iwaniec method. 
In our treatment of the case  `$H\ll M T^{-1/3}$' we use nothing more than some of the bounds for the exponential sum $S$ that were already obtained in [W10]. 
It is convenient to get this case out of the way before beginning any work on the proof of the case `$H\gg M T^{-1/3}$'. Therefore we include in this section the following lemma, from which (via a sequence of straightforward corollaries) we obtain a proof of the case `$H\ll M T^{-1/3}$' of Theorem~2. 

\noindent
{\bf Lemma~6.1.}
{\sl
 Let the hypotheses of Theorem~2 concerning $\varepsilon$, $C_2 , \ldots , C_5$, $F(x)$, $g(x)$, $G(x)$, $M$, $T$ and $H$ be satisfied. 
Then one has the following, in which $B_0$ is a small positive constant constructed from $C_2 , \ldots , C_5$.
 \hfill\break\smallskip

\noindent(A)\ If $H$, $M$ and $T$ satisfy (6.4), (6.5) and the condition  
$$
H\leq B_0 M T^{-{49\over 164}}(\log T)^{{969\over 22960}},\eqno{(6.15)}
$$
then one has 
$$
S\ll_{\varepsilon}\, H\left({H\over M}\right)^{\!\!{1\over 25}} T^{{131\over 400}+\varepsilon}\;.\eqno(6.16)
$$ 
\medskip

\noindent(B)\ If $H$, $M$ and $T$ satisfy the two conditions 
$$
T^{{1\over 3}}\leq M\leq T^{{7\over 16}} (\log T)^{{57\over 448}}\;,\eqno{(6.17)}
$$
$$
H\!\leq\min\!\left\{  B_0 M^{{35\over 69}} T^{-{6\over 69}}  (\log T)^{{969\over 9660}}\;,\;  
B_0 M^{{3\over 2}} T^{-{1\over 2}} \;,\; M^{-9} T^4 (\log T)^{{171\over 140}} \right\}\!,\eqno{(6.18)}
$$
\noindent then one has 
$$
S\ll_{\varepsilon}  T^{{17\over 80}+\varepsilon} M^{{7\over 32}} H^{{171\over 160}}  
+T^{{79\over 204}+\varepsilon} M^{-{11\over 68}} H^{{191\over 204}} \;.\eqno{(6.19)}
$$}

\smallskip
\noindent
{\bf Proof.}\ 
What is stated in this lemma is a slightly weakened and specialized version of what follows immediately  
from [W10, Proposition~1, Parts~(A) and~(B)] if one assumes the case $(\kappa ,\lambda) =(3/10 , 57/140)$ of a certain `Hypothesis $H(\kappa ,\lambda)$' (formulated in [H03, Section~1]): 
note, in particular, that we may assume $B_0\leq 1$, so that the conditions (6.4), (6.5) and (6.15) 
will imply that one has $T^{141/328} (\log T)^{1083/9184}\leq M\leq T^{187/328} (\log T)^{-1083/9184}$, 
which is the case $\kappa =3/10$, $\lambda = 57/140$, $C_6=1$ of [W10, Condition~(1.7)]. 
The lemma therefore follows by virtue of it having been shown, in [W10, Theorem~1], 
that [W10, Proposition~1] remains valid if the first of its hypotheses (to the effect that one has $(\kappa , \lambda)\in\{
(K,L)\in [1/4 , 1/3]\times [0 , \infty) \;:\; {\text Hypothesis}\  H(K,L)\ {\text is\ valid}\}\,$) 
is replaced by the hypotheses that one has $\kappa = 3/10$ and $\lambda =57/140.\quad\blacksquare $ 

\noindent
{\bf Corollary~6.1.1.}
{\sl Let the hypotheses of Theorem~2, up to and including the condition (6.3), 
be satisfied. Suppose moreover that $H$, $M$ and $T$ satisfy the conditions (6.4), (6.5) and (6.6) of 
Part~(A) of that theorem (in which $B_5$ is a certain small positive constant constructed from $C_2 , \ldots , C_5\,$). 
Then the bound (6.16) holds. 

\noindent
{\bf Proof.}\ Since $\nu\geq 6$, we have both 
$${149\nu -400\over 16(29\nu -75)}\geq 
{247\over 792}=0.311\dot 8\dot 6 > 0.29878\ldots\ ={49\over 164}\qquad {\text and}\qquad  
{969\nu\over 2240(29\nu -75)}\leq {969\over 36960}<{969\over 22960}\;,$$ 
and so (assuming, as we may, that $\log T\geq 1$, and that $B_5$ in Theorem~2 is not greater than 
the constant $B_0$ in Lemma~6.1) it follows from (6.6) that the condition (6.15) is satisfied. 
Therefore (given the hypotheses of the corollary) it follows from Part~(A) of Lemma~6.1 that 
we obtain the bound (6.16) for $S.\quad\blacksquare$}

\noindent
{\bf Corollary~6.1.2.} {\sl Let the hypotheses of Theorem~2, up to and including the condition (6.3), be satisfied. 
Suppose moreover that $H$, $M$ and $T$ satisfy the conditions (6.10) and (6.11) in 
Part~(B) of that theorem (in which $B_4$ is a certain small positive constant constructed from 
$C_2,\ldots , C_6$). Then the bound (6.19) holds.} 

\noindent
{\bf Proof.}\ 
Since $\nu\geq 6$, we have ${5\over 13}<{5\nu -12\over 13\nu -36}\leq {3\over 7}$, and so (given that $H\geq 1\geq B_4$) the
condition (6.11) implies that we have both 
$$
M\geq B_4^{-{7\over 9}} T^{{1\over 3}}\geq T^{{1\over 3}}\eqno{(6.20)}
$$ 
and 
$$
H\leq B_4 \left( {M^3\over T}\right)^{\!\!{3\over 7}}\;.\eqno{(6.21)}
$$ 
Furthermore, assuming (as we may) that $\log T\geq 1$, 
it follows by a calculation that if $M\leq T^{1/2}$ then the term 
$( {M^{155\nu -480} (\log T)^{969\nu /140} 
/ T^{46\nu -160}} )^{1/(189\nu -480)}$ that occurs in (6.11) is monotonic 
decreasing, as a function of $\nu\in [6,\infty)$. By this, (6.10), and the point noted in (6.21), 
it follows from (6.11) that we have: 
$$
\eqalignno{ 
H &\leq B_4 \min\left\{ C_6^{{155\over 189}} M^{{75\over 109}} T^{-{58\over 327}} (\log T)^{{969\over 15260}} 
\;,\, \left( {M^3 \over T}\right)^{\!\!{3\over 7}}\right\} \cr 
 &= B_4 \left( {M^3 \over T}\right)^{\!\!{1\over 3}} \min\left\{ C_6^{{155\over 189}} 
(\log T)^{{969\over 15260}} \left( {T\over M^2}\right)^{\!\!{17\over 109}} 
\;,\, \left( {M^3 \over T}\right)^{\!\!{2\over 21}}\right\} &(6.22)\cr 
 &\leq B_4 M T^{-{1\over 3}} C_6^{{16895\over 43092}} 
(\log T)^{{17\over 560}} T^{{17\over 684}} \;.
}
$$
Since $T$ is large, since ${17\over 684}-{1\over 3}=-{211\over 684} = -0.308\ldots\ 
< -0.298\ldots\ =-{49\over 164}$, 
and since we may assume that $B_4 \leq B_0 /C_6$ (where $B_0$ is the constant in Lemma~6.1), 
the above shows that
$$
\hbox{the condition (6.15) is satisfied.}\eqno{(6.23)}
$$ 
Given that  $M^{3/2} T^{-1/2} = (M^3 /T)^{1/2}$, that 
$M^{35/69} T^{-6/69} = (T/M^2)^{17/69} (M^3/T)^{1/3}$,  
and that $1/2>3/7$ and $17/69>17/109$, we are similarly able to deduce from (6.10), (6.21) and (6.22) 
that  
$$
\hbox{the condition (6.18) is satisfied if}\ \,H\leq M^{-9} T^4 (\log T)^{{171\over 140}}\;.\eqno{(6.24)}
$$ 

We observe that if $M\leq T^{7/16} (\log T)^{57/448}$ and $H\leq M^{-9} T^4 (\log T)^{171/140}$  
then it follows by (6.20) and (6.24) that the conditions (6.17) and (6.18) 
are satisfied, so that Part~(B) of Lemma~6.1 yields the bound (6.19). 
If instead $M > T^{7/16} (\log T)^{57/448}$ and $H\leq M^{-9} T^4 (\log T)^{171/140}$  
then, by (6.10) and the corollary of results of Kusmin (or Landau) and Van der Corput that is 
noted in [W10, Equation~(4.8)], one has 
$$
S\ll H\left( H T M^{-2}\right)^{{1\over 2}} M^{{1\over 2}} 
= T^{{1\over 2}} M^{-{1\over 2}}  H^{{3\over 2}} 
<T^{{17\over 80}} M^{{7\over 32}} H^{{171\over 160}}
$$
(as a short calculation shows), and so it is again the case that the bound (6.19) holds. 
These observations show that (6.19) holds whenever $H\leq M^{-9} T^4 (\log T)^{171/140}$.
Therefore we may assume (for the remainder of this proof) that one has 
$$
H > M^{-9} T^4 (\log T)^{{171\over 140}}\;.\eqno{(6.25)}
$$ 

By (6.25) and (6.10) it follows that the conditions (6.4) and (6.5) are satisfied 
(here we assume, as we may, that one has 
$T^{1/16} (\log T)^{-57/448} > C_6\,$). 
By (6.23), we have also (6.15). Furthermore (assuming, again, that $T$ is sufficiently large in terms of $C_6$) 
it follows from (6.10) that the upper bound on $M$ in (6.17) is 
satisfied, while by 
(6.25) and (6.15) (in which we assume $B_0\leq 1$), it follows that we have 
$$M^{-9} T^4 (\log T)^{{171\over 140}} < M T^{-{49\over 164}}(\log T)^{{969\over 22960}}\;,$$ 
which enables us to deduce that the lower bound on $M$ in (6.17) is satisfied, and so to conclude that both of the bounds on $M$ in (6.17) hold. Since each one of the conditions 
(6.4), (6.5) and (6.15) is satisfied, it therefore follows by Part~(A) of Lemma~6.1 that we obtain the 
bound (6.16) for $S$. A calculation shows that (6.16) and (6.25) imply the bound 
$S\ll_{\varepsilon}\, T^{\varepsilon +17/80} M^{7/32} H^{171/160}$, 
which (in turn) implies (6.19). \ $\blacksquare$ 

\noindent
{\bf Corollary~6.1.3.} 
{\sl Let $C$ be a positive constant. 
Suppose that the hypotheses of Theorem~2, up to and including (6.3), be satisfied. 
Then there exist small positive constants $B_5$ and $B_4\,$ (constructed from $C_2 , \ldots , C_6$) 
such that 
the results stated in Parts~(A) and~(B) of Theorem~2 are valid 
whenever $H$, $M$ and $T$ satisfy 
$$
T H^3\leq C M^3\;.\eqno{(6.26)}
$$}

\noindent
{\bf Proof.} \
Let the constants $B_5$ and $B_4$ be the same as in Corollary~6.1.1 and 
Corollary~6.1.2, respectively. 

In considering Part~(A) of Theorem~2, we may assume that the conditions 
(6.4)-(6.6) are all satisfied (otherwise there is nothing to prove). By Corollary~6.1.1 it follows 
that the bound (6.16) holds. Our hypothesis (6.26) implies that 
$(H/M)^{1/25} T^{131/400}\ll (M/H)^{277/600} T^{397/2400}$. By this and (6.16) it follows that  
we obtain the bound (6.8) for $S$.  Since ${149\over 464}=0.32112\ldots\ <{1\over 3}$, and since $T$ is 
assumed to be large, the hypothesis (6.26) also implies that the condition (6.7) is satisfied. 
Both (6.7) and (6.8) have been shown to hold. We therefore find that the result stated in Part~(A) of 
Theorem~2 is valid subject to the condition in (6.26) and the stated hypotheses. 

We now have only to consider Part~(B) of Theorem~2. We assume that the relevant conditions, (6.10) and 
(6.11), are satisfied. It follows by Corollary~6.1.2 that the bound (6.19) holds. 
Since $\nu\geq 6$, it is moreover the case that the remarks concerning the implications of (6.11) 
that were made at the beginning of the proof of Corollary~6.1.2 are still valid in the present context, and 
so, by (6.10) and (6.20), we have also: 
$$
M^2\ll T\leq M^3\;.\eqno{(6.27)}
$$ 
In order to complete this proof we need only show that the bounds for $H$ and $S$ 
in (6.12) and (6.13) hold. Since (6.26) implies $H\ll (M^3 /T)^{1/3}$, and since we have already 
established that (6.19) holds, we are able to verify that (6.12) and (6.13) hold by 
observing that (6.26) and (6.27) imply that one has: 
$$
M^{155\over 189} T^{-{46\over 189}} =\left( {M^3\over T}\right)^{\!\!{1\over 3}} 
\left( {T\over M^2}\right)^{\!\!{17\over 189}}\gg 
\left( {M^3\over T}\right)^{\!\!{1\over 3}}\gg H\;,\ \quad    
M^{15\over 13} T^{-{5\over 13}} =\left( {M^3\over T}\right)^{\!\!{5\over 13}} 
\geq \left( {M^3\over T}\right)^{\!\!{1\over 3}} \gg H\;,$$
$${T^{17\over 80} M^{7\over 32} H^{171\over 160}\over   
T^{13\over 160} M^{125\over 192} H^{179\over 320}} =\left( {M^2\over T}\right)^{\!\!{37\over 960}} \left( {T H^3\over
M^3}\right)^{\!\!{163\over 960}} \ll 1 \ \ {\text and}\ \ 
{T^{79\over 204} M^{-{11\over 68}} H^{191\over 204}\over  
T^{32\over 153} M^{19\over 51} H^{283\over 612}} = \left( {T\over M^3}\right)^{\!\!{37\over 1836}} \left( {T H^3\over
M^3}\right)^{\!\!{145\over 918}} \ll 1.\ \blacksquare$$

\medskip
\noindent
{\bf Remark.}\ 
Note that it is only in Corollary~6.1.3 that we have assumed that $H\ll M T^{-1/3}\,$ 
(that being what the condition (6.26) effectively states). Indeed, we shall later 
make use of Lemma~6.1.1 and Lemma~6.1.2 in dealing with certain cases in 
which one does not have $H\ll M T^{-1/3}\,$ (see below (12.46) for where this occurs). 

\medskip
\beginsection
{7. Initial steps in the application of the Bombieri-Iwaniec method}  

The only cases that need concern us, in completing our proof Theorem~2, 
are those in which $H$, $M$ and $T$ satisfy the condition 
$$H > 32 C_3^{{1\over 3}} M T^{-{1\over 3}}\;.\eqno(7.1)$$
Indeed, whenever this condition is not satisfied we obtain the results of Theorem~2 by virtue of 
the case $C=2^{15} C_3$ of Corollary~6.1.3.  The factor $32 C_3^{1/3}$ that occurs in (7.1) is 
put there in order to ensure that we get the final upper bound seen in (11.8), below. 
We shall assume henceforth that (7.1) holds. 

We shall bound $S$ by applying the Bombieri-Iwaniec method; we follow [W10] in using results of Huxley [H03] on the `Second Spacing Problem' associated with this method, 
but shall modify the approach taken in [W10] in order order to make use of  new results on the `First Spacing Problem' obtained
in Proposition~10 and Proposition~10$'$.  
For the sake of brevity we shall have occasion to refer to steps and intermediate results from the proofs given in [H03], [W04] and [W10] (this seems preferable to repeating the relevant calculations). 

By partial summation it shall suffice to consider the case in which one has 
$$
S=S_F\left( T ; H , H_1 ; M , M_1\right) 
=\sum_{H_1<h\leq H}\sum_{M_1<m\leq M}{\text e}\left( TF\left({m+h\over M}\right) -TF\left({m-h\over
M}\right)\right)\,,\eqno{(7.2)}
$$
with some $H_1\in [H/2,H)$ and some $M_1\in [M/2,M)$. If $\overline{S}$ is substituted for $S$, 
then the function $F$ is effectively replaced by $-F$. By this device we are freed from having to 
consider any case in which $F^{(3)}(x)$ is negative valued on the interval $[1/3,3]$. In the cases where one has $F^{(3)}(x)>0>F^{(2)}(x)$ for $1/3\leq x\leq 3$, we may divide the 
sum $S$ up into 13 similar sums, $S_1,\ldots ,S_{13}$ (say), and can do this in such a way that, within the sum $S_j$, the variable of summation $m$ is constrained to lie in an interval $(M_j',M_j]$ of length not exceeding $M/26$. We may then rewrite $S_j$ by means of a substitution of the form $m=M_j^{\diamond} -m'$. The effect of this is that $F(x)$ is replaced by the function $F_j(x)=-F( (M_j^{\diamond} - M_j^{*} x)/M)$, where $M_j^{*}=M_j^{\diamond} -M_j'+O(1)$. Provided that $M_j^{\diamond}\in{\Bbb Z}$ is suitably chosen, we  
will then have both $F_j^{(3)}(x)>0$ and $F_j^{(2)}(x)>0$ for $1/3\leq x\leq 3$. 
It is moreover possible to ensure that each of the sums $S_1,\ldots ,S_{13}$, when 
rewritten in the way just indicated, will satisfy all of the same conditions as are attached to the sum $S$ in
the statement of Theorem~2 (albeit with each $C_r$ there possibly having to be increased by a certain factor $\Phi_r\in(1,24^r]$). If one or more of the conditions (6.4), (6.6) or 
(6.11) should cease to be satisfied when $M$ is replaced by $M_j^{*}$, then this can be remedied by means of 
the substitutions $F_j =\delta_j F_j^{*}$, $T=\delta_j^{-1} T_j^{*}$, where $\delta_j$ is a suitable constant 
satisfying $1>\delta_j\geq 24^{-4}$. Therefore the only   
cases of Theorem~2 requiring further attention are those 
in which both $F^{(3)}(x)$ and $F^{(2)}(x)$ are positive valued on the interval $[1/3 , 3]$, and so we 
may assume this henceforth. 

In applying the Bombieri-Iwaniec method to $S$  we repeat, with one exception (that being the utilization of [W04, Equations~(2.32) and~(2.33)]), 
the steps described in [W04, Sections~2-5]. These steps assume (from the outset) a fixed choice of parameters $N,R\in{\Bbb N}$ 
satisfying: 
$$
{1\over (R-1)^2}> {2NT\over C_3M^3}\geq {1\over R^2}\;,\eqno{(7.3)}
$$
$$
{HN^2\over MR^2} =O(1)\eqno{(7.4)}
$$
and 
$$
\max\left\{ {H\over R^2} , {R\over H} , {H\over N}, {N\over M}\right\}\leq B_1\;,
\eqno{(7.5)}
$$ 
where the constant $B_1\in (0,1)$ is assumed to be sufficiently small in terms of $C_3$. 

The initial step is a partitioning of the range of the variable of summation $m$ that is achieved 
through a covering of the interval $[M_1-2N,M]$ by a minimal set of disjoint intervals $I_0,I_1,\ldots ,I_{\ell}$, each of length $N$. To each interval $I_i$ there corresponds an `arc' $J_i\subset{\Bbb R}$ that is the image of $I_i$ under the mapping  
$x\mapsto 2T M^{-2} F^{(2)}(x/M)$. These arcs are classified as `major' or `minor', according to the case $r=1$ of the rules set out in [W04]. Some arcs are then fused, so that some minor arcs (and all major arcs) become parts of `long major arcs'. 
For each $i\geq 2$ such that $J_i$ is not part of any long major arc we choose $a/q\in J_{i-2}\cap{\Bbb Q}$ such that $q$ is minimized, subject to the constraints $q\geq R$ and $(a,q)=1$ (the arc $J_{i-2}$ in these cases being minor, though it may, at the same time, be part of a long major arc) and we put $I(a/q)=I_i$. Each such $I(a/q)$ can then be classified as either `bad', or else `good', according to how well $a/q$ can be approximated by rationals of a smaller denominator (see [H03, Page~600] or [W10, Section~2] for details). As a consequence of the results of [W04, Section~3] concerning long major arcs, 
the results of [W04, Equations~(2.16) and~(6.21)] concerning `$q$', 
the bound [W04, Equation~(2.30)] and the case $\eta =0$ of  [W04, Equation~(4.5), (4.22) and~(5.1)],  one has 
$$
|S|\leq  O\left( {MR\log N\over B_2 H^{{1\over 2}} N^{{1\over 2}}}\right) + \sum_{\cal C}\sum_{R\leq Q\leq B_2 H} S\left( {\cal
C}_Q\right),\eqno{(7.6)}
$$
where each ${\cal C}$ is a subset of the set ${\cal I}_{\pi}=\{ I_i\;:\; 2\leq i\leq\ell\ {\text and}\ J_i\ {\text is\ not\ part\ of\ a\ long\ major\ arc}\}$ and each $Q$ is an integer of the form $2^{b-1} R$ (with $b\in{\Bbb N}$), while 
$$
{\cal C}_Q=\{ I(a/q)\in{\cal C}\;:\; Q\leq q<2Q\}\;,
$$ 
$$
S\left( {\cal C}_Q\right)=\sum_{I(a/q)\in{\cal C}_Q}\left| \sum_{H_1<h\leq H}\ \sum_{k\in I(a/q)\cap (M_1,M]} {\text e}\left(
TF\left({k+h\over M}\right) -TF\left({k-h\over M}\right)\right)\right|\eqno{(7.7)}
$$
and $B_2$ denotes a positive constant that is chosen to be sufficiently small (in terms of $C_3$).  
To clarify this we remark that the variable of summation ${\cal C}$ in (7.6)  is subject to the condition 
$$
{\cal C}\in\{ {\cal G}^{[A(j)]}\;:\; j\in{\Bbb N}\}\cup\{ {\cal B} \} ,\eqno{(7.8)}
$$
where  
$$
\eqalign{ A(j) &=2^{j-1}\log N\;,\cr 
{\cal G}^{[A]} &=\left\{ I(a/q)\in{\cal I}_{\pi}\;:\; I(a/q)\ {\text is\ good\ and}\ A>\alpha(a/q;Q')\geq A\min\{ 1/2 , A-\log N\}\right\}\cr 
}
$$
and 
$$
{\cal B}=\left\{ I(a/q)\in{\cal I}_{\pi}\;:\; I(a/q)\ {\text is\ bad}\right\}\;,
$$
while $Q'$ is a certain parameter (to be specified later) and $\alpha(a/q,Q')$ is given by [W10, Equation~(2.1)]. 
Note that, by the relevant definitions, the sets ${\cal B}$ and ${\cal G}^{[A(j)]}\,$ ($j=1,2,\ldots\ $) are pairwise disjoint, and so, in light of our remarks preceding (10.3) (below), 
it follows by [H03, Lemma~2.3]  that one has 
$$
\left| {\cal B}_Q\right|+\sum_{j=1}^{\infty}\left|{\cal G}^{[A(j)]}_{Q}\right| 
=\left|\left\{ I(a/q)\in {\cal I}_{\pi}\;:\; Q\leq q<2Q\right\}\right| 
\ll {M R^2\over N Q^2}\qquad\hbox{($Q\geq R$).}\eqno{(7.9)}
$$ 

In [W04, Equation~(2.30)] it is found (by partial summation) that for each $I(a/q)\in{\cal I}_{\pi}$ one has a bound of the form 
$$\sum_{H_1<h\leq H}\sum_{\scriptstyle k\in I(a/q)\atop\scriptstyle M_1<k\leq M} {\text e}\left( TF\left({k+h\over M}\right) -TF\left({k-h\over M}\right)\right)
\ll\left| 
\sum_{h=H_2}^{H_3}\sum_{n=N_2}^{N_3} {\text e}\left( {(an+b+\kappa)h\over q}+\mu n^2 h +{\mu h^3\over 3}\right)\right|\;,$$
in which $H_2$, $H_3$, $N_2$ and $N_3$ are certain positive integers (dependent on $I(a/q)$) satisfying 
$$H/2\leq H_2\leq H_3\leq H\quad\hbox{and}\quad N\leq N_2\leq N_3\leq 3N\;,$$ 
while $b=b(a/q)\in{\Bbb Z}$, $\kappa=\kappa(a/q)\in [0,1)$ and $\mu =\mu(a/q)>0$ are given by: 
$$b+{\kappa}=2qT M^{-1} F^{(1)}(m/M)\quad\hbox{and}\quad \mu = T M^{-3} F^{(3)}(m/M)\;,$$ 
with $m$ being a nearest integer to the number $x\in (M_1 -2N,M]$ satisfying $2T M^{-2} F^{(2)}(x/M)=a/q$. 
By this, (7.7), (7.9) and the case $r=1$, $\chi=\chi_0$, $\eta =0$ of [W04, Equations~(4.5), Equation~(4.22) and Lemma~5.4] we may conclude that, for $Q\in [R , B_2 H]$ and ${\cal C}$ satisfying (7.8), one has  
$$
S\left( {\cal C}_Q\right) \ll \left( {R\over Q}\right)\left( {H\over N}\right)^{\!\!{1\over 2}} M\log N + \sum_{I(a/q)\in {\cal
C}_Q}\left|\sigma(a/q)\right|\;, \eqno{(7.10)}
$$
where 
$$\sigma(a/q)=\sum_{L_2 <\ell <L_3}\sum_{\scriptstyle k\atop\scriptstyle K_4(\ell)\leq\phi(k-\kappa ,\ell)\leq K_5(\ell)} 
{{\text e}\left( \text{\bf x}^{(a/q)}\cdot \text{\bf y}^{(k,\ell)}(\kappa)\right)\over (4\mu q)^{1/2} \left( (k-\kappa)^2 -\ell^2\right)^{1/4}}\;, \eqno{(7.11)}$$ 
with: 
$$L_j =2\mu q H_j N_j\qquad\hbox{($j=2,3$),}\qquad\qquad\phi(u,\ell)=u-\sqrt{u^2 -\ell^2}\;,$$
$$K_4(\ell)=2\mu q\max\left\{ H_2^2 , \left( \ell /2\mu q N_3\right)^2\right\} ,\qquad 
K_5(\ell)=2\mu q\min\left\{ H_3^2 , \left( \ell/2\mu q N_2\right)^2\right\} ,$$
$$\text{\bf y}^{(k,\ell)}(\kappa)=\left( k\ell , \ell , { ((k-\kappa)+\ell)^{{3\over 2}} -((k-\kappa)-\ell)^{{3\over 2}}\over 3}\right) = \left( k\ell , \ell , \omega(k-\kappa ,\ell)\right)\qquad\hbox{(see (1.1))}$$
and 
$$
\text{\bf x}^{(a/q)}=\left(\left\{ \overline{a\over q}\right\} , \left\{ {\overline{a}b\over q}\right\} , {1\over\sqrt{\mu q^3}}\right)  ,$$
where we have, in the last line above, $\{\beta\} = \beta -\max\{ j\in{\Bbb Z}\;:\; j\leq\beta\}$ (the `fractional part' of $\beta\in{\Bbb R}$), and take 
$\overline{a}$ to be any integer satisfying $a\overline{a}\equiv 1\!\!\!\pmod{q}$. 
Note that we have here implicitly corrected an erroneous statement made in [W04, (4.12)], but  
not propagated to any subsequent part of [W04]. 
For future reference note that, by (6.1), (6.2), (7.3) and (7.5), one has 
$$
{1\over 2 N R^2}\leq\mu < {2 C_3^2\over N R^2}\;.\eqno{(7.12)}
$$

\medskip
\beginsection
{8. Preparations for the the modified double large sieve}

The sum over $I(a/q)$ in (7.10)-(7.11) is not suitable for an immediate application of the Bombieri-Iwaniec double large sieve (for which see [H96, Lemma~5.6.6], for example). 
The principal reason for this is the dependence of the ranges of summation of both $k$ and $\ell$ upon the interval $I(a/q)$. The same problem occurs in [W04, Section~6], and we shall make use of 
one part of the solution given there. However, since we do not indulge here in the averaging over $\eta\in (-1/2,1/2)$ that was found useful in [W04], 
and since it is the modified form of the double large sieve from Section~5 that we seek to apply, 
the preparations the we shall make for its application have to differ in certain other respects from the preparatory steps  undertaken in [W04, Section~6]. 
In particular we shall deal in a different way with terms depending on the variable $\kappa =\kappa(a/q)$. 

Our first concern is with the dependence of the condition $K_4(\ell)\leq\phi(k-\kappa ,\ell)\leq K_5(\ell)$ upon $\kappa$. 
Given any $\ell\in (L_2 ,L_3)$, the set $\{ u\in{\Bbb R}\;:\; K_4(\ell)\leq\phi(u ,\ell)\leq K_5(\ell)\}$ is one subinterval of $(0,\infty)$, and so, for $0\leq\kappa <1$, 
the sums $\sigma(a/q)$ and 
$$\sigma_1 (a/q)
=\sum_{L_2 <\ell <L_3}\sum_{\scriptstyle k\atop\scriptstyle K_4(\ell)\leq\phi(k,\ell)\leq K_5(\ell)} 
{{\text e}\left( \text{\bf x}^{(a/q)}\cdot \text{\bf y}^{(k,\ell)}(\kappa)\right)\over (4\mu q)^{{1\over 2}} \left( (k-\kappa)^2 -\ell^2\right)^{{1\over 4}}}
$$ 
contain less than $2 L_3$ summands that are not common to both. One can show furthermore that within either one of the sums $\sigma(a/q)$,  $\sigma_1(a/q)$ one has 
$$
{HQ\over R^2}\ll \ell\ll {C_3^2 HQ\over R^2}\quad\hbox{and}\quad {N\over H} \ll {k\over \ell}\ll {N\over H}
$$ 
(where the implicit constants are absolute). Therefore, given (7.5) (where one may take $B_1$ to be arbitrarily small) and (7.12), the summands of $\sigma(a/q)$ or  $\sigma_1(a/q)$ have 
absolute values that are bounded above by 
$$
O\left( {1\over\sqrt{\mu q k}}\right) \ll {R^2\over Q}\;,
$$
and so we have: 
$$
\sigma(a/q)-\sigma_1 (a/q) \ll L_3 R^2 Q^{-1} \ll H\;.\eqno{(8.1)}
$$

By elementary calculus, we have also 
$$
\left( (k-\kappa)^2 -\ell^2\right)^{-{1\over 4}}=\left( k^2 -\ell^2\right)^{-{1\over 4}}\left( 1 +O
\left( k^{-1}\right)\right)\qquad\hbox{($0\leq\kappa <1$)}
$$ 
within the sum $\sigma_1(a/q)$. We therefore have: 
$$
\sigma_1 (a/q)-\sigma_2 (a/q)\ll\, L_3\sum_{NQ/R^2\ll k\ll NQ/R^2}\,{R^2\over Qk}\ll H\;,\eqno{(8.2)}
$$ 
where 
$$
\sigma_2 (a/q) 
=\sum_{L_2 <\ell <L_3}\sum_{\scriptstyle k\atop\scriptstyle K_4(\ell)\leq\phi(k,\ell)\leq K_5(\ell)} 
{{\text e}\left( \text{\bf x}^{(a/q)}\cdot \text{\bf y}^{(k,\ell)}
(\kappa)\right)\over (4\mu q)^{{1\over 2}} \left( k^2 -\ell^2\right)^{{1\over 4}}}\;.
$$

Next we work to replace ${\text{\bf y}}^{(k,\ell)}
(\kappa)=\text{\bf y}^{(k,\ell)}(\kappa(a/q))$ by a (higher dimensional) vector that is independent of the interval $I(a/q)$. 
Using the binomial theorem and some elementary estimates, we find that within the sum $\sigma_2(a/q)$ one has  
$$
y^{(k,\ell)}_3 (\kappa) =\omega(k-\kappa ,\ell) 
=O\left( \ell k^{-{5\over 2}}\right) + {1\over 3}\sum_{j=0}^2 \pmatrix 3/2\\ j\endpmatrix (-\kappa)^j 
\left( (k+\ell)^{{3\over 2}-j}-(k-\ell)^{{3\over 2}-j}\right)  ,$$
so that  
$$
x^{(a/q)}_3 y^{(k,\ell)}_3 (\kappa) 
=O\left(  \left( \mu q^3\right)^{-{1\over 2}}\ell k^{-{5\over 2}}\right) + \sum_{j=0}^2 \, {(-\kappa)^j\over j! \left(\mu q^3\right)^{{1\over 2}}} \, \omega_j (k,\ell)\;,
$$
where 
$$\omega_j(u,\ell)={\partial^j\over\partial u^j} \, \omega(u,\ell)\qquad\hbox{($u>\ell >0$).}$$ 
After noting that we will have here
$$\left( \mu Q^3\right)^{-{1\over 2}}\ell k^{-{5\over 2}}   \ll (N/R)^{{1\over 2}} 
(H/N) ( N/R)^{-{3\over 2}}=HR/N^2\leq B_1^2\;,$$
we are able to deduce that 
$$
\sigma_2 (a/q) - \sigma_3 (a/q) \ll \sum_{0<\ell\ll HQ/R^2} \ \sum_{k\gg NQ/R^2} \mu^{-1} Q^{-2} \ell k^{-3} 
\ll N R^2 Q^{-2} (H/N)^2\leq H^2 /N<H\;,\eqno{(8.3)}
$$
where 
$$
\sigma_3 (a/q) 
=\sum_{L_2 <\ell <L_3}\sum_{\scriptstyle k\atop\scriptstyle K_4(\ell)\leq\phi(k,\ell)\leq K_5(\ell)} 
{{\text e}\left( \widetilde {\text{\bf x}}^{(a/q)}\cdot\widetilde 
{\text {\bf y}}^{(k,\ell)}\right)\over (4\mu q)^{{1\over 2}} \left( k^2 -\ell^2\right)^{{1\over 4}}}
$$ 
with 
$$
\widetilde {\text{\bf x}}^{(a/q)} 
=\left(\left\{ \overline{a\over q}\right\} , \left\{ {\overline{a}b\over q}\right\} , {1\over\sqrt{\mu q^3}} , 
{-\kappa\over\sqrt{\mu q^3}} , {\kappa^2\over 2\sqrt{\mu q^3}} \right) 
$$
and 
$$ 
\widetilde {\text{\bf y}}^{(k,\ell)} = \left( k\ell , \ell , \omega(k,\ell) , \omega_1 (k,\ell) , \omega_2 (k,\ell)\right)\;.
$$

The sum $\sigma_3(a/q)$ is now almost suitable for the application of the double large sieve: the sole remaining problem is the dependence upon $I(a/q)$ of the ranges of summation for $k$ and $\ell$. 
By virtue of our elimination of $\kappa$ from the conditions of summation, the problem  just mentioned is essentially the special case $\kappa =0$ of the problem that is addressed in the first two thirds of [W04, Section~6], and can 
therefore be dealt with by employing exactly the same method as is described there. We begin the process by observing that the conditions on the pair $(k,\ell)\in {\Bbb Z}^2$ in the sum $\sigma_3 (a/q)$ are satisfied if and only if 
one has both 
$$
2 H_2 \sqrt{\mu q}\leq \sqrt{k+\ell} -\sqrt{k-\ell}\leq 2 H_3 \sqrt{\mu q}\eqno{(8.4)}
$$ 
and 
$$
2 N_2 \sqrt{\mu q}\leq \sqrt{k+\ell} +\sqrt{k-\ell}\leq 2 N_3 \sqrt{\mu q}\;.\eqno{(8.5)}
$$ 
Similarly to what is found in [W04, Equations~(6.8)-(6.10)], we have now  
$$
\sigma_3 (a/q) 
=\sum_{K_0}\sum_{L_0}\sigma_4\left( a/q ; K_0 , L_0 \right) \;,\eqno{(8.6)}
$$ 
where 
$$
\sigma_4\left( a/q ; K_0 , L_0 \right) 
=\sum_{L_0 <\ell \leq 2L_0}\,{\sum_{K_0 <k\leq 2K_0}}^{\!\!\!\!\!\!\!\!\!\!(a/q)}  
\ {{\text e}\left(\widetilde {\text\bf x}^{(a/q)}\cdot\widetilde {\text{\bf y}}^{(k,\ell)}\right) 
\over (4\mu q)^{{1\over 2}} \left( k^2 -\ell^2\right)^{{1\over 4}}}
$$ 
(with the superfix $(a/q)$ attached to the inner summation indicating that $k$ and $\ell$ are constrained to satisfy (8.4) and (8.5)),  
while $(K_0 , L_0 )$ runs over the pairs of integer powers of $2$ that satisfy
$$
{K_0\over K} , {L_0\over L}\in\left[ {1\over 8} , 144 C_3^2\right]\quad {\text and}\quad {12N\over H}\geq {K_0 \over L_0}\geq
{N\over 4H}\geq 16\;,\eqno{(8.7)}
$$
with 
$$
K={NQ\over R^2}\quad {\text and}\quad L={HQ\over R^2}\;.
$$ 
We may deal with the sum $\sigma_4 ( a/q ; K_0 , L_0)$ in the same way that the sum ${\cal B}_i^{**} (K_0 , L_0 ; \kappa)$, defined in [W04, Equation~(6.11)], is dealt with in [W04, Pages~342-344].
In particular, by means of an application of [W04, Lemma~6.1], it can be shown that one has  
$$
\sigma_4\left( a/q ; K_0 , L_0 \right)  
\ll \int_{-1}^1 \int_{-1}^1 {\left|  \sigma_5 \left( a/q ; K_0 , L_0 , \text{\bf w}\right) \right| 
{\text d}w_1\; {\text d}w_2 \over\Delta\!\left( K_0 , \text{\bf w}\right)  } +H\log N\;,\eqno{(8.8)}
$$
where 
$$
\text{\bf w}=(w_1 , w_2)\;, \qquad 
\Delta\!\left( K_0 , \text{\bf w}\right) 
= \left( K_0^{-2} + |w_1|\right) \left( K_0^{-2} + |w_2|\right) >0
$$ 
and 
$$
\sigma_5 \left( a/q ; K_0 , L_0 , \text{\bf w}\right) 
=\sum_{L_0 <\ell \leq 2L_0}\,\sum_{K_0 <k\leq 2K_0} { 
{\text e}\left(\widetilde {\text{\bf x}}^{(a/q)}\cdot\widetilde {\text{\bf y}}^{(k,\ell)} - K_0^{3/2} \text{\bf w}\cdot
\text {\bf c}^{(k,\ell)}\right)\over (4\mu q)^{{1\over 2}} \left( k^2 -\ell^2\right)^{{1\over 4}}} \;,
$$
with
$$
\text{\bf c}^{(k,\ell)} =\left( \sqrt{k+\ell} -\sqrt{k-\ell} \;,\, \sqrt{k+\ell} +\sqrt{k-\ell}\right)\;.
$$ 

By (7.10), (8.1)-(8.3), (8.6), (8.8) and (7.9) and (7.5), it follows that for $Q\in [R , B_2 H]$ and ${\cal C}$ satisfying (7.8) one has either 
$$
S\left( {\cal C}_Q\right) 
\ll \left( {R\over Q}\right)\left( {H\over N}\right)^{\!\!{1\over 2}} \!\!M\log N 
\leq \left( {H\over N}\right)^{\!\!{1\over 2}} \!\!M\log N \;,\eqno{(8.9)}
$$ 
or else 
$$
S\left( {\cal C}_Q\right) \ll\sum_{K_0} \sum_{L_0} 
\int_{-1}^1 \int_{-1}^1 \left(\sum_{I(a/q)\in {\cal C}_Q} \left|  \sigma_5 \left( a/q ; K_0 , L_0 , \text{\bf w}\right) \right| \right) 
\left(\Delta\!\left( K_0 , \text{\bf w}\right)\right)^{-1} {\text d}w_1\; {\text d}w_2  \;,
$$
where $(K_0 , L_0 )$ runs over the pairs of integer powers of $2$ that satisfy (8.7). 
In the latter case one should observe that the integral $\int_{-1}^1 \int_{-1}^1 (\Delta (K_0 , \text{\bf w}))^{-1} {\text d}w_1
{\text d}w_2$ 
is equal to $(2\log (1+K_0^2 ))^2$, and that 
the relevant number of pairs $(K_0 , L_0)$ does not exceed $O( \log C_3 )$. 
It may therefore be deduced that in that latter case one will have 
$$
S\left( {\cal C}_Q\right)\ll \left( \log K_0\right)^2 \sum_{I(a/q)\in {\cal C}_Q} \left| 
 \sigma_5 \left( a/q ; K_0 , L_0 , \text{\bf w}\right) \right|\;,
$$
for some $\text{\bf w}\in {\Bbb R}^2$ and some pair $( K_0 , L_0)$ satisfying (8.7). 

In order to present our conclusions (just reached) in a form slightly more convenient for the work in the next section, we remark that they trivially imply 
that, if $Q\in [R , B_2 H]$ and ${\cal C}$ satisfies (7.8), then either (8.9) holds, or else there exists some pair $( K_0 , L_0)$ satisfying (8.7) 
and some $\text{\bf W}\in {\Bbb R}^2$ such that 
$$
S\left( {\cal C}_Q\right) \ll (R\log N)^2 Q^{-1} \sum_{I(a/q)\in {\cal C}_Q} \left|  \widetilde\sigma_6 \left( a/q ; K_0 , L_0 , 
\text{\bf W}\right) \right| 
=(R\log N)^2 Q^{-1} \widetilde S\left( {\cal C}_Q\right)\quad\hbox{(say),}\eqno{(8.10)}
$$ 
where 
$$
\widetilde\sigma_6 \left( a/q ; K_0 , L_0 , \text{\bf W}\right) 
=\sum_{L_0 <\ell \leq 2L_0}\,\sum_{K_0 <k\leq 2K_0} \psi_{k,\ell} (\text{\bf W}) {\text e}
\left(\widetilde {\text{\bf x}}^{(a/q)}\cdot\widetilde {\text{\bf y}}^{(k,\ell)}\right) \;,
$$
with 
$$
\psi_{k,\ell} (\text{\bf W}) = \left( {K_0^2 -4L_0^2 \over k^2 -\ell^2 }\right)^{\!\!{1\over 4}} {\text e}
\left( - \text{\bf W}\cdot \text{\bf c}^{(k,\ell)}\right)
$$ 
(so that $\psi_{k,\ell} (\text{\bf W})$ is here independent of $I(a/q)$ and is such that $|\psi_{k,\ell} 
(\text{\bf W})|<1$ when $k/K_0 , \ell /L_0\in (1,2]$). 

It is worth noting that, by (7.5) and (11.11) and (11.12) (below), the rightmost bound in (8.9) is stronger than the bounds  for 
$S({\cal C}_Q)$ that we shall ultimately obtain in (12.11) and (12.14) (below). Therefore, in the course of our proof of (12.11) and (12.14) 
(spanning Sections~9-12, below) we may suppose it to be the case that (8.9) does not hold, and, on the basis of that supposition, may infer from the preceding paragraph that one does have 
the bound (8.10) (with $K_0$, $L_0$ and $\text{\bf W}$ as just described above).

\bigskip
\beginsection
{9. The application of the modified double large sieve}  

Let 
$$
{\cal X}={\cal X}({\cal C} , Q)=\left\{\widetilde {\text{\bf x}}^{(a/q)} \;:\; I(a/q)\in{\cal C}_Q\right\}
$$
and 
$$
{\cal Y}={\cal Y}\left( K_0 , L_0\right) =\left\{\widetilde {\text{\bf  y}}^{(k,\ell)} \;:\; 
K_0 <k\leq 2K_0 ,\, L_0 < \ell\leq 2L_0\ {\text and}\ k,\ell\in{\Bbb Z}\right\}\;.
$$
By virtue of the definitions of 
$\widetilde x^{(a/q)}_1$, $\widetilde x^{(a/q)}_3$, $\widetilde y^{(k,\ell)}_1$ and $\widetilde y^{(k,\ell)}_2$, 
the mappings $I(a/q)\mapsto \widetilde {\text{\bf x}}^{(a/q)}\in{\cal X}$ and $(k,\ell)\mapsto \widetilde {\text{\bf  y}}
^{(k,\ell)}\in{\cal Y}$ 
are injective functions on the domains ${\cal C}_Q$ and $((K_0 ,2K_0 ]\times (L_0 ,2L_0 ])\cap{\Bbb Z}^2$, respectively.
Therefore we may deduce from (8.10) that,  
for a certain pair of functions $\alpha : {\cal X}\rightarrow{\Bbb C}$ and $\beta : {\cal Y}\rightarrow{\Bbb C}$,  
determined by ${\cal X}$, ${\cal Y}$, $K_0$, $L_0$ and $\text{\bf W}$, and satisfying   
$$
\left|\beta(\text{\bf y})\right|<1=\left|\alpha(\text{\bf x})\right|\qquad\hbox{($\text{\bf x}\in{\cal X}$, 
$\text{\bf y}\in{\cal Y}$),}\eqno{(9.1)}
$$
one has:  
$$
0\leq\widetilde S\left( {\cal C}_Q \right)
=\sum_{{\text{\bf x}}\in{\cal X}}\sum_{{\text{\bf y}}\in{\cal Y}} {\text e}(\text{\bf x}\cdot \text{\bf y}) 
\alpha(\text{\bf x})\beta(\text{\bf y})\;.\eqno{(9.2)}
$$ 

In Section~5 it is shown how a modified form of the double large sieve may be be used to obtain useful upper bounds for the
absolute value of a sum similar to the above sum over $(\text{\bf x},\text{\bf y})\in{\cal X}\times{\cal Y}$. 
Only the case in which $\alpha(\text{\bf x})$ and $\beta(\text{\bf y})$ both have range $\{ 1\}$ is treated there, so our next task is to show that this restriction does not prevent us from using the large sieve of Section~5 
to get upper bounds for $\widetilde S({\cal C}_Q)$.

A helpful observation is that one has: 
$$\alpha(\text{\bf x})=\sum_{m=1}^4 (-i)^m\max\left\{ 0 , {\text Re}\left( i^m \alpha(\text{\bf
x})\right)\right\}\qquad\hbox{($\text {\bf x}\in{\cal X}$),}$$ 
and (of course) a similar formula for each $\beta(\text{\bf y})$ occurring in (9.2). By using these formulae, and a change in the order of summation, one rewrite the sum on the right-hand side of (9.2) 
as a sum of sixteen sums that are each similar to the original sum, but have a product of the form $(-i)^{m+n}\max\{ 0 , {\text
Re}( i^m \alpha(\text{\bf x}))\}\max\{ 0 , {\text Re}( i^n \beta(\text{\bf y}))\}$ in 
place of $\alpha(\text{\bf x})\beta(\text{\bf y})$. From this and (9.1) we may deduce that there exist functions 
$$\alpha_1 : {\cal X}\rightarrow [0,1]\;,\qquad \beta_1 : {\cal Y}\rightarrow [0,1]$$
such that 
$$0\leq \widetilde S\left( {\cal C}_Q\right)\leq 16\left|\sum_{\text{\bf x}\in{\cal X}}\sum_{\text{\bf y}\in{\cal Y}} {\text e}
(\text{\bf x}\cdot \text{\bf y}) \alpha_1(\text{\bf x})\beta_1(\text{\bf y})\right|\;.$$ 
Within the last sum over $(\text{\bf x},\text{\bf y})\in{\cal X}\times{\cal Y}$ we may apply the substitutions 
$$\alpha_1 (\text{\bf x})=\int_0^{\alpha_1 (\text{\bf x})}{\text d}\theta\quad{\text and}\quad \beta_1 (\text
{\bf x})=\int_0^{\beta_1 (\text {\bf x})}{\text d}\phi\;.$$
Then, via a change in the order of summation and integration, we are able to deduce that 
$$
0\leq \widetilde S\left( {\cal C}_Q\right)\leq 16\int_0^1 \int_0^1\Biggl|\sum_{\scriptstyle{\text{\bf x}}\in{\cal X}\atop\scriptstyle 
\alpha_1 (\text{\bf x})\geq\theta} 
\sum_{\scriptstyle {\text{\bf y}}\in{\cal Y}\atop\scriptstyle \beta_1 (\text{\bf y})\geq\phi} {\text e}(\text{\bf x}\cdot 
\text{\bf y})\Biggr|{\text d}\theta\,{\text d}\phi\;.
$$ 
Since the integral here does not exceed the maximum value attained by its integrand, it follows that there exist subsets 
$$
{\cal X}_1\subseteq{\cal X}\quad{\text and}\quad {\cal Y}_1\subseteq{\cal Y}\eqno{(9.3)}
$$ 
such that one has 
$$
0\leq \widetilde S\left( {\cal C}_Q\right)\leq 16\left| S^* \left( {\cal X}_1 , {\cal Y}_1\right)\right|\;,\eqno{(9.4)}
$$ 
with
$$
S^* \left( {\cal X}_1 , {\cal Y}_1\right) 
=\sum_{{\text{\bf x}}\in{\cal X}_1}\sum_{{\text{\bf y}}\in{\cal Y}_1} {\text e}(\text{\bf x}\cdot \text{\bf y})\;.
$$ 

Let $V$ be chosen so that one has 
$$
V\geq 1\eqno{(9.5)}
$$ 
(this choice will ultimately be determined by our use of Lemma~10.1 and Lemma~10.2, below, and so $V$ will depend on the set ${\cal C}$, and on which of (10.14) or (10.15) is taken as the definition of $N$). 
Then, by (7.12), (8.7), (9.3) and the definitions of the sets ${\cal C}_Q$, ${\cal X}$, ${\cal Y}$ and their elements, we have:
${\cal X}_1,{\cal Y}_1\subset{\Bbb R}^5$ and  
$$
\left| x_i\right| <D_i\quad{\text and}\quad\left| y_i\right| < E_i^{[V]}\qquad\hbox
{($\text{\bf x}\in{\cal X}_1$, $\text{\bf y}\in{\cal Y}_1$, $i=1,\ldots ,5$),}\eqno{(9.6)}
$$ 
where 
$$
\text{\bf D}=\left( 1 , 1 , \sqrt{ 2NR^2 \over Q^3} , \sqrt{ 2NR^2 \over Q^3} , \sqrt{ NR^2 \over 2Q^3} \right)\quad{\text and}\quad  
\text{\bf E}^{[V]}=\left( 5K_0 L_0 V , 3L_0 , 3L_0 \sqrt{K_0} , {2L_0 \over \sqrt{K_0}} , {L_0 \over\sqrt{ K_0^3}}\right)
.\eqno{(9.7)}
$$ 
Given we assume $R\leq Q\leq B_2 H$, it follows by (7.5), (8.7), (9.5) and (9.7) that 
$$
\prod_{i=1}^5 \left( 1+D_i E_i^{[V]}\right) \ll\left( {HNQ^2 V\over R^4}\right)\left( {HQ\over R^2}\right)\left( {HN\over R^2}\right)\left( {H\over Q}\right) (1)
={H^4 N^2 Q^2 V\over R^8}\;.
$$ 
By this and (9.6), an application of the large sieve inequality (5.5) shows that we have 
$$
S^* \left( {\cal X}_1 , {\cal Y}_1\right) 
\ll \, \left(   A_p\left( {\cal Y}_1 , \text{\bf D}\right)  \breve b\bigl( {\cal X}_1 , \text{\bf E}^{[V]}\bigr) H^4 N^2 Q^2 V R^{-8}\right)^{{1\over p}} \left| {\cal X}_1\right|^{1-{2\over p}}\qquad 
\hbox{($2<p<\infty$),}
$$
where 
$$
A_p\left( {\cal Y}_1 , \text{\bf D}\right) ={1\over 2^5 D_1 D_2\,\cdots\, D_5}\int_{-D_1}^{D_1}\ \cdots\ \int_{-D_5}^{D_5} 
\left| \sum_{{\text{\bf y}}\in{\cal Y}_1} {\text e}(\text{\bf x}\cdot\text{\bf y})\right|^p {\text d}x_5\,\cdots \,{\text d}x_1
$$ 
and 
$$
\breve b\bigl( {\cal X}_1 , \text{\bf E}^{[V]}\bigr) 
=\left|\left\{ \left( \text{\bf x} , \text{\bf x}'\right)\in{\cal X}_1 \times{\cal X}_1 \;:\; \left| 
x_j -x_j'\right| < 1/E^{[V]}_j\ \,\hbox{($j=1,\ldots ,5$)}\right\}\right|\;.
$$
Now, by (9.3) and the definition of the set ${\cal X}$, we have 
$|{\cal X}_1|\leq |{\cal X}| = | {\cal C}_Q |$. Therefore, by the bound just obtained for $S^* \left( {\cal X}_1 , {\cal Y}_1\right)$, in combination with (8.10) and (9.4), one has: 
$$
S\left( {\cal C}_Q\right) 
\ll \left(   A_p\left( {\cal Y}_1 , \text{\bf D}\right) \breve b\bigl( {\cal X}_1 , \text{\bf E}^{[V]}\bigr) V H^4 N^2 R^{-6}\right)^{{1\over p}} \left(\left| {\cal C}_Q\right| R/Q\right)^{1-{2\over p}} R (\log N)^2\qquad 
\hbox{($2<p<\infty$).}\eqno{(9.8)}
$$

\bigskip
\beginsection
{10. The second spacing problem}  

Since there is a one-to-one correspondence between the elements of ${\cal C}_Q$ and those of ${\cal X}$ , 
it therefore follows by (9.3) and the definitions of both ${\cal X}$ and $\breve b\bigl( {\cal X}_1 , \text
{\bf E}^{[V]}\bigr)$ that one has 
$$
\eqalignno{\breve b\bigl( {\cal X}_1 , \text{\bf E}^{[V]}\bigr) &\leq 
\left|\left\{ \left( \text{\bf x} , \text{\bf x}'\right)\in{\cal X} \times{\cal X} \;:\; \left| x_j -x_j'\right| < 1/E^{[V]}_j\ \,\hbox{($j=1,\ldots ,4$)}\right\}\right| \cr 
 &=\left|\left\{ \left( I(a/q) , I(a'/q')\right)\in {\cal C}_Q\times {\cal C}_Q \;:\; \left| \widetilde x_j^{(a/q)} -\widetilde x_j^{(a'/q')}\right| < 1/E^{[V]}_j\ \,\hbox{($j=1,\ldots ,4$)}\right\}\right| \cr 
 &\leq B\bigl( {\cal C}_Q ; V\bigr)\;, &(10.1) }
$$
where (in light of (9.7), (9.5), (8.7), (7.12) and (7.5)) one may take $B\bigl( {\cal C}_Q ; V\bigr)$ to be the number of pairs of intervals $I(a/q),I(a'/q')\in {\cal C}_Q$ satisfying a 
system of inequalities of the form 
$$
\left\| \overline{a\over q} - \overline{a'\over q}\right\|\leq\Delta_1\;,\qquad 
\left| {\mu' {q'}^3\over\mu q^3} -1\right|\leq\Delta_2\;,\qquad \left\| {\overline{a} b\over q} - 
{\overline{a'} b'\over q'}\right\|\leq\Delta_3\quad {\text and}\quad \left| \kappa - \kappa'\right| \leq\Delta_4\;,
$$
in which $\| \alpha\|=\min\{ |\alpha - n|\;:\; n\in{\Bbb Z}\}$, while the numbers $\Delta_1 ,\ldots ,\Delta_4$ are determined by $V$, $K_0$, $L_0$, $R$, $N$ and $C_3$, and satisfy: 
$$
\Delta_1\asymp {R^4\over HNQ^2 V}\;, \qquad 
\Delta_2\asymp {R^2\over HN} ,\qquad 
\Delta_3\asymp {R^2\over HQ}\quad {\text and}\quad 
\Delta_4\asymp {Q\over H}\eqno(10.2)
$$
(with the notation $X\asymp Y$ signifying that one has both $X\ll Y$ and $Y\ll X$). 

The problem of obtaining good upper bounds for $B\bigl( {\cal C}_Q ; V\bigr)$ is essentially the same `Second Spacing Problem' as that referred to in [H03, Section~3] 
(see also [H05, Section~3] for a somewhat generalized definition of this problem). Indeed, the only difference between the two that is of any significance is that, whereas the 
function $2 F^{(1)} (x)$ and its derivatives play a certain part in determining the second spacing problem in this present paper (i.e. they play their part in determining the set ${\cal X}$), the 
corresponding part in [H03] is played instead by the function there named $F(x)$, and its derivatives.
 The only consequence of this difference is that, where  a condition of (for example) the form $C\geq |F^{(r)}(x)|\geq 1/C$  is assumed in [H03], we shall instead need only an assumption implying that  
one has $C'\geq 2|F^{(r+1)}(x)|\geq 1/C'$, for some constant $C'\geq 1$. Therefore each of the results on the second spacing problem that are stated in 
[H03, Section~3] implies a similar result for $B\bigl( {\cal C}_Q ; V\bigr)$ (differing only in that the hypotheses concerning derivatives of $F(x)$ are modified in the way that our preceding remark indicates). 
For the same reason we are able to infer from [H03, Lemmas~2.3, 2.4 and~2.5] certain bounds for the number of elements in each set ${\cal C}_Q$ that we need consider. One of these bounds is (7.9) (above). 
The other two assume more about how the classification of elements of ${\cal I}_{\pi}$ (as being either good or bad) is done. That classification 
is dependent on a pair of chosen parameters, $\eta$ and $Q'$. If these chosen parameters satisfy 
$$
1+{M^2\over T}\leq Q' <\eta R<R\eqno{(10.3)}
$$ 
then one may infer from [H03, Lemma~2.3, Lemma~2.4 and Lemma~2.5] that one has 
$$
\left| {\cal B}_Q\right|\ll {\eta M R^2\over N Q^2 } + {M Q'\over N Q}\qquad\hbox{($R\leq Q\leq B_2 H$)}\eqno{(10.4)}
$$ 
and 
$$
\left| {\cal G}^{[A]}_Q\right|\ll {M R^2 \log N \over A N Q^2}\qquad\hbox{($R\leq Q\leq B_2 H$ and $A\geq\log
N$).}\eqno{(10.5)}
$$ 

Huxley's unconditional results in [H03, Lemma~3.4] are an outcome of his work in [H05] on `resonance curves'. 
It is assumed in [H03, Lemma~3.4] that one has  
$$
V_0^2\ll N\sqrt{ V_1 V_2}\;, \eqno{(10.6)}
$$ 
with 
$$
V_0 =\left( {H\over R}\right)^{\!\!{18\over 17}},\qquad V_1 ={R^4\over HN}\quad {\text and}\quad V_2 ={M^2\over H N^3}\;.\eqno{(10.7)}
$$
Note that (10.6) is [H03, Condition~(3.12)]. The parameters $V_1$ and $V_2$ have a significance that is explained below [H03, Equation~(3.5)]. 
There is also another parameter $\Delta'$ that plays a part within certain calculations of [H05]. In [H03, Lemma~3.4] it is assumed that 
one may assign $\Delta'$ a value such that 
$$
\left( {R\over H}\right)^{\!\!{4\over 17}}\ll {\Delta'\over\Delta_2} 
\ll\left( {R\over H}\right)^{\!\!{4\over 119}}\;,\eqno{(10.8)}
$$ 
$$
{Q'\over R}\ll \left( {\Delta' R\over\Delta_2 H}\right)^{\!\!{1\over 3}} \eta\;, \eqno{(10.9)}
$$
$$
V_0\gg {\Delta_2 R^2 \eta\over\Delta' {Q'}^2} \eqno{(10.10)}
$$ 
and 
$$
\left( {\Delta'\over\Delta_2}\right)\left( {Q'\over R}\right)^{\!\!6} 
\gg \left( {R\over H}\right)^{\!\!{70\over 17}}\;.\eqno{(10.11)}
$$ 

>From [H03, Lemma~3.4] we infer, as an immediate corollary, the following lemma.

\medskip
\noindent
{\bf Lemma~10.1.  (Huxley).} {\sl Let those of the hypotheses of Theorem~2 that concern $F(x)$ be satisfied. Suppose also that (7.3)-(7.5), (10.2), (10.3) and (10.6)-(10.11) hold. 
Put 
$$ 
Q_2= {R^{{80\over 119}} H^{{39\over 119}}\over (\log (H/R))^{{3\over 4}}}\;.\eqno{(10.12)}
$$ 
Then one has 
$$
B\bigl( {\cal G}^{[A]}_Q ; V\bigr) \ll {M\over N} 
+\left( 1 +{Q\over Q_2}\right)\left( {M^2 R^4 (\log N)^{{1\over 4}} A^{{11\over 70}} 
\over H^2 N^4 V_0^{{2\over 3}} V}\right) 
\qquad\hbox{($R\leq Q\leq B_2 H$ and $A\geq\log N$)}\eqno{(10.13)}
$$ 
in the following three independent cases: \hfill\break 
\vskip -2mm
\indent (A) when one has $V=V_0\ll\min\{ V_1 , V_2\}$; \hfill\break 
\indent (B) when $M^2\ll T$, with $V=\max\{ V_1 , 1\}$;\hfill\break 
\indent (C) when $M^2\gg T$, with $V=\max\{ V_2 , 1\}$. }
\medskip

\noindent
{\bf Remark.}\  
In the present work (where our main concern is with bounds for exponential sums that are of use in estimating $E(T)$) we shall employ (10.13) only in the cases~(A) and~(B) 
(hence the lack of any `Case~(C)' in our statement of Theorem~2). 

\bigskip 

The choice of $N$  must (of course)  be made prior to the application of Lemma~10.1, and the case of the lemma that applies will depend on that choice. 
Nevertheless one can take account of the form of the bound (10.13), and the definitions of the cases of the lemma in optimizing that choice. We shall restrict 
our choice of $N\in{\Bbb N}$ so as to have either 
$$
N\asymp H\left( {M\over H}\right)^{\!\!{41\over 25}} T^{-{49\over 100}} (\log T)^{{969\over 14000}}\;,\eqno{(10.14)}
$$ 
or else 
$$
N\asymp \min\left\{ {M^{{7\over 8}} (\log T)^{{969\over 5600}} \over T^{{3\over 20}} H^{{29\over 40}}} 
\; , \, {B_3 M^2\over T^{{2\over 3}} H^{{1\over 3}}}\right\}\;,\eqno{(10.15)}
$$ 
where $B_3$ is a sufficiently small positive constant (constructed from $C_2,\ldots ,C_6$). These are essentially the same choices for $N$ that are described in 
in [H03, Equations~(3.19)-(3.21)] (except that we specialize to the case $\kappa =3/10$, $\lambda =(1/4)+(11/70)=57/140$ of what is stated there). 
We allow the choice (10.14) only if that choice results in Case~(A) of Lemma~10.1 being applicable; the choice (10.15) is similarly associated with Case~(B) of the lemma 
(and so is permissible only when $M^2\ll T$). 
For certain combinations of values of $T$, $M$ and $H$, both the options (10.14) and (10.15) may be available (we shall then consider what is the outcome from each of the two choices of $N$). 

These choices for $N$ (and the associated restrictions on the use of Lemma~10.1) are exactly what is required in order to ensure that we never find the term $M/N$ on the right-hand side of (10.13) dominating the other term there. 
For this reason we obtain, every time, bounds of the form 
$$
B\bigl( {\cal G}^{[A]}_Q ; V\bigr) V \ll\left( 1 +{Q\over Q_2}\right)\left( 
{M^2 R^4 (\log N)^{{1\over 4}} A^{{11\over 70}}\over H^2 N^4 V_0^{{2\over 3}}}\right) ,\eqno{(10.16)}
$$ 
with 
$$
V=V\bigl( {\cal G}^{[A]}\bigr) =\cases V_0 &\text {if (10.14) is assumed}; \\ \max\{ V_1 , 1\} &\text{if instead (10.15) is
assumed}\endcases.
\eqno{(10.17)}
$$ 

The association of (10.15) with Case~(B) of Lemma~10.1 also has the effect of ensuring that we  have 
$$\min\left\{ V_1 , V_2\right\}\gg 1\;.\eqno(10.18)$$
This bound is immediate when Case~(A) of the lemma applies (for (7.5) and (10.7) imply that one has $V_0 >1$); when Case~(B) of the lemma applies one obtains (10.18) by virtue of  
(10.7), (10.15), (7.3) and (7.5). 

We require, in addition to the bound (10.13), some sufficiently strong bounds for the numbers $B\bigl( {\cal B}_Q ; V\bigr)\,$ ($R\leq Q\leq B_2 H$). 
For this we fall back on the following immediate corollary of [H03, Lemma~3.2].

\medskip
\noindent
{\bf Lemma~10.2. (Huxley).}  
{\sl Let those of the hypotheses of Theorem~2 that concern $F(x)$ be satisfied. Suppose also that (7.3)-(7.5), (10.3) and (10.7) hold. 
Put 
$$
Q_3= {R^{{2\over 3}} H^{{1\over 3}}\over \log (H/R) }\;.\eqno{(10.19)}
$$ 
Then one has 
$$
B\bigl( {\cal B}_Q ; V\bigr) \ll {M\over N} 
+\left( 1 +{Q\over Q_3}\right)\left( {M^2 R^4 \over H^2 N^4 (H/R)^{{2\over 3}} V}\right) 
\qquad\hbox{($R\leq Q\leq B_2 H$)}\eqno(10.20)
$$ 
in the following three independent cases: \hfill\break 
\vskip -2mm
\indent (A) when one has $V=H/R\ll\min\{ V_1 , V_2\}$; \hfill\break 
\indent (B) when $M^2\ll T$, with $V=\max\{ V_1 , 1\}$;\hfill\break 
\indent (C) when $M^2\gg T$, with $V=\max\{ V_2 , 1\}$. 
}
\medskip 

Note that the conditions defining Case~(B) in this lemma are the same as those 
defining Case~(B) in Lemma~10.1. Since $V_0$ is greater than $H/R$ it may also be observed that, whenever the choice $V=V_0$ would cause   
the conditions of Case~(A) of Lemma~10.1 to be satisfied, the alternate choice $V=H/R$ would  ensure (instead) that the conditions of Case~(A) of Lemma~10.2 are satisfied. 
Consequently,  whenever a bound of the form (10.16)-(10.17) is obtained (in the manner indicated above) it will follow from Lemma~10.2 that one also obtains the bound 
$$
B\bigl( {\cal B}_Q ; V\bigr) V \ll\left( 1 +{Q\over Q_3}\right)\left( {M^2 R^4 (\log N)^{{57\over 140}} 
\over H^2 N^4 (H/R)^{{2\over 3}}}\right) ,\eqno{(10.21)}
$$ 
with 
$$
V=V( {\cal B} ) =\cases H/R & {\text {if (10.14) is assumed}}; \\ \max\{ V_1 , 1\} & {\text {if instead (10.15) is assumed}.}
\endcases \eqno{(10.22)}
$$ 
\bigskip

\beginsection
{11. The first spacing problem: bounds from decoupling for perturbed cones}

>From (9.7) and the definition of $A_p ( {\cal Y}_1 , \text{\bf D})$  
in Section~9 it may be deduced that, for some point $(\xi_3 , \xi_4 , \xi_5)\in{\Bbb R}^3$ satisfying 
$-D_j\leq\xi_j\leq D_j\,$ ($j=3,4,5$), one has 
$$
0\leq A_p \bigl( {\cal Y}_1 , \text{\bf D} \bigr) 
\leq {1\over D_3}\int_{-1}^1 \int_{-1}^1  \int_{\xi_3}^{\xi_3 + \left( {1\over 4}\right) D_3} 
\left|\sum_{{\text{\bf y}}\in {\cal Y}_1} {\text e}\left(\bigl( x_1 , x_2 , x_3 , \xi_4 , \xi_5\bigr)\cdot \text{\bf y}\right)
\right|^p {\text d}x_3\,{\text d}x_2\,{\text d}x_1\;.
$$ 
By this, (9.3) and the definitions relating to the set ${\cal Y}$, one obtains bounds of the form 
$$
\eqalign{ 
 &{1\over D_3}\int_{-1}^1 \int_{-1}^1  \int_{-\left( {1\over 8}\right) D_3}^{\left( {1\over 8}\right) D_3} 
\left|\sum_{K_0 < k\leq 2K_0}\,\sum_{L_0 <\ell\leq 2L_0} a_{k , \ell}\, {\text e}
\left( x_1 k\ell + x_2 \ell + x_3 \omega (k , \ell)\right)\right|^p {\text d}x_3\,{\text d}x_2\,{\text d}x_1 \cr  
 &\quad\geq A_p \bigl( {\cal Y}_1 , \text{\bf D} \bigr) \geq 0\;,} 
$$ 
with certain complex coefficients $a_{k,\ell}$ that are independent of $(x_1 , x_2 , x_3 )\in{\Bbb R}^3$, and that satisfy  
$$
\left| a_{k,\ell} \right|\leq 1\qquad\hbox{($K_0 < k\leq 2K_0$ and $L_0 <\ell\leq 2L_0$).}\eqno{(11.1)}
$$ 
Therefore, upon recalling the notation from Section 4 (and renumbering $x_1$ and $x_2$), we arrive at
$$
\left( A_p \bigl( {\cal Y}_1 , \text{\bf D} \bigr) \right)^{\!{1\over p}} 
\leq \biggl\| \sum_{k\sim K_0 ,\,\ell\sim L_0} a_{k , \ell}\, {\text e}\left( x_1 \ell + x_2 k\ell + x_3 \omega (k , \ell)\right) 
\biggr\|_{L^p_{\#}\bigl[ |x_1|<1 , |x_2|<1 , |x_3|<\left( {1\over 8}\right) D_3\bigr]}\;.\eqno{(11.2)}
$$

By (9.7) one has 
$$
{D_3\over 8} =\left( {NR^2\over 32 Q^3}\right)^{\!\!{1\over 2}} ={1\over\sqrt{K_0} L_0\eta}\quad\hbox{(say),}\eqno{(11.3)}
$$ 
where, by (8.7), 
$$
\eta =\left( {32 Q^3\over NR^2}\right)^{\!\!{1\over 2}} K_0^{-{1\over 2}} L_0^{-1} 
= {(Q/R)^2\over\left( K K_0 /32\right)^{{1\over 2}}  L_0}\;.
$$
By (8.7) (again) it follows that 
$$
\eta\asymp {(Q/R)^2 \over K_0 L_0}\eqno{(11.4)}
$$ 
and, in particular, that one has: 
$$\eqalignno{\eta &> {1\over K_0 L_0}\qquad\hbox{($Q\geq R$),} &(11.5)\cr 
\eta &< {1\over K_0}\qquad\hbox{($Q\leq B_2 H\leq H/64$, say)} &(11.6)}$$ 
and 
$${1\over K_0}< {1\over 8 L_0}\;.\eqno(11.7)$$ 
Moreover, by (8.7), (7.3), (7.5) and (7.1) , one has also: 
$$
\eqalignno{\left( K_0 / L_0\right)^2 \eta &=2^{{5\over 2}} (Q/R)^2 K^{-{1\over 2}} 
\left( K_0 /L_0\right)^{{3\over 2}} L_0^{-{3\over 2}} \cr 
 &\leq 2^{{5\over 2}} (Q/R)^2 K^{-{1\over 2}} (12K/L)^{{3\over 2}} (8/L)^{{3\over 2}} \cr 
 &=3072\sqrt{3} NR^2 /H^3 \cr 
 &<6144\sqrt{3} C_3 (M/H)^3 T^{-1} <1\;. &(11.8)}
$$ 

Assuming that we have 
$$
R\leq Q\leq B_2 H\leq H/64\;,\eqno(11.9)
$$
it follows by (11.1) and (11.5)-(11.8) that an application of  Proposition~10 yields the upper bound
$$\eqalign{ &\biggl\| \sum_{k\sim K_0 ,\,\ell\sim L_0} a_{k , \ell}\, {\text e}\left( x_1 \ell + x_2 k\ell + x_3 \omega (k , \ell)\right) 
\biggr\|_{L^{q_{\nu}}_{\#}\bigl[ |x_1|<1 , |x_2|<1 , |x_3|<{1\over\sqrt{K_0} L_0 \eta}\bigr]} \cr 
&\ \ll_{\nu ,\varepsilon}\, 
\left( 1 +{\eta K_0^2\over L_0}\right)^{3(\nu -1)\over 13\nu -12} \left( 1 +{L_0^3\over K_0^2}\right)^{\nu-1\over 2(13\nu -12)} 
\left( 1+{L_0^{2\nu -3}\over K_0^{\nu}} +{\left(\eta K_0 L_0\right) L_0^{2\nu -6} \over 
K_0^{\nu -2}}\right)^{1\over 2(13\nu -12)} K_0^{\varepsilon +{1\over 2}} L_0^{{1\over 2}}}\eqno{(11.10)}
$$ 
for all $\varepsilon >0$ and all pairs $(\nu , q_{\nu})\in{\Bbb R}^2$ such that 
$$
\nu\in{\Bbb Z}\;,\qquad \nu\geq 3\quad{\text and}\quad q_{\nu}={2(13\nu -12)\over 6\nu -5}\;.\eqno{(11.11)}
$$

We assume henceforth that $\nu$ and $q_{\nu}$ are as in (11.11), and that one also has both 
$$
\nu\geq 6\eqno{(11.12)}
$$ 
and 
$$
H\ll N^{\nu -2\over 2\nu -6} R^{\nu -4\over 2\nu -6}\;.\eqno{(11.13)}
$$ 
Note that, by our assumptions in (7.5), the right-hand side of (11.13) is a monotonic decreasing function of $\nu$, and so (11.12) and (11.13) certainly imply that one has $H \ll N^{2/3} R^{1/3}$. 
By (11.4), (8.7) and (11.9) we have also  
$$
1 +{\eta K_0^2\over L_0}\asymp 1+{(Q/R)^2 K\over L^2} = 1+{NQ\over H^2}\;,
$$
$$
1 +{L_0^3\over K_0^2}\asymp 1 +{L^3\over K^2}= 1+ {H^3 Q\over N^2 R^2}
$$ 
and, since $H^3\ll N^2 R$,  
$$
\eqalign{1+{L_0^{2\nu -3}\over K_0^{\nu}} +{\left(\eta K_0 L_0\right) L_0^{2\nu -6} \over K_0^{\nu -2}}
 &\asymp 1+\left( {(Q/R)^2 L^{2\nu -6} \over K^{\nu -2}}\right)\left( 1+{(R/Q)^2 L^3\over K^2}\right) \cr 
 &\leq 1+\left( {H^{2\nu -6} Q^{\nu -2}\over R^{2\nu -6} N^{\nu -2}}\right) \left( 1 + {H^3\over N^2 R}\right) \cr 
 &\asymp 1+{H^{2\nu -6} Q^{\nu -2}\over R^{2\nu -6} N^{\nu -2}}\;.}
$$
Therefore, subject to the assumptions made, the bound (11.10) implies that one has, for $\varepsilon >0$, 
$$
\eqalign{ &\biggl\| \sum_{k\sim K_0 ,\,\ell\sim L_0} a_{k , \ell}\, {\text e}\left( x_1 \ell + x_2 k\ell + x_3 \omega (k , \ell)\right) 
\biggr\|_{L^{q_{\nu}}_{\#}\bigl[ |x_1|<1 , |x_2|<1 , |x_3|<{1\over\sqrt{K_0} L_0 \eta}\bigr]} \cr 
&\ \ll_{\nu ,\varepsilon}\, 
\left( 1+{NQ\over H^2}\right)^{3(\nu -1)\over 13\nu -12} \left( 1+ {H^3 Q\over N^2 R^2}\right)^{\nu-1\over 2(13\nu -12)} 
\left( 1+{H^{2\nu -6} Q^{\nu -2}\over R^{2\nu -6} N^{\nu -2}}\right)^{1\over 2(13\nu -12)} 
K_0^{\varepsilon +{1\over 2}} L_0^{{1\over 2}}\;.}
\eqno{(11.14)}
$$ 

In order to describe our use of (11.14) (and of the different bound (11.23), below) it is helpful to distinguish between certain cases. We shall find it convenient 
to consider two main cases:

\bigskip 

$\underline{\hbox{Case~I}\,}\,$, in which (11.9) holds and one has 
$$
24C_3 H>\sqrt{NR}\;;\eqno{(11.15)}
$$

$\underline{\hbox{Case~II}\,}\,$, in which (11.9) holds and (11.15) does not hold. 

\bigskip 

\noindent We shall also find it useful to split the latter case up into two (more specialized) cases: 

\bigskip 

$\underline{\hbox{Case~II(i)}}\,$, in which one has 
$$
1\leq {NR\over\left( 24 C_3 H\right)^2}<{Q\over R}\leq {B_2 H\over R}\;;\eqno{(11.16)}
$$

$\underline{\hbox{Case~II(ii)}}\,$, in which one has 
$$
1\leq {Q\over R}\leq\min\left\{ {NR\over\left( 24 C_3 H\right)^2} , {B_2 H\over R}\right\}\;.\eqno{(11.17)}
$$ 

\bigskip 

Note that Case~I, Case~II(i) and Case~II(ii) are mutually exclusive cases. Note also 
that if (11.9) holds then one of the above three conditions (i.e. one of (11.15), (11.16), or (11.17)) will be satisfied. 
Therefore we may complete the work of this section by obtaining,  in each one of the cases (i.e. in Case~I, in Case~II(i) and in Case~II(ii)), 
a sufficiently strong upper bound for some $A_p \bigl( {\cal Y}_1 , \text{\bf D} \bigr)$ 
with $p>2$ (in fact we shall always have $p\geq q_6=132/31>4$, since it is in that range that the bounds on  $A_p \bigl( {\cal Y}_1 , 
\text {\bf D} \bigr)$   are optimal, for our purposes). 

In Case~I we note that the conditions (11.9), (11.15) and assumptions (7.5), (11.12) and (11.13) imply: 
$$
1+{NQ\over H^2}\ll 1+{Q\over R}\ll {Q\over R}\;,\qquad 
1+ {H^3 Q\over N^2 R^2}\ll 1+{Q\over R}\quad {\text and}\quad 
1+{H^{2\nu -6} Q^{\nu -2}\over R^{2\nu -6} N^{\nu -2}}\ll 1+{Q^{\nu -2}\over R^{\nu -2}}\ll_{\nu}\,\left( {Q\over R}\right)^{\nu -2}\;.$$ 
Therefore we find that, by (11.2), (11.3), (11.14) and (11.11), one has
$$
\left( A_{q_{\nu}} \bigl( {\cal Y}_1 , \text{\bf D} \bigr) \right)^{1/q_{\nu}}\ll_{\nu ,\varepsilon}\, 
\left( {Q\over R}\right)^{\!\!{8\nu -9\over 2(13\nu -12)}} \!K_0^{\varepsilon +{1\over 2}} L_0^{{1\over 2}} 
=\left( {Q\over R}\right)^{\!\!{1 - 3 q_{\nu}^{-1}}} \!K_0^{\varepsilon +{1\over 2}} L_0^{{1\over 2}} \qquad\hbox{(in Case~I).}\eqno{(11.18)}
$$

With regard to  Case~II(i), we may note that (11.16) implies $H<N^{1/2}R^{1/2}$, and so (given (7.5)) it follows that (11.13) will hold for all $\nu \geq 6$. 
We also have (in Case~II(i) : 
$$
R< {N R^2\over H^2}\ll Q <H\;,
$$ 
and so: 
$$
1+{NQ\over H^2}< {NR\over H^2}+{NQ\over H^2}\ll {NQ\over H^2}\;,\qquad 1+ {H^3 Q\over N^2 R^2}<1+
\left( {H^2\over NR}\right)^{\!\!2} <2
$$
and 
$$
1+{H^{2\nu -6} Q^{\nu -2}\over R^{2\nu -6} N^{\nu -2}} = 1+\left( {H^2 Q\over R^2 N}\right)^{\nu -2} \left( {R^2\over H^2}\right) 
\ll_{\nu}\,\left( {H^2 Q\over R^2 N}\right)^{\nu -2}\ll \left( {H^2 Q\over R^2 N}\right)^{\nu -1}\;.
$$
Therefore, bearing in mind that (11.12) implies $1/(13\nu -12)<1/(13\nu -13)$,   we may deduce from (11.14) that in Case~II(i) one has 
$$
\eqalign{ &\biggl\| \sum_{k\sim K_0 ,\,\ell\sim L_0} a_{k , \ell}\, {\text e}\left( x_1 \ell + x_2 k\ell + x_3 \omega (k , \ell)\right) 
\biggr\|_{L^{q_{\nu}}_{\#}\bigl[ |x_1|<1 , |x_2|<1 , |x_3|<{1\over\sqrt{K_0} L_0 \eta}\bigr]} \cr 
&\ \ll_{\nu ,\varepsilon}\, 
\left( {NQ\over H^2}\right)^{\!\!{3\over 13}} \left( {H^2 Q\over R^2 N}\right)^{\!\!{1\over 26}} 
K_0^{\varepsilon +{1\over 2}} L_0^{{1\over 2}} 
=\left( {Q\over N R^2 /H^2}\right)^{\!\!{7\over 26}} \left( {NR\over H^2}\right)^{\!\!{6\over 13}} K_0^{\varepsilon +{1\over 2}} L_0^{{1\over 2}} }\eqno{(11.19)}
$$ 
for all integers $\nu\geq 6$. 

In order to simplify the application of (11.19) we observe now that, by (11.11) and (11.12), 
$$
q_{\nu}=\left( {13\over 3}\right)\left( 1 -\delta_{\nu}\right)\quad{\text with}\quad 
\delta_{\nu}={7\over 13(6\nu -5)}\in \biggl( 0 \,, {7\over 403}\biggr] \subset \biggl( 0 \,, {1\over 13}\biggr)\;.
$$ 
Therefore, given (11.1), a trivial bound on the relevant sums over $k$ and $\ell$ is enough to show that one has
$$
\eqalign{ &\biggl\| \sum_{k\sim K_0 ,\,\ell\sim L_0} a_{k , \ell}\, {\text e}\left( x_1 \ell + x_2 k\ell + x_3 \omega (k , \ell)\right) 
\biggr\|_{L^{{13\over 3}}_{\#}\bigl[ |x_1|<1 , |x_2|<1 , |x_3|<{1\over\sqrt{K_0} L_0 \eta}\bigr]}^{{13\over 3}} \cr 
&\ \ll \left( K_0 L_0\right)^{\left( {13\over 3}\right)\delta_{\nu}} 
\biggl\| \sum_{k\sim K_0 ,\,\ell\sim L_0} a_{k , \ell}\, {\text e}\left( x_1 \ell + x_2 k\ell + x_3 \omega (k , \ell)\right) 
\biggr\|_{L^{q_{\nu}}_{\#}\bigl[ |x_1|<1 , |x_2|<1 , |x_3|<{1\over\sqrt{K_0} L_0 \eta}\bigr]}^{q_{\nu}}\;.}
$$ 
By this and (11.19) we are able to conclude that, for $\varepsilon >0$ and all integers $\nu\geq 6$, one has: 
$$
\eqalignno{ &\biggl\| \sum_{k\sim K_0 ,\,\ell\sim L_0} a_{k , \ell}\, {\text e}\left( x_1 \ell + x_2 k\ell + x_3 \omega (k , \ell)\right) 
\biggr\|_{L^{{13\over 3}}_{\#}\bigl[ |x_1|<1 , |x_2|<1 , |x_3|<{1\over\sqrt{K_0} L_0 \eta}\bigr]} \cr 
&\ \ll_{\nu ,\varepsilon}\,\left( K_0 L_0\right)^{\delta_{\nu}} 
\left( \left( {Q\over N R^2 /H^2}\right)^{\!\!{7\over 26}} \left( {NR\over H^2}\right)^{\!\!{6\over 13}} K_0^{{\varepsilon +1\over 2}} L_0^{{1\over 2}}\right)^{\!\!{1-\delta_{\nu}}} \cr 
 &\ \ll \left( K_0 L_0\right)^{\left( {1\over 2}\right)\delta_{\nu}} 
\left( {Q\over N R^2 /H^2}\right)^{\!\!{7\over 26}} \left( {NR\over H^2}\right)^{\!\!{6\over 13}} K_0^{{\varepsilon +1\over 2}} L_0^{{1\over 2}}\;. &(11.20)}
$$ 

We now choose to put 
$\nu =6+[1/\varepsilon ]\,$ (where $[x]=\max\{j\in{\Bbb Z}\;:\; j\leq x\}$). This, given (8.7), (7.5) and the definition of $\delta_{\nu}$, is easily enough to ensure that one has 
$$
\left( K_0 L_0\right)^{\left( {1\over 2}\right)\delta_{\nu}}\ll K_0^{\delta_{\nu}}\leq K_0^{{\varepsilon\over 2}}\;.
$$ 
Therefore it follows by (11.2), (11.3) and the case $\nu =6+[1/\varepsilon ]$ of (11.20) that we have, for $\varepsilon >0$,  
$$
\left( A_{{13\over 3}} \bigl( {\cal Y}_1 , \text{\bf D} \bigr) \right)^{{3\over 13}}\ll_{\varepsilon}\, 
\left( {Q\over N R^2 /H^2}\right)^{\!\!{7\over 26}} \left( {NR\over H^2}\right)^{\!\!{6\over 13}} K_0^{\varepsilon +{1\over 2}} 
L_0^{{1\over 2}}\qquad\hbox{(in Case~II(i)).}\eqno{(11.21)}
$$ 

In Case~II(ii) it follows  by (8.7) and (11.17) that, in addition to (11.1) and (11.5)-(11.8), one has  
$$
{L_0^2\over K_0}\leq {\left( 24 C_3 L\right)^2\over K}={\left( 24 C_3\right)^2 H^2 Q\over R^2 N}\leq 1\;,\eqno{(11.22)}
$$ 
and so it follows from Proposition~10$'$ that we have, for $\varepsilon >0$, 
$$
\eqalign{ &\biggl\| \sum_{k\sim K_0 ,\,\ell\sim L_0} a_{k , \ell}\, {\text e}\left( x_1 \ell + x_2 k\ell + x_3 \omega (k , \ell)\right) 
\biggr\|_{L^{{48\over 11}}_{\#}\bigl[ |x_1|<1 , |x_2|<1 , |x_3|<{1\over\sqrt{K_0} L_0 \eta}\bigr]} \cr 
&\ \ll_{\varepsilon}\, 
\left( 1 + \eta K_0  L_0 \right)^{{1\over 48}} \left( 1 +{\eta K_0^2\over L_0}\right)^{\!\!{5\over 24}}  K_0^{\varepsilon +{1\over 2}} L_0^{{1\over 2}}}\;.
\eqno{(11.23)}
$$ 
By (11.4), (11.17), (11.22) and (8.7), we have here 
$$
1+\eta K_0 L_0\asymp 1+{Q^2\over R^2}\leq {2 Q^2 \over R^2}\quad {\text and}\quad 
1+{\eta K_0^2\over L_0}\asymp 1+{Q^2 K_0\over R^2 L_0^2}\leq {2 Q^2 K_0\over R^2 L_0^2}\asymp {Q^2 K\over R^2 L^2}={NQ\over H^2}\;.
$$ 
It therefore follows from (11.2), (11.3) and (11.23)  that, for $\varepsilon >0$, one has: 
$$
\left( A_{{48\over 11}} \left( {\cal Y}_1 , \text {\bf D}\right)\right)^{{11\over 48}} 
\ll_{\varepsilon} \left( {Q\over R}\right)^{\!\!{1\over 4}} \left( {NR\over H^2}\right)^{\!\!{5\over 24}}  
K_0^{\varepsilon +{1\over 2}} L_0^{{1\over 2}}\quad\hbox{(in Case~II(ii)).}\eqno{(11.24)}$$ 

\bigskip 
\beginsection
{12.  Results from the application of the Bombieri-Iwaniec method}  

Let $N$ is given either by (10.14), or else by (10.15). We seek a bound for the sum, over ${\cal C}$ and $Q$, one the right-hand side of (7.6). 
The bound (10.5) implies that one has 
$$
{\cal G}^{[A]}_Q\neq\emptyset\quad {\text only\ if}\quad A\ll {M R^2 \log N\over N Q^2}\;.
$$ 
Therefore the sum on the right-hand side of (7.6) has (given (7.5)) no more than $O( (\log M)^2 )$ terms. For this reason it will be enough that we obtain 
bounds for $S({\cal C}_Q)$ that are uniform, in the sense of being independent of the indices of summation, ${\cal C}$ and $Q$. 

Suppose now that $S({\cal C}_Q)$ is one of the terms occurring in the sum on the right-hand side of (7.6). It follows that we have 
$R\leq Q\leq\varepsilon H/2$ 
and either ${\cal C}={\cal B}\,$  (the set of `bad' intervals $I(a/q)$), or else ${\cal C}={\cal G}^{[A]}$ for some $A\geq\log N$. 
In the latter case Huxley's bounds  (10.5) and (10.16)-(10.17) imply: 
$$
\left( B\bigl( {\cal G}^{[A]}_Q ; V\bigr) V\right)^{\!{1\over p}} 
\left( \bigl| {\cal G}^{[A]}_Q \bigr| R/Q\right)^{\!{1-{2\over p}}} 
\ll \left( 1+{Q\over Q_2}\right)^{\!\!{1\over p}}\left( {M^2 R^4 (\log N)^{{1\over 4}} A^{{11\over 70}} 
\over H^2 N^4 V_0^{{2\over 3}}}\right)^{\!\!{1\over p}} 
\left( {M R^3 \log N\over A N Q^3}\right)^{\!\!{1-{2\over p}}} \;,
$$
when $p>2$. The right-hand side here is a decreasing function of $A$ for each fixed choice of $p$ satisfying $p>2+{11\over 70}$. 
Therefore, if $p>4\,$  (say), then we will certainly have 
$$
\eqalignno{\left( B\bigl( {\cal G}^{[A]}_Q ; V\bigr) V\right)^{\!{1\over p}} 
\left( \bigl| {\cal G}^{[A]}_Q \bigr| R/Q\right)^{\!{1-{2\over p}}} 
 &\ll \left( 1+{Q\over Q_2}\right)^{\!\!{1\over p}}\left( {M^2 R^4 (\log N)^{{57\over 140}} 
\over H^2 N^4 V_0^{{2\over 3}}}\right)^{\!\!{1\over p}} 
\left( {M R^3\over N Q^3}\right)^{\!\!{1-{2\over p}}} \cr 
 &=\left( 1+{Q\over Q_2}\right)^{\!\!{1\over p}}\left( {R^4 (\log N)^{{57\over 140}} 
\over H^2 N^2 V_0^{{2\over 3}}}\right)^{\!\!{1\over p}} 
\left( {R\over Q}\right)^{\!\!{3-{6\over p}}} M N^{-1}\;. &(12.1)}
$$

If we have instead ${\cal C}={\cal B}$, then we choose to observe that (10.4), (10.8) and (10.9) imply the bound 
$$
\bigl| {\cal B}_Q\bigr| \ll \left( {\eta M R^2\over N Q^2}\right)\left( 1 +{Q\over Q^*}\right) ,
$$ 
where 
$$
Q^*={\eta R^2\over Q'}\gg \left( {\Delta' R\over\Delta_2 H}\right)^{-{1\over 3}} R\gg 
\left( {H\over R}\right)^{\!\!{41\over 119}} R = Q_4\quad\hbox{(say).}\eqno{(12.2)}
$$
By this and (10.21)-(10.22), we obtain: 
$$
\eqalignno{ &\left( B\bigl( {\cal B}_Q ; V\bigr) V\right)^{{1\over p}} 
\left( \bigl| {\cal B}_Q \bigr| R/Q\right)^{1-{2\over p}} \cr 
 &\ \ll\left( 1+{Q\over Q_3}\right)^{\!\!{1\over p}}\left( {M^2 R^4 (\log N)^{{57\over 140}} 
\over H^2 N^4 (H/R)^{{2\over 3}}}\right)^{\!\!{1\over p}} 
\left( 1+{Q\over Q_4}\right)^{\!\!{1-{2\over p}}} \left( {\eta M R^3\over N Q^3}\right)^{\!\!{1-{2\over p}}} \cr 
 &\ =\left( 1+{Q\over Q_4}\right)^{\!\!{1-{2\over p}}} \left( 1+{Q\over Q_3}\right)^{\!\!{1\over p}} 
\left( {\eta^{p-2} R^4 (\log N)^{{57\over 140}} \over H^2 N^2 (H/R)^{{2\over 3}}}\right)^{\!\!{1\over p}} 
\left( {R\over Q}\right)^{\!\!{3-{6\over p}}} M N^{-1}\;. &(12.3)}
$$

By (10.12), (10.19), (12.2) and (7.5), we have 
$$
Q_2\ll Q_3 < Q_4\;.
$$ 
Therefore, subject to the condition that $p$ and $\eta$ satisfy 
$$
{\eta^{p-2}\over (H/R)^{{2\over 3}}}\ll {1\over V_0^{{2\over 3}}}\;,\eqno{(12.4)}
$$ 
we may conclude from (12.1) and (12.3) that one has 
$$
\left( B\bigl( {\cal C}_Q ; V\bigr) V\right)^{\!{1\over p}} 
\left( {\bigl| {\cal C}_Q \bigr| R\over Q}\right)^{\!\!{1-{2\over p}}} 
\ll\left( 1+{Q\over Q_2}\right)^{\!\!{1\over p}} 
\left( {R^4 (\log N)^{{57\over 140}}\over H^2 N^2 V_0^{{2\over 3}}}\right)^{\!\!{1\over p}} 
\!\left( {R\over Q}\right)^{\!\!{3-{6\over p}}}\left( {M\over N}\right)\quad{\text if}\ \,Q\leq Q_4 .\eqno{(12.5)}
$$ 
Assuming that $p>4$, it follows from (10.7) that the condition (12.4) will be satisfied if and only if one has $\eta\ll (H/R)^{-2/(51(p-2))}$. Therefore, given (7.5), we certainly obtain (12.5) subject 
to the conditions that 
$$
4<p<\infty\quad {\text and}\quad \eta\ll \left( {H\over R}\right)^{\!\!-{1\over 51}}\;.\eqno(12.6)
$$ 

If however, we have $Q\geq Q_4$ and ${\cal C}={\cal B}$, then we choose not to use the bound (10.4) for $|{\cal B}_Q|$, and instead simply recall (7.9). 
Since $Q_4 > Q_3$, this single change enables us to replace (12.3) with the alternative bound: 
$$
\left( B\bigl( {\cal B}_Q ; V\bigr) V\right)^{{1\over p}} 
\left( \bigl| {\cal B}_Q \bigr| R/Q\right)^{1-{2\over p}} 
\ll \left( {Q\over Q_3}\right)^{\!\!{1\over p}} 
\left( { R^4 (\log N)^{{57\over 140}}\over H^2 N^2 (H/R)^{{2\over 3}}}\right)^{\!\!{1\over p}} 
\left( {R\over Q}\right)^{\!\!{3-{6\over p}}} M N^{-1}\;,
$$
when $Q\geq Q_4$. This last upper bound exceeds that on the right-hand side of (12.5) by a factor $\Phi^{1/p}\,$ (say), where, by  (10.7), (10.12) and (10.19), $\Phi$ satisfies 
$$
\Phi < \left( {Q_2\over Q_3}\right) \left( {V_0\over H/R}\right)^{\!\!{2\over 3}}= 
\left( {H\over R}\right)^{\!\!{4\over 119}} \left( \log{H\over R}\right)^{\!\!{1\over 4}}\;.\eqno{(12.7)}
$$ 
By this, (12.1) and the trivial conditional inequality 
$$
1\leq (Q/ Q_4)^{\theta}\qquad\hbox{($Q\geq Q_4$ and $\theta\geq 0$),}
$$ 
we may conclude that if one chooses to put
$$
\theta_p=\cases p-4 & \text {if $4<p<13/3$}, \\ {1\over 2} &\text {if $p=13/3$}, \\  {7\over 11} &
\text {if $13/3<p<\infty $}\endcases\eqno{(12.8)}
$$ 
(for example), then the condition $Q\leq Q_4$ can be omitted from (12.5) if the bound appearing there is weakened through multiplication by 
$(1+(Q/Q_4)^{\theta_p}\Phi)^{1/p}$. That is, we have (given (12.7)): 
$$
\eqalign{&\left( B\bigl( {\cal C}_Q ; V\bigr) V\right)^{{1\over p}} 
\left( \bigl| {\cal C}_Q \bigr| R/Q\right)^{1-{2\over p}} \cr 
 &\quad\ll \left( 1+\left( {Q\over Q_4}\right)^{\!\!\theta_p} 
\left( {H\over R}\right)^{\!\!{4\over 119}} 
\left( \log {H\over R}\right)^{\!\!{1\over 4}}\right)^{\!\!{1\over p}} 
\left( 1+{Q\over Q_2}\right)^{\!\!{1\over p}} 
\left( {R^4 (\log N)^{{57\over 140}}\over H^2 N^2 V_0^{{2\over 3}}}\right)^{\!\!{1\over p}} 
\!\left( {R\over Q}\right)^{\!\!{3-{6\over p}}}\!\left( {M\over N}\right)\;,}
$$ 
subject to (12.6) and (12.8) holding. 

By the bound just obtained, in combination with (9.8), (10.1), (7.5) and (10.7), we find that
$$
\eqalignno{S\left( {\cal C}_Q\right) 
 &\ll \left( 1+\left( {Q\over Q_4}\right)^{\!\theta_p}\!\left( {H\over R}\right)^{\!\!{4\over 119}}\right)^{\!\!{1\over p}} 
\!\left( 1+{Q\over Q_2}\right)^{\!{1\over p}} \!\left( {R\over Q}\right)^{\!\!3- {6\over p}}\!\left( {(H/R)^2\over V_0^{2/3}}\right)^{{1\over p}} 
\!\left( {M R (\log N)^3\over N}\right)\left(   A_p\left( {\cal Y}_1 , \text{\bf D}\right)\right)^{{1\over p}}  \cr 
 &\ll_{\varepsilon}\,\left( 1+\left( {Q\over Q_4}\right)^{\!\theta_p}\!\left( {H\over R}\right)^{\!\!{4\over 119}}\right)^{\!\!{1\over p}} 
\left( 1+{Q\over Q_2}\right)^{\!{1\over p}}\!\left( {R\over Q}\right)^{\!\!3- {6\over p}}\!\left( {H\over R}\right)^{\!\!{22\over17 p}} 
\!\left( {M R\over N^{1-\varepsilon}}\right)\left(   A_p\left( {\cal Y}_1 , \text{\bf D}\right)\right)^{{1\over p}}\;, &(12.9)}
$$
subject to (12.6), (12.8) and the condition $\varepsilon >0$. 
We shall bound the factor $(  A_p ( {\cal Y}_1 , \text{\bf D} ))^{1/p}$ in (12.9) through an appeal to the results of the previous section. 
Our choice of $p$ therefore depends, in each individual case, on which one of the results (11.18),  (11.21), (11.24) is applied, and will (in all cases) 
satisfy 
$$
p\in\{ q_{\nu} \;:\; \nu\in{\Bbb Z}\ {\text and}\ \nu\geq 6\}\cup\{ \textstyle{13\over 3} , 
\textstyle{48\over 11} \}\subset {\Bbb Q}\cap [ \textstyle{132\over 31} , 
\textstyle{48\over 11} ]\subset (4.258 , 4.\dot 3\dot 6]\;.\eqno{(12.10)}
$$
It is 
helpful to note that, by (8.7) and (7.5), we have 
$$
K_0^{\varepsilon +{1\over 2}} L_0^{{1\over 2}}\ll_{\varepsilon} K^{\varepsilon +{1\over 2}} L^{{1\over 2}} 
\leq\left( {Q^2 NH\over R^4}\right)^{\!\!{1\over 2}} \left( {NH\over R^2}\right)^{\!\!\varepsilon} \leq 
{Q H^{{1\over 2}} N^{\varepsilon +{1\over 2}}\over R^2} \qquad\hbox{($R\leq Q\leq B_2 H$)}
$$ 
in (11.18), (11.21) and (11.24). Using this we deduce from (11.18), (11.21) and (11.24) three corresponding upper bounds for
$(  A_p ( {\cal Y}_1 , \text{\bf D} ))^{1/p\,} Q^{((\theta_p +7)/p) -3}$ 
that, by virtue of (12.8),  are each independent of $Q$. By these bounds, combined with (12.9), we obtain upper bounds for $S( {\cal C}_Q )$ that are, in each case, monotonic decreasing functions of $Q$. 
In particular, by (12.9), (11.18) and the lower bound on $Q$ in (11.9), we obtain (for $\varepsilon >0$): 
$$
S\left( {\cal C}_Q\right)\ll_{\varepsilon ,\nu}\,
\left( {H\over R}\right)^{\!\!\left( {22\over 17}\right) q_{\nu}^{-1}} 
\!\!\left( {M R\over N^{1-{\varepsilon \over 2}}}\right) 
\left( {H^{{1\over 2}} N^{{\varepsilon +1\over 2}}\over R}\right) 
=\left( {H\over R}\right)^{\!\!\left( {22\over 17}\right) q_{\nu}^{-1}} 
\!\!\left( {H\over N}\right)^{\!\!{1\over 2}} 
M N^{\varepsilon}\qquad 
\hbox{(in Case~I),}\eqno{(12.11)}
$$
subject to the final part of (12.6) holding (and with $\nu$, $q_{\nu}$ and `Case~I' being as described in Section~11). Note that the expression on the right-hand side of (12.11) involves fewer factors than that in (12.9). This is owing to  
the fact that, by (10.12), (11.11), (11.12),  (12.2), (12.8) and (7.5),  one has both 
$R/Q_2<1$ and 
$$
\left( {R\over Q_4}\right)^{\!\!\theta_p} \left( {H\over R}\right)^{\!\!{4\over 119}}
=\left( {H\over R}\right)^{\!\!{4-41 (p -4)\over 119}}
\leq \left( {H\over R}\right)^{\!\!{4\over 119}-\left( {41\over 119}\right)\left( q_6 -4\right)} 
=\left( {H\over R}\right)^{\!\!{4\over 119}-\left( {41\over 119}\right)\left( {8\over 31}\right)} 
=\left( {H\over R}\right)^{\!\!-{12\over 217}}<1\;,
$$
when $p=q_{\nu}$. 
Since $8/31<1/2<7/11$, it follows from (12.8) that one has an even stronger bound for \hfill\break
$( {R/ Q_4})^{\theta_p} ( H/ R)^{4/119}$ when $p\in\{ 13/3 , 48/11\}$. 

We postpone discussion of Case~II(i) until after dealing with Case~II(ii). In that latter case we may note that (11.17) implies $Q\geq R$. Therefore,  reasoning similar to that which produced (12.11)  
enables it to be deduced  from (12.9) and (11.24) that one has, for $\varepsilon >0$,  
$$
S\left( {\cal C}_Q\right) 
\ll_{\varepsilon}\, \left( {NR\over H^2}\right)^{\!\!{5\over 24}} 
\left( {H\over R}\right)^{\!\!\left( {22\over17 }\right)\left( {11\over 48}\right) } 
\left( {H\over N}\right)^{\!\!{1\over 2}} M N^{\varepsilon}
=\left( {H\over R}\right)^{\!\!{3\over 34}} \left( {H\over N}\right)^{\!\!{7\over 24}} M 
N^{\varepsilon}\qquad\hbox{(in Case~II(ii)),}\eqno{(12.12)}
$$
subject to the final part of (12.6) holding. 

In Case~II(i) of Section~11 we apply (12.9) in combination with (11.21). 
The conditions (11.17) defining Case~II(i) certainly imply $Q>R$, and so, for reasons given in the paragraph containing  (12.11), we find that 
$$
\eqalign{ S\left( {\cal C}_Q\right) 
 &\ll_{\varepsilon}\,\left( {H\over R}\right)^{\!\!\left( {22\over17}\right)\left( {3\over 13}\right)} 
\!\left( {M R\over N^{1-{\varepsilon \over 2}}}\right)
\left( {R\over N R^2 /H^2}\right)^{\!\!{7\over 26}} \left( {NR\over H^2}\right)^{\!\!{6\over 13}}
\left( {H^{{1\over 2}} N^{{\varepsilon +1\over 2}}\over R}\right) \cr 
 &=\left( {NR\over H^2}\right)^{\!\!{5\over 26}}
\!\left( {H\over R}\right)^{\!\!{66\over221}} 
\left( {H\over N}\right)^{\!\!{1\over 2}} M N^{\varepsilon} \cr 
&=\left( {H\over R}\right)^{\!\!{47\over442}} \left( {H\over N}\right)^{\!\!{4\over 13}} M N^{\varepsilon}\;,}
$$
in Case~II(i). 

We observe that this last upper bound for $S( {\cal C}_Q)$ will exceed the corresponding Case~II(ii) upper bound (given in (12.12)) only when one has 
$(H/R)^{4/221}(H/N)^{5/312}>1$, and so only when 
$$
{N\over H}<\left( {H\over R}\right)^{\!\!{96\over 85}}\;.\eqno{(12.13)}
$$ 
Subject to (12.13) holding, one has, by (10.12), 
$$
{N R^2 /H^2\over Q_2}=\left( {N\over H}\right) \left( {R\over H}\right)^{\!\!{158\over 119}}\left(\log 
{H\over R}\right)^{\!\!{3\over 4}} 
< \left( {H\over R}\right)^{\!\!-{118\over 595}}\left(\log {H\over R}\right)^{\!\!{3\over 4}}\ll 1
$$ 
and, by (12.2) and (12.8), 
$$
\left( {N R^2 /H^2\over Q_4}\right)^{\!\!\theta_{13/3}}\!\left( {H\over R}\right)^{\!\!{4\over 119}} 
=\left( \left( {N\over H}\right) \left( {R\over H}\right)^{\!\!{160\over 119}}\right)^{\!\!{1\over 2}} \left( {H\over R}\right)^{\!\!{4\over 119}} 
< \left( {H\over R}\right)^{\!\!{48\over 85}-{76\over 119}}=\left( {H\over R}\right)^{\!\!-{44\over 595}}<1\;.
$$
It therefore follows from (12.9), (11.21) and the lower bound estimate $Q\gg N R^2 /H^2$ (implied by (11.16)) that, when the conditions for Case~II(i) are satisfied and (12.13) holds, one will have (for $\varepsilon >0$): 
$$
\eqalign{ S\left( {\cal C}_Q\right) 
  &\ll_{\varepsilon}\,\left( {R\over N R^2 /H^2}\right)^{\!\!2- {18\over 13}}
\!\left( {H\over R}\right)^{\!\!{66\over 221}} 
\!\left( {M R\over N^{1-{\varepsilon \over 2}}}\right)
\left( {NR\over H^2}\right)^{\!\!{6\over 13}}
\left( {H^{{1\over 2}} N^{{\varepsilon +1\over 2}}\over R}\right) \cr 
 &=\left( {H^2\over N R}\right)^{\!\!{2\over 13}}\!\left( {H\over R}\right)^{\!\!{66\over 221}} 
\left( {H\over N}\right)^{\!\!{1\over 2}} M N^{\varepsilon}\cr 
 &=\left( {H\over R}\right)^{\!\!{100\over 221}} 
\left( {H\over N}\right)^{\!\!{17\over 26}} M N^{\varepsilon}\;,}
$$ 
subject to the final condition of (12.6) being satisfied. 

By this last finding, allied with both (12.12) and the observation made at the start of the preceding paragraph, one has (for $\varepsilon >0$): 
$$
S\left( {\cal C}_Q\right) 
 \ll_{\varepsilon}\,\left( \left( {H\over R}\right)^{\!\!{3\over 34}} \left( {H\over N}\right)^{\!\!{7\over 24}} 
+\left( {H\over R}\right)^{\!\!{100\over 221}} 
\left( {H\over N}\right)^{\!\!{17\over 26}}\right) M N^{\varepsilon}\qquad\hbox{(in Case~II),}\eqno{(12.14)}
$$
subject to the final condition of (12.6) being satisfied. 

Given the definitions (in Section~11) of Cases~I and ~II, and bearing in mind the point noted in the first paragraph of this section, it follows from (7.5), (7.6), (12.11) and (12.14) that we have, for each 
$\varepsilon >0$, a bound of the form 
$$
|S|\leq\Psi_{\varepsilon , \nu}  ( H/R\,, H/N) M^{1+\varepsilon}\;,\eqno{(12.15)}
$$ 
where  
$$
\Psi_{\varepsilon , \nu} (\Delta , \delta) =
\cases  C_{II}(\varepsilon)\left(\Delta^{3\over 34} \delta^{7\over 24} + \Delta^{{100\over 221}} \delta^{{17\over 26}}\right) 
&\text {if $\,576 C_3^2 \Delta\delta\leq 1$}, \\ 
C_{I}(\varepsilon , \nu) \Delta^{\left( {22\over 17}\right) q_{\nu}^{-1}} \delta^{{1\over 2}}  &\text{otherwise},
\endcases\eqno{(12.16)}
$$ 
with $C_{II}(\varepsilon)$ denoting a positive constant constructed from 
$C_2 , \ldots , C_6$ and $\varepsilon$, while $C_I(\varepsilon , \nu)$ denotes a positive constant 
constructed from $C_2 , \ldots , C_6$, $\varepsilon$ and $\nu$. 
This, of course, assumes that the values of $T$, $H$ and $M$, and that of our chosen of parameter $N$ (an integer satisfying either (10.14), or else  (10.15)),  
are consistent with being able to satisfy all of the conditions (7.3), (7.4), (7.5), (10.3), (10.6), (10.8)-(10.11), (11.12), (11.13) and (12.6), as well as 
the condition attached to the relevant case, `(A)' or `(B)', of Lemma~10.1. 
We devote the remainder of the section to obtaining a resolution of this issue that 
will complete our proof of Theorem~2. Therefore we shall no longer be assuming that all of 
the conditions just mentioned are satisfied (for our goal, in what follows, is a concise description of   
the circumstances in which certain, quite specific, choices of $N\in{\Bbb N}$ and $R\in{\Bbb N}$ 
will satisfy all those conditions). We shall, however, find it convenient to assume that the 
positive constants $C_I (\varepsilon , \nu)$, $C_{II}(\varepsilon)$ that occur in (12.16)  
satisfy 
$$
C_I (\varepsilon , \nu)/C_{II}(\varepsilon)\geq 24 C_3 \;.\eqno{(12.17)}
$$ 
This causes no loss of generality, for we are effectively able to ensure that 
(12.17) will hold by substituting $\max\{ C_I (\varepsilon , \nu) , 24 C_3 C_{II}(\varepsilon)\}$ for $C_I (\varepsilon , \nu)$ in (12.16).

As a first step, we specify the parameters $\eta$, $Q'$ and $\Delta'$ (occurring in (10.3), (10.8)-(10.11) and (12.6)) by putting: 
$$\eta =\left( {R\over H}\right)^{\!\!{1\over 51}}\;,\qquad 
\Delta' =\left( {R\over H}\right)^{\!\!{4\over 119}} \Delta_2 \quad\ {\text and}\quad\  
Q' =\left( {\Delta ' R\over \Delta_2 H}\right)^{\!\!{1\over 3}} R\eta 
=\left( {R\over H}\right)^{\!\!{41\over 119}} R\eta = \left( {R\over H}\right)^{\!\!{130\over 357}} R\;.\eqno(12.18)$$ 
Assuming that we have $0<R/H<1\,$ (as (7.5) would imply), these specifications, along with that of $V_0$ in (10.7), 
can be shown to ensure that the conditions (10.8)-(10.11) and the condition on $\eta$ in (12.6) are satisfied, and that 
$Q'$ and $\eta$ satisfy $Q'<\eta R < R\,$  (as stated in (10.3)). 
We are therefore able to reduce the set of conditions on $\eta$, $Q'$ and $\Delta'$ to the combination 
of (12.18) and the single condition $Q'\geq 1+M^2 /T\,$ (seen in (10.3)); given (12.18), this single condition 
(on $Q'$) will be satisfied if one has 
$$
R \left( {R\over H}\right)^{{130\over 357}}\geq 1+{M^2\over T}\;.\eqno{(12.19)}
$$ 

\medskip
\noindent
{\bf Lemma~12.1.}
{\sl  Let $E_1,E_2\in [1/16 , 16]$. Let $\widetilde N , \widetilde R \in (0 , \infty)$ 
and $\nu\in [6 , \infty)$ be such that one has both 
$$
E_1 \widetilde N \widetilde R^2 = {C_3 M^3 \over 2 T}\eqno{(12.20)}
$$ 
and  
$$
E_2 B_1 \widetilde N^{{\nu-2\over 2\nu -6}} \widetilde R^{{\nu-4\over 2\nu -6}} \geq H\;.\eqno{(12.21)}
$$ 
Let $\widetilde V_0 , \widetilde V_1 , \widetilde V_2\in (0,\infty)$ be determined by the constraint 
that the equalities in (10.7) should hold if $\widetilde V_0$, 
$\widetilde V_1$, $\widetilde V_2$, $\widetilde N$ 
and $\widetilde R$ are substituted for $V_0$, $V_1$, $V_2$, $N$ and $R$, respectively. 
Suppose that one has either 
$$
E_2\min\bigl\{ \widetilde V_1 \,,\, \widetilde V_2\bigr\} \geq \widetilde V_0\eqno{(12.22)}
$$ 
or 
$$
M\leq C_6 T^{{1\over 2}}\quad\ {\text and}\quad\ \widetilde N^3 \leq {E_2 B_3^3 M^6 \over T^2 H}\;. 
\eqno{(12.23)}
$$ 
Then the conditions (7.4), (7.5), (10.6) and (12.19) will hold if 
$\widetilde V_0$, 
$\widetilde V_1$, $\widetilde V_2$, $\widetilde N$ 
and $\widetilde R$ are substituted for $V_0$, $V_1$, $V_2$, $N$ and $R$, respectively.}

\noindent
{\bf Proof.}\ 
Let (7.4)$'$, (7.5)$'$, (10.6)$'$ and (12.19)$'$ denote the conditions that (7.4), (7.5), (10.6) and (12.19) (respectively) become when $\widetilde N$, $\widetilde R$, $\widetilde V_0$, $\widetilde V_1$ and $\widetilde V_2$ are substituted for $N$, $R$, $V_0$, $V_1$ and $V_2\,$ (respectively). We are required to show that it follows from the hypotheses of the lemma that the conditions (7.4)$'$, (7.5)$'$, (10.6)$'$ and (12.19)$'$ 
are satisfied. 

Since $E_1^{-1}\leq 16$ and $E_2\leq 16$, it follows from (7.1), (12.20) and (12.21) that we have 
$$
\widetilde N \widetilde R^2 \leq 8 C_3 M^3 T^{-1} < 2^{-12} H^3 \leq B_1^3 
\widetilde N^{{3\nu-6\over 2\nu -6}} \widetilde R^{{3\nu-12\over 2\nu -6}} \;,\eqno{(12.24)}
$$ 
and so 
$$
\left( {\widetilde R\over \widetilde N}\right)^{\!\!2} < B_1^3 \left( {\widetilde R\over \widetilde N}\right)^{\!\!{3\nu -12 \over 2\nu -6}} \quad\ 
{\text and}\quad\ {H\over \widetilde N}\leq 16 B_1 \left( {\widetilde R\over \widetilde N}\right)^{\!\!{\nu -4 \over 2\nu -6}}\;. 
$$
Since $\nu \geq 6$, and since we may assume here that $0< B_1\leq 1/16$, 
the last two inequalities above imply that one has both 
$$
{\widetilde R\over \widetilde N}<B_1^{{6(\nu -3)\over \nu}}\leq B_1^3\leq 2^{-12}\quad\ {\text and}\quad\ 
{H\over \widetilde N}\leq 16 B_1^{{4(\nu -3)\over \nu}}\leq 16 B_1^2 \leq B_1\;.\eqno{(12.25)}
$$ 
Note that the first two inequalities of (12.24) imply that one has $H^3 > 2^{12} \widetilde N \widetilde R^2$, and so $H/\widetilde N > 2^{12} \widetilde R^2/H^2$. 
By this and the final two inequalities of (12.25), one can deduce that 
$$
{\widetilde R\over H}<{B_1\over 16}\leq 2^{-8}\;.\eqno{(12.26)}
$$ 

Next we observe that (12.22) and (12.26) would imply: 
$$
\min\bigl\{ \widetilde V_1 \,,\, \widetilde V_2\bigr\}\geq E_2^{-1} \widetilde V_0 = E_2^{-1} 
\left( {H\over \widetilde R}\right)^{\!\!{18\over 17}} 
>16^{-1} \left( {16\over B_1}\right)^{\!\!{18\over 17}}>B_1^{-{18\over 17}}\geq B_1^{-1}\geq 16\;.
$$
If (12.22) does not hold, then (by hypothesis) we have instead the inequalities in (12.23), and can combine these 
with (12.20) so as to obtain: 
$$
\widetilde V_1 ={\widetilde R^4 \over H\widetilde N}={\left( E_1 \widetilde N \widetilde R^2 \right)^{\!2} 
\over E_1^2 H \widetilde N^3} 
\geq {\left( E_1 \widetilde N \widetilde R^2 \right)^{\!2} \over E_1^2 E_2 B_3^3 M^6 T^{-2} } 
={C_3^2 \over 4 E_1^2 E_2 B_3^3}\geq 2^{-14} B_3^{-3}
$$
and 
$$
\widetilde V_2 = {M^2 \over H \widetilde N^3} \geq {T^2 \over E_2 B_3^3 M^4}\geq {1\over 16 C_6^4 B_3^3}\;.
$$ 
Since we may assume (for example) that $0<B_3\leq (32 C_6)^{-4/3}\leq 2^{-20/3}$, it follows from (12.26) and the points just noted that, regardless 
of whether it is (12.22) or (12.23) that holds, we are certain to have 
$$
\widetilde V_j\geq 16\qquad\hbox{($j=0,1,2$).}\eqno{(12.27)}
$$ 
The inequality $\widetilde V_1\geq 16$ implies $H\widetilde N\leq \widetilde R^4 /16$. This, together with (12.25) gives: 
$$
{H\over \widetilde R^2}\leq {4 B_1 \bigl( H \widetilde N \bigr)^{\!{1\over 2}} \over \widetilde R^2} \leq B_1\;.\eqno{(12.28)}
$$
The inequality $\widetilde V_2\geq 16$ implies $H \widetilde N^3\leq M^2 /16$, and so $H \widetilde N^2 < M^2 /\widetilde N$. This, together with 
(12.25), (12.26), (12.28) and the hypothesis $H\geq 1$, gives 
$\widetilde N^2\leq H \widetilde N^2\leq M^2 /\widetilde N = M^2 (\widetilde R/\widetilde N)(\widetilde R/H)H/\widetilde R^2 < M^2 B_1^5$, so that one has 
$$
{\widetilde N\over M}< B_1^{{5\over 2}}\leq B_1\;.\eqno{(12.29)}
$$ 
By (12.25), (12.26), (12.28) and (12.29), we conclude that the condition (7.5)$'$ is satisfied. 
\par 
By (12.27), one has  
$$
1>{1\over \bigl( \widetilde V_1 \widetilde V_2\bigr)^{\!{1\over 2}}} ={H \widetilde N^2 \over M \widetilde R^2}\;,
$$ 
which implies (7.4)$'$. With regard to (10.6), we note that (12.27) and (12.28) imply: 
$$
{\widetilde V_0^2 \over \widetilde N \bigl( \widetilde V_1 \widetilde V_2\bigr)^{\!{1\over 2}}}<{\widetilde V_0^2 \over \widetilde N}<{\widetilde V_0^{{85\over 36}} \over \widetilde N} 
={H^{{5\over 2}}\over \widetilde N \widetilde R^{{5\over 2}}} 
\leq B_1 \left( {H^3\over \widetilde N^2 \widetilde R}\right)^{\!\!{1\over 2}}\;.
$$
Since it moreover follows from (12.21), the hypothesis $\nu\geq 6$ and the first three inequalities of (12.25) 
that one has $H\leq E_2 B_1 \widetilde N^{2/3} \widetilde R^{1/3}\leq 16 B_1 (\widetilde N^2 \widetilde R)^{1/3}$, we may deduce that 
$$
{\widetilde V_0^2 \over \widetilde N \bigl( \widetilde V_1 \widetilde V_2\bigr)^{\!{1\over 2}}} 
< B_1 \left( \left( 16 B_1\right)^3\right)^{\!{1\over 2}} 
<\left( 16 B_1\right)^{{5\over 2}}\leq 1\;,
$$ 
so that (10.6)$'$ is satisfied. 

Finally, with regard to the condition (12.19)$'$, we note that (12.20) implies  
$$
\left( {\widetilde V_1 \over \widetilde V_2}\right)^{\!\!{1\over 2}} 
={\widetilde R^2 \widetilde N\over M} ={C_3 M^2 \over 2 E_1 T}\geq {M^2 \over 32 T}\;,
$$ 
so that, by (12.27), one has: 
$$
1+{M^2\over T}\leq 1 + 32 \left( {\widetilde V_1 \over \widetilde V_2}\right)^{\!\!{1\over 2}} 
\leq 9 \widetilde V_1^{{1\over 2}}\;.
$$ 
This proof may therefore be completed by observing that   
(12.26) and the first three inequalities of (12.25) imply that one has 
$$
{9 \widetilde V_1^{{1\over 2}} \over \widetilde R (\widetilde R/H)^{{130\over 357}}} 
={9\widetilde R (H/\widetilde R)^{{130\over 357}}\over (H \widetilde N)^{{1\over 2}}}   
=9\left( {\widetilde R\over H}\right)^{\!\!{97\over 714}} 
\left( {\widetilde R\over \widetilde N}\right)^{\!\!{1\over 2}} <1. \quad\blacksquare 
$$ 

\medskip
\noindent
{\bf Corollary~12.1.1.} Suppose that the hypotheses of Lemma~12.1 are satisfied, and that  
$E_1 = E_2 =1$.  Put  
$$
N=\left\lceil \widetilde N \right\rceil\eqno{(12.30)}
$$ 
and 
$$
R=\left\lceil \left( {C_3 M^3 \over 2 N T}\right)^{{1\over 2}} \right\rceil\eqno{(12.31)}
$$ 
(where $\lceil x\rceil =\min\{ n\in{\Bbb Z} \;:\; n\geq x \}$, the `ceiling' function). 
Then the conditions (7.3), (7.4), (7.5), (10.6)-(10.7), (11.13) and (12.19) are satisfied, and one has:  
$$
\widetilde N\leq N < 2\widetilde N\eqno{(12.32)}
$$
and 
$$
{\widetilde R \over\sqrt{2}} < R < 2\widetilde R\;.\eqno{(12.33)}
$$

\noindent{\bf Proof.}\quad 
Assume (12.30) and (12.31). 
Then, by Lemma~12.1, 
$$
\widetilde N +1 > N\geq \widetilde N\geq B_1^{-1} H \geq B_1^{-1}\geq 1\;,
$$ 
and so we obtain (12.32). 
By (12.20) and Lemma~12.1, we have also: 
$$
{C_3 M^3 \over 2 \widetilde N T} =\widetilde R^2 \geq B_1^{-1} H \geq B_1^{-1} \geq 16\;.
$$ 
By this, (12.31) and (12.32), we have 
$$
1<\left( {C_3 M^3 \over 2 N T}\right)^{\!\!{1\over 2}} \leq R 
< \left( {C_3 M^3 \over 2 N T}\right)^{\!\!{1\over 2}} +1 
< 2\left( {C_3 M^3 \over 2 N T}\right)^{!\!{1\over 2}}\;,\eqno{(12.34)}
$$
and so (with the help of (12.32) and (12.20)) the result (12.33) follows. 

By (12.34), we have 
$E N R^2 = C_3 M^3 /(2 T)$ for some $E\in [1,4)$. 
Given that $\nu \geq 6$, it moreover follows from (12.32), (12.33) and the 
hypotheses of Lemma~12.1 stated between (12.21) and (12.23) that one has 
$H\leq 2^{1/4} B_1 N^{(\nu -2)/(2\nu -6)} R^{(\nu -4)/(2\nu -6)}$ 
and either 
$M\leq C_6 T^{1/2}$ and $N^3 <8B_3^3 M^6 /( T^2 H)$, or else 
$V_0 < 16\min\{ V_1 , V_2\}$. Therefore we have (11.13), 
and it follows by the case $(E_1 , E_2) = (E , 16)$ of 
Lemma~12.1 that the conditions (7.4), (7.5), (10.6)-(10.7) and (12.19) are satisfied. 
By (12.34), we have also the inequalities stated in (7.3).\quad$\blacksquare$ 

\medskip
\noindent
{\bf Lemma~12.2.} 
{\sl Let $E'\in [\sqrt{2} , \infty)$. Let the hypotheses of Theorem~2, up to and including (6.3), be satisfied. 
Let the function $\Psi_{\varepsilon , \nu} (\Delta , \delta)$ 
be as described in (12.16), and let (10.14)$'$ and (10.15)$'$ denote the conditions that 
(10.14) and (10.15) (respectively) become when, in both of those two conditions, 
one substitutes $\widetilde N$ and the relation of equality 
for $N$ and the relation `$\asymp$' (respectively). 
Suppose moreover that the hypotheses of Lemma~12.1 are satisfied, and that  $E_1 = E_2 =1$. 
Then one has 
$$
|S|\leq \Psi_{\varepsilon , \nu} \bigl( E' H / \widetilde R , H / \widetilde N \bigr) M^{1+\varepsilon}\eqno{(12.35)}
$$
in each of the following two independent cases: \hfill\break 
\vskip -2mm 
\indent (A)$'$\ when the conditions (12.22) and (10.14)$'$ are satisfied;\hfill\break 
\indent (B)$'$\ when $M\leq C_6 T^{1/2}$ and the condition (10.15)$'$ is satisfied.
}

\noindent
{\bf Proof.}\ 
Let $N$ and $R$ be given by (12.30) and (12.31). Since $\widetilde N >0$, it follows that we have 
$N,R\in{\Bbb N}$. By Corollary~12.1.1, it is moreover the case that the conditions  (7.3), (7.4), (7.5), 
(10.6)-(10.7), (11.13) and (12.19) are satisfied, and that the inequalities in (12.32) and (12.33) hold. 
By hypothesis, we have (11.12). By choosing $\eta$, $Q'$ and $\Delta'$ to be as stated in (12.18), 
we are able to ensure that the conditions  (10.3), (10.8)-(10.11) and (12.6) are satisfied 
(regarding this point see the remarks preceding (12.10), and those preceding (12.19)). 
Therefore we obtain the result stated in (12.15)-(12.16) if both (10.14) and the conditions 
attached to Case~(A) of Lemma~10.1 are satisfied, or if both (10.15) and the conditions 
attached to Case~(B) of Lemma~10.1 are satisfied. 
This occurs in Case~(A)$'$: for (10.14)$'$ and (12.32) imply (10.14), while (12.22), (12.32) and (12.33) imply the 
bound $V_0\ll\min\{ V_1 , V_2\}$. It also occurs in Case~(B)$'$, for (10.15)$'$ and (12.32) imply (10.15), while 
the inequality $M\leq C_6 T^{1/2}$ implies $M^2 \ll T$. 
Therefore we may complete this proof by showing that one has 
$\Psi (E' H / \widetilde R , H / \widetilde N)\geq \Psi (H/R , H/N)$. 

By (12.32), (12.33) and our hypothesis concerning $E'$, we have 
$H/\widetilde N \geq H/N$ and $(E' H / \widetilde R)/(H/R)\geq E' /\sqrt{2}\geq 1 $. 
By Lemma~12.1 we have, moreover, $H/\widetilde N\leq B_1\leq 1$, and by (7.5) we have 
$H/R\geq B_1^{-1}\geq 1$. 
Given the points just noted and the definition of $\Psi(\Delta , \delta)$ in (12.16),  
we find that the desired inequality $\Psi (E' H / \widetilde R , H / \widetilde N)\geq \Psi (H/R , H/N)$ 
is a consequence of the 
observations that, for $\sigma =\pm 1$, one has both 
$(\partial /\partial \Delta)\Psi(\Delta , \delta) >0$ and $(\partial /\partial \delta)\Psi(\Delta , \delta) >0$
on the subset ${\cal R}(\sigma)=\{ (\Delta , \delta)\in [1,\infty)\times(0,1] \;:\; 576 C_3^2 \Delta \delta\sigma <\sigma\}$ of the $\Delta,\delta$-plane, and that at all points $(\Delta ,\delta)$ 
lying on the part of the hyperbola $576 C_3^2 \Delta \delta =1$  that is the common boundary 
of ${\cal R}(1)$ and ${\cal R}(-1)$ one has  both 
$$
\Delta^{{3\over 34}} \delta^{{7\over 24}} + \Delta^{{100\over 221}} \delta^{{17\over 26}}  
= \Delta^{-{83\over 408}} \left( 24 C_3\right)^{-{7\over 12}} + \Delta^{-{89\over 442}} \left( 24 C_3\right)^{-{17\over 13}} 
< \Delta^{-{89\over 442}}
$$
and
$$
\Delta^{\left( {22\over 17}\right) q_{\nu}^{-1}} \delta^{{1\over 2}} 
=\Delta^{{11(6\nu -5)\over 17(13\nu -12)}-{1\over 2}} \left( 24 C_3\right)^{-1} 
\geq \Delta^{{66\over 221}-{1\over 2}} \left( 24 C_3\right)^{-1}
\geq \Delta^{-{89\over 442}} C_{II}(\varepsilon) / C_I (\varepsilon , \nu)
$$
(with the latter following by virtue of (11.11), (11.12) and our assumption in (12.17)).\quad$\blacksquare$ 

\bigskip 

\noindent{\bf Completion of the proof of Theorem~2.}\ 
A number of essentially very straightforward calculations suffice show that Theorem~2 is a corollary of Lemma~12.2. These calculations are not of much 
interest in themselves, so we shall give a sketchy account of them that 
covers the key points, but omits much of the (purely computational) detail. 

We observe firstly that, when $\widetilde N$ and $\widetilde R$ are given by 
$$
\widetilde N = H\left( {M\over H}\right)^{\!\!{41\over 25}} T^{-{49\over 100}} 
(\log T)^{{969\over 14000}}\eqno{(12.36)}
$$ 
and 
$$
\widetilde R = \left( {C_3 M^3 \over 2 \widetilde N T}\right)^{\!\!{1\over 2}} \;,\eqno{(12.37)}
$$ 
the hypotheses of Case~(A)$'$ of Lemma~12.2 will be satisfied if one has both 
$$
H\leq B_1 \widetilde N^{{\nu -2\over 2\nu -6}} \widetilde R^{{\nu -4\over 2\nu -6}}\eqno{(12.38)}
$$ 
and 
$$
\min\left\{ {\widetilde R^4\over H\widetilde N} \,,\, {M^2\over H \widetilde N^3}\right\} \geq 
\left( {H\over \widetilde R}\right)^{\!\!{18\over 17}}\;.\eqno{(12.39)}
$$ 
Moreover, subject to (12.37) holding, the conditions (12.38) and (12.39) are satisfied 
if and only if one has both 
$$
H\leq B_1 \left( {C_3 M^3\over 2 T}\right)^{\!\!{\nu -4\over 4 (\nu -3)}} \widetilde N^{{\nu\over 4(\nu -3)}} 
\eqno{(12.40)}
$$ 
and 
$$
\min\left\{ \left( {C_3 M^3 \over 2 T}\right)^{\!\!{43\over 60}} \,,\, 
\left( {C_3 M^3 \over 2 T}\right)^{\!\!{3\over 20}} M^{{17\over 30}}\right\} \geq H^{{7\over 12}} \widetilde N\;.\eqno{(12.41)}
$$ 
A calculation shows that, when $C_3 =2$ and $\widetilde N$ is given by (12.36), the conditions (12.40) and (12.41) 
become conditions on $T$, $M$ and $H$ that are effectively equivalent to the conditions (6.4)-(6.6) 
of Theorem~2 (albeit with $B_1^{25(\nu -3)/(29\nu -75)}$ in place of the constant $B_5$): 
note in particular that, although neither (6.4) nor (6.5) applies when $M$ lies in the interval 
${\cal E}=(T^{7/16} (\log T)^{57/448} , T^{9/16} (\log T)^{-57/448})$, it does nevertheless 
follow directly from (7.1) that we have 
$H\geq  (\log T)^{171/140}\max\left\{ T^4 M^{-9} \,,\, T^{-6} M^{11}\right\}$  
when $M$ lies in the interval $(T^{13/30} (\log T)^{171/1400} , T^{17/30} (\log T)^{-171/1400})\supset {\cal E}$. Therefore, given that we have $B_1\in (0,1]$, $C_3\geq 2$ and $\nu\geq 6$, it follows that 
the  conditions (6.4)-(6.6), with $B_5 = B_1^{25/29}\,$ (say), are sufficient to ensure that if one 
puts $E'=16/9$ and chooses $\widetilde N$ and $\widetilde R$ to be as stated in (12.36) and (12.37) then the 
hypotheses of Case~(A)$'$ of Lemma~12.2 will be satisfied. Consequently it follows from Lemma~12.2 that, if $B_5\leq B_1$ and the conditions in (6.4)-(6.6) are satisfied, then one has 
$$
{|S|\over H M^{\varepsilon}}\leq 
\cases
\left( {16\over 9}\right)^{{100\over 221}} C_{II}(\varepsilon) 
\left(\left(\displaystyle{M\over H}\right)^{\!\!{277\over 600}} 
T^{397\over 2400} + \left( \displaystyle{H\over M}\right)^{\!\!{19\over 50}} T^{1133\over 2600}\right) 
 & \text {if $\,H\leq \displaystyle{(\log T)^{{969\over 64960}} M 
\over (128 C_3)^{{75\over 116}} T^{{149\over 464}}}$,} \cr 
\left( {16\over 9}\right)^{{31\over 102}} C_{I}(\varepsilon , \nu) 
\left(\displaystyle{H\over M}\right)^{\!\left( {22\over 25}\right) q_{\nu}^{-1}-{9\over 50}} 
T^{\left( {33\over 100}\right) q_{\nu}^{-1}+{49\over 200}} &\text {otherwise.}\endcases
\eqno{(12.42)}
$$ 

We observe also that, in order for the hypotheses of Case~(B)$'$ of Lemma~12.2 to be satisfied, 
it is enough that one have (12.37), (12.38), (6.10) and 
$$
\widetilde N = \min\left\{ {M^{{7\over 8}} (\log T)^{{969\over 5600}} 
\over T^{{3\over 20}} H^{{29\over 40}}} \;,\, {B_3 M^2\over T^{{2\over 3}} H^{{1\over 3}}}\right\}
\eqno{(12.43)}
$$ 
(note, in particular, that (6.10) and (12.43) imply that the case $E_2 =1$ of the condition (12.23) is satisfied). Subject to $\widetilde R$ being given by (12.37), the condition (12.38) becomes equivalent 
to the inequality in (12.40), and so (given that $C_3 \geq 2$ and $\nu\geq 6$) we may deduce that, when  
$\widetilde N$ and $\widetilde R$ are as stated in (12.37) and (12.43),  the condition (12.38) will hold if one has
$$
\left( {H\over B_1}\right)^{\!\!{4(\nu -3)\over \nu}} 
\leq \left( {M^3 \over T}\right)^{\!\!{\nu -4\over \nu}} 
\min\left\{ {M^{{7\over 8}} (\log T)^{{969\over 5600}} 
\over T^{{3\over 20}} H^{{29\over 40}}} \;,\, {B_3 M^2\over T^{{2\over 3}} H^{{1\over 3}}}\right\}\;.
$$ 
A calculation shows that this last inequality is satisfied if one has the upper bound 
on $H$ in (6.11), with $B_4 = B_1^{12/13} B_3^{3/7},$ (say). Therefore, subject to 
(6.10) and (6.11) holding (with $B_4$ as just stated), we find that by applying 
Lemma~12.2, with $\widetilde N$ and $\widetilde R$ given by (12.43) and (12.37), and with 
$E'=16/9$, one is able to obtain the bounds 
$$
\align 
{|S|\over M^{1+\varepsilon}} &\leq 
\cases
\left( {16\over 9}\right)^{{100\over 221}} C_{II}(\varepsilon) 
\left( \left( \displaystyle{H\over \widetilde R}\right)^{\!\!{3\over 34}} 
\left( \displaystyle{H\over \widetilde N}\right)^{\!\!{7\over 24}} 
+ \left( \displaystyle{H\over \widetilde R}\right)^{\!\!{100\over 221}} 
\left( \displaystyle{H\over \widetilde N}\right)^{\!\!{17\over 26}}\right) 
 & \text {if $\,\displaystyle{2^{10} C_3^2 H^2 \over  \widetilde R \widetilde N }\leq 1$},\\  
\left( {16\over 9}\right)^{{31\over 102}} C_{I}(\varepsilon , \nu) 
\left( \displaystyle{H\over \widetilde R}\right)^{\!\left( {22\over 17}\right) q_{\nu}^{-1}} 
\left( \displaystyle{H\over \widetilde N}\right)^{\!\!{1\over 2}}  & \text{ otherwise},
\endcases\cr  
 &\quad \cr 
&\leq \cases \left( {4\over 3}\right) C_{II}(\varepsilon) 
\left( \left( \displaystyle{H^3 T\over M^3}\right)^{\!\!{3\over 68}} 
\left( \displaystyle{H\over \widetilde N}\right)^{\!\!{101\over 408}} 
+ \left( \displaystyle{H^3 T\over M^3}\right)^{\!\!{50\over 221}} 
\left( \displaystyle{H\over \widetilde N}\right)^{\!\!{189\over 442}}\right) 
 &\text {if $\,\displaystyle{2^{21} C_3^3 H^4 T\over M^{3}}\leq \widetilde N$,} \\ 
\left( {4\over 3}\right) C_{I}(\varepsilon , \nu) 
\left( \displaystyle{H^3 T\over M^3}\right)^{\!\left( {11\over 17}\right) q_{\nu}^{-1}} 
\left( \displaystyle{H\over \widetilde N}\right)^{\!\!{1\over 2}-\left( {11\over 17}\right) q_{\nu}^{-1}}    &
\text{otherwise}.
\endcases
\endalign
$$ 
By these bounds, in which ${1\over 2}-\left( {11\over 17}\right) q_{\nu}^{-1}>{1\over 2} -{11\over 68}>0$,  while $\widetilde N$ is as stated in (12.43), we are able to conclude that, subject to 
the conditions (6.10) and (6.11) both being satisfied, one will have   
$$
S\ll_{\varepsilon} \, \left( 
T^{13\over 160} M^{125\over 192} H^{179\over 320} +
T^{32\over 153} M^{19\over 51} H^{283\over 612}  
+T^{151\over 520} M^{-{11\over 208}} H^{1473\over 1040} +
T^{113\over 221} M^{-{118\over 221}} H^{276\over 221} \right) B_3^{-{1\over 2}} M^{\varepsilon}  \eqno{(12.44)}
$$ 
if 
$$
H\leq\min\left\{ 
\left( 128 C_3\right)^{-{40\over 63}} M^{155\over 189} T^{-{46\over 189}} (\log T)^{323\over 8820} \,,\,  
\left( 128 C_3\right)^{-{9\over 13}} M^{15\over 13} T^{-{5\over 13}}\right\}\;,\eqno{(12.45)}
$$
and will otherwise have  
$$
S\ll_{\varepsilon , \nu} \, \left( T^{\left( {11\over 20}\right) q_{\nu}^{-1}+{3\over 40}} 
M^{{9\over 16}-\left( {11\over 8}\right) q_{\nu}^{-1}} 
H^{\left( {33\over 40}\right) q_{\nu}^{-1} +{69\over 80}} 
+ T^{\left( {11\over 51}\right) q_{\nu}^{-1} +{1\over 3}} 
M^{-\left( {11\over 17}\right) q_{\nu}^{-1}} 
H^{\left( {55\over 51}\right) q_{\nu}^{-1} +{2\over 3}}\right) 
B_3^{-{1\over 2}} M^{\varepsilon}\;.\eqno{(12.46)}
$$  

Since $B_1$ and $B_3$ are both positive constants constructed from $C_2,\ldots ,C_6$, 
the results of Theorem~2 are an immediate consequence of the combination of Corollary~6.1.1,  
Corollary~6.1.2 and our conclusions reached in (12.42) and (12.44)-(12.46): note, in particular, 
that (6.4)-(6.6) imply 
$$M\leq B_5^{{1\over 10}} T^{{527\nu-1360\over 32(29\nu -75)}} (\log T)^{-{57(55\nu -144)\over 896(29\nu -75)}}\;,$$ 
while (6.10) and (6.11) imply $M\leq C_6 T^{1/2}$,  
and so (given that $\nu\geq 6$, that $0<B_5,B_4\leq 1\leq C_6$, and that $T$ is large) it follows that neither Part~(A) nor Part~(B) of Theorem~2 will apply unless one has 
$M\leq C_6 T^{901/1584}$, so that $M^{\varepsilon}\ll_{\varepsilon}\,T^{\varepsilon}$. 
This completes our proof that Theorem~2 is valid when the condition (7.1) is satisfied;  
given what was noted below (7.1), it has therefore been shown that Theorem~2  
is valid in all cases$.\quad\blacksquare$   

\bigskip
\beginsection
{13. Applications to the mean square of $|\zeta({1\over 2}+it)|$}  

\medskip
\noindent
{\bf Theorem~3.}
{\sl  Let the function $I : [0,\infty)\times (0,\infty)\rightarrow{\Bbb R}$ be given by: 
$$
I(t,U)={1\over 2U} \int_{t-U}^{t+U} \left|\zeta\left(\textstyle{1\over 2}+i\tau\right)\right|^2\,{\text d}\tau\;.
$$  
Suppose that $\varepsilon$ is a positive constant. 
Then one has 
$$
I\left( t,t^{{1273\over 4053}+\varepsilon}\right) =O(\log t)\quad\ {\text as}\ \,t\rightarrow\infty\;.
$$ }
\medskip

\noindent
{\bf Theorem~4.}
{\sl  Let the function $E(T)$ be defined on the interval $[1,\infty)$ by: 
$$E(T)=\int_0^T \left|\zeta\left(\textstyle{1\over 2}+it\right)\right|^2\,{\text d}t 
-\left(\log\left( {T\over 2\pi}\right) +2\gamma -1\right) T\qquad\hbox{($T\geq 1$),}$$ 
where $\gamma$ denotes the Euler-Mascheroni constant.
Suppose that $\varepsilon$ is a positive constant. Then one has 
$$
E(T) = O\left( T^{{1515\over 4816}+\varepsilon}\right)\quad\,{\text as}\ \,T\rightarrow\infty\;.
$$} 

\bigskip 

\noindent
{\bf Remarks.}\ 

\smallskip

\noindent
{\bf (i)} \ The proofs of these two theorems can be found at the end of this section. They depend crucially on Lemma~13.1 (below), which is a corollary of Theorem~2, and involve the use of three further lemmas 
(one aiding the proof of Theorem~4, while the other two aid the proof of Theorem~3). 

\smallskip 

\noindent
{\bf (ii)}\  Note that Theorem~4 improves upon the estimate $E(T)=O(T^{131/416} (\log T)^{32587/8320})$ 
obtained in [W10]: for one has $131/416=0.314903...\ $, while $1515/4816=0.314576...\ $. Theorem~4 is, moreover, 
quite close to the conditional bound $E(T)=O_{\varepsilon}(T^{\varepsilon +39/124})=O_{\varepsilon}(T^{\varepsilon +0.314516...\ })$ which, as is noted in [W10], would hold subject to the validity of the case $\kappa =1/4$, $\lambda =0$ of the `Hypothesis $H(\kappa ,\lambda )$' of Huxley [H03], [H05]. 
\par 
To have included, in this paper, a proper discussion of the (conditional) consequences of Huxley's Hypothesis $H(\kappa ,\lambda)$ would have led to an unwanted degree of complexity in our results and their proofs: we have (in any case) nothing certain to report regarding progress on this matter. It may nevertheless be 
worth mentioning that, on the basis of certain calculations, 
we do expect that, subject to the validity of the hypothesis $H(133/457,0)$, the number 
$1273/4053=0.314088...\ $ occurring in Theorem~3 might be replaced by $2811/8951=0.314043...\ $ (this 
would require using also the methods of the present paper). It is more complicated to determine what consequences of this sort would follow from the validity of the hypothesis $H(\kappa , \lambda)$ in cases where $133/457>\kappa\geq 1/4$, and we have done no work on that.

\smallskip 

\noindent
{\bf (iii)}\  In Lemma~13.1 (below) the bound $S\ll M/\log T$ is obtained whenever one has $U\geq T^{\varepsilon +1273/4053}\,$ (with some constant $\varepsilon >0$), whereas the stronger bound 
$S\ll UH/\log T$ is obtained only when $U$ satisfies the more restrictive condition $U\geq T^{\varepsilon +1515/4816}$. This is the reason for the differing exponents, $\varepsilon +1273/4053$ 
and $\varepsilon +1515/4816$, that occur in Theorem~3 and Theorem~4, respectively. 
The cases within the proof of Lemma~13.1 that are crucial in determining the limit $1515/4816$  (in (13.8)) 
are, perhaps surprisingly, not those in which $H=M T^{O(\varepsilon)-1515/4816}$: they are instead those cases 
in which one has either 
$$T^{O(\varepsilon )+{1041\over 2408}} \leq M\ll T^{{1\over 2}}\quad\,{\text and}\quad\, 
H=M T^{O(\varepsilon )-{389\over 1204}}\;,$$ 
or 
$$
T^{O(\varepsilon )+{1041\over 2408}}\geq M\gg T^{{807\over 1940}}\quad\,{\text and}\quad\, 
H=M^{{425\over 489}} T^{O(\varepsilon )-{42\over 163}}\;.
$$ 
The corresponding cases in Section~11 fall within the scope of Case~II (which is defined in that section): they lie along the boundary that separates those cases within Case~II that are best dealt with 
by the application of Proposition~10$'$ from those cases in which a stronger bound on $S$ is obtained 
by appealing instead to either Corollary~6.1.1 or Corollary~6.1.2. Since these crucial cases are 
quite far from being in Case~I, which is the only case in which one is left with a 
free choice of $\nu\in\{ 6,7,8,\ldots\ \}$, it seems likely that, in our proof of 
the case $c>1515/4816$ of Lemma~13.1, we might have been able to put $\nu$ equal to 
an arbitrary element of the set $\{ 6 , 7, 8 , \ldots\ \}$, instead of making 
the specific choice $\nu =7$ indicated in (13.9). However (as we hope is made clear by Remark~(iii) below Theorem~2) we really do need to put $\nu =7$ in our proof of the case $c\leq 1515/4816$ of Lemma~13.1. 

\smallskip

\noindent
{\bf (iv)} \ Given that $\zeta(\overline{s})=\overline{\zeta(s)}$ for all $s\in{\Bbb C}-\{ 1\}$, it is a direct consequence of the definitions in Theorem~3 and Theorem~4 that, when $U=t^c$, one has: 
$$
I(t,U)=\cases 
\displaystyle{E(t+U)-E(t-U)\over 2U} +{1\over 2U}\int_{t-U}^{t+U}\left( \log\left( {x\over 2\pi}\right) +2\gamma\right) 
{\text d}x &\text {if $0<c<1$ and $t\geq 2^{1/(1-c)}$;}\cr
 &\quad \cr  
\displaystyle{E(2U)\over 2U}+\log\left( {U\over\pi}\right) +2\gamma -1 
&\text{if $c=1$ and $t\geq 1/2$;} \cr
 &\quad \cr 
\displaystyle{E(U+t)+E(U-t)\over 2U}+\log\left( {U\over 2\pi}\right) +2\gamma -1 +O\left( {t^2\over U^2}\right) 
&\text {if $c>1$ and $t\geq 2^{1/(c-1)}$}.
\endcases
$$ 
By this it is readily be seen that Theorem~4 contains (i.e. implies immediately) those cases of  
Theorem~3 in which $\varepsilon$ exceeds the difference 
between $1273/4053=0.314088...\ $ and $1515/4816=0.314576...\ $, which is 
$0.00048808...\ <2^{-11}$.

\medskip
\noindent
{\bf Lemma~13.1.}
{\sl  Let $\delta_0$ and $c$ be constants satisfying 
$0<\delta_0 <1$ and 
$$
{1273\over 4053}<c<{1\over 3}\;.\eqno{(13.1)}
$$ 
Suppose that $T$ is a large positive parameter, that $U\in{\Bbb R}$ satisfies  
$$
T^c\leq U\leq 3 T^c\;,\eqno(13.2)
$$
and that $H,H_1,M,M_1\in (0,\infty)$ satisfy:  
$$
{H\over 2}\leq H_1\leq H,\qquad {M\over 2}\leq M_1\leq M\;,\eqno(13.3)
$$ 
$$
UH\leq M\leq 4 T^{{1\over 2}}\;,\eqno{(13.4)}
$$ 
and either 
$$
M>\delta_0 T^{{1\over 2}}\;,\eqno{(13.5)}
$$ 
or else 
$$
H>{U^{{7\over 2}}\over T}\;.\eqno{(13.6)}
$$ 
Suppose moreover that $b\in {\Bbb Z}$, that $F : [1/3,3]\rightarrow{\Bbb R}$ is the function given by 
$$
F(x)=\log(x) - {b M^2 x^2 \over 8 T}\qquad\hbox{($1/3\leq x\leq 3$),}\eqno{(13.7)}
$$ 
and that $S=S_F(T;H,H_1;M,M_1)$ is the exponential sum defined in (7.2). 
Then one has: 
$$
S\ll\cases\displaystyle{UH\over\log T} &\text{ if $\,c>\displaystyle{1515\over 4816}\,$,}\cr 
 &\quad \cr 
\displaystyle{M\over\log T} &\text{otherwise}.
\endcases
\eqno{(13.8)}
$$} 
\medskip

\noindent{\bf Proof.}\ 
Since $S=0$ if $H<1$, we may assume throughout that $H\geq 1$. 
\par 
We shall complete this proof by showing that the bound (13.8) is a corollary of the case $\nu =7$ of 
the results of Theorem~2. As a first step towards this we verify that the sum $S$ is such that 
the relevant hypotheses of Theorem~2 are satisfied. Given (7.2) and (13.3), 
the present sum $S$ is similar in form to the sum $S$ occurring in Theorem~2: it corresponds to the special case in which the functions $g(x)$ and $G(x)$ of the theorem are the step functions defined on the interval $[1/2,1]$ by $g(x)=|(H_1 /H , \infty)\cap\{ x\}|$ and $G(x)=|(M_1 /M , \infty)\cap\{ x\}|$. 
By (13.7), we have the cases $r=3,4,5$ of (6.1) and the cases $r=3,4$ of (6.2) for any choice 
of $C_3$, $C_4$ and $C_5$ satisfying $C_r\geq (r-1)!  3^r\,$ ($r=3,4,5$). 
Before considering (6.3) and the case $r=2$ of the conditions (6.1) and (6.2) it should be 
noted that the definitions (7.2) and (13.7) imply that the sum $S$ depends on the integer $b$ 
only insofar as it depends on whether $b$ is even or odd (indeed, each term of the sum $S$ is of the form 
$\phi(h,m){\text e}(-bmh/2)$, where $\phi(h,m)$ is a factor that is independent of $b$). 
The integers $1$ and $325$ are both odd, whereas $0$ is an even number. Therefore we may assume 
that either it is the case that $M\leq T^{1/2} /3$ and $b\in\{ 0 , 1\}$ or else 
it is the case that $M> T^{1/2} /3$ and $b\in\{ 0 , 325\}$.
In either of these two cases one has 
$$
F^{(2)}(x)=-\left( {1\over x^2} + {b M^2 \over 4T}\right)
$$ 
and 
$$
F^{(2)}(x)F^{(4)}(x)-3 F^{(3)}(x)^2 
= \left( {1\over x^2} + {b M^2 \over 4T}\right)\left( {6\over x^4}\right) 
-3\left( {2\over x^3}\right)^{\!\!2} 
={3b M^2 \over 2T x^4}-{6\over x^6}\;,
$$ 
for all $x$ lying in the interval $[1/3,3]$. 
Given that we have $0\leq b\leq 325$ and (by (13.4)) $\,M^2 /T\leq 16$, 
it follows that we have the case $r=2$ of both (6.1) and (6.2) 
for any choice of $C_2$ satisfying $C_2\geq 1309$. 
If $b=0$ then we have (6.3) for any choice of $C_5 \geq 3^5 /2$. 
If $b=1$ and $M\leq T^{1/2} /3$, then we have (6.3) for any choice of $C_5 \geq (4/3)(3^5 /2)=162$. 
In the remaining cases, where $b=325$ and $M>T^{1/2} /3$, we have (6.3) for any choice of 
$C_5 \geq  54 T /M^2$, and so for any $C_5 \geq 486$. We conclude that, in all 
the cases under consideration,  
the conditions (6.1)-(6.3) of Theorem~2 will hold  
if one puts $C_r = (r-1)!  3^r\,$ ($r=3,4,5$) and $C_2 =1309$. It follows that the 
results of Part~(A) and Part~(B) of Theorem~2 will be applicable to the sum $S=S_F(T;H,H_1 ;M,M_1 )$, 
provided only that it can be shown that the relevant additional conditions (i.e. (6.4)-(6.6) for Part~(A); 
(6.10)-(6.11) for Part~(B)) are satisfied. 
\par 
Given the upper bound on $M$ in (13.4), we may assume that the condition (6.5) is satisfied (for we shall have $T^{9/16} (\log T)^{-57/448}>4 T^{1/2}$, provided only that $T$ is large enough); the same bound on $M$ trivially implies that the condition (6.10) will be satisfied if we put $C_6 =4$. We choose now to put: 
$$
\nu =7\eqno(13.9)
$$ 
and $q=q_{\nu}$, with $q_{\nu}$ as defined in Theorem~2, so that 
$$
q=q_7={158\over 37}=4+{10\over 37}=4.\dot 27\dot 0\;.\eqno{(13.10)}
$$
By (13.9), the condition (6.6) will be satisfied if and only if  
$$
H\leq B_5 M T^{-{643\over 2048}} (\log T)^{{969\over 40960}}\;.\eqno{(13.11)}
$$ 
Since our hypotheses in (13.1), (13.2) and (13.4) imply 
that we have $H/M\leq 1/U\leq T^{-c}$, where $c$ is a constant satisfying $c>1273/4053=0.31408...\ >643/2048=0.31396...\ $, 
it follows that the condition (13.11) will be satisfied if $T$ is large enough (in terms of the 
small positive constant $B_5$).  We may therefore assume that (13.11) does hold, so that the 
case $\nu =7$ of the condition (6.6) is satisfied. 
\par 
To complete the data concerning our application of Theorem~2 (the implications of which 
are discussed below) we now specify $\varepsilon$ by putting  
$\varepsilon =\varepsilon(c)$, where 
$\varepsilon(c)$ is equal to $(c-(1515/4816))/200$
if $c>1515/4816$,  and is otherwise equal to  
$(c-(1273/4053))/5$.  
Given (13.1), this ensures that $\varepsilon$ is a constant satisfying the following three conditions: 
$$
\eqalignno{ 0<\varepsilon &\leq \left( \textstyle{1\over 3}\right)\left( c-\textstyle{1515\over 4816}\right)\quad\,{\text if}\ \,c>\textstyle{1515\over 4816}\;, &(13.12)\cr 
 &\quad\cr 
0<\varepsilon &\leq \left( \textstyle{1\over 3}\right)\left( c-\textstyle{1273\over 4053}\right)\;, &(13.13)\cr 
&\quad\cr 
0<\varepsilon &\leq 0.0001\;. &(13.14)}
$$
\par 
Our next steps depend on whether or not it is the case that the first inequality occurring in (6.4) is satisfied. Suppose, firstly, that one does have 
$$
H\geq M^{-9} T^4 (\log T)^{{171\over 140}}\;.\eqno{(13.15)}
$$ 
Then, recalling the points noted in the previous paragraph, we are able to conclude that 
all three of the conditions ((6.4), (6.5) and (6.6)) attached to Part~(A) of Theorem~2 are satisfied, and 
so it follows by (13.9), (13.10) and Part~(A) of that theorem that either  
$$
S\ll H \left( {H\over M}\right)^{\!\!{103\over 3950}} T^{{1273\over 3950}+\varepsilon}\;,\eqno{(13.16)}
$$ 
or else one has the bounds stated in (6.7) and (6.8). 
By (13.2), (13.4), (13.13) and (13.4) (again), the bound (13.16) would imply  
$$
S\ll {H T^{{1273\over 3950}+\varepsilon} \over U^{{103\over 3950}} } 
=\left( {T^{{1273+3950\varepsilon\over 4053}} \over U}\right)^{\!\!{4053\over 3950}} U H 
\leq {U H \over T^{\left( {4053\over 3950}\right)\left( \left( c -{1273\over 4053}\right) -\varepsilon\right)}} 
< {UH\over T^{2\varepsilon}}\leq {M\over T^{2\varepsilon}}\;,
$$
and so would yield the result (13.8) of the lemma. If we do not have (13.16), then we have instead (6.7) and (6.8), which imply that one has   
$$
{S\over H}\ll T^{\varepsilon} \min\{ X , Z\} + T^{\varepsilon}Y\;,
$$ 
with 
$$
X = \left({M\over H}\right)\left({H\over M}\right)^{\!\!{323\over 600}} 
T^{397\over 2400}  
\ll \left({M\over H}\right)\left( {(\log T)^{969\over 64960}\over T^{{149\over 464}}}\right)^{\!\!\!{323\over 600}} T^{{397\over 2400}} 
= {M (\log T)^{{104329\over 12992000}}\over H T^{{2705\over 278400}}} 
\ll {M\over T^{2\varepsilon} H}\;,
$$
$$
Y = \left( {H\over M}\right)^{\!\!{19\over 50}} T^{1133\over 2600} 
\ll \left( {(\log T)^{969\over 64960}\over T^{{149\over 464}}}\right)^{\!\!\!{19\over 50}} T^{1133\over 2600} 
=T^{{3785\over 12064}} (\log T)^{{18411\over 3248000}}
$$
and $Z=(H/M)^{1/25} T^{131/400}$, so that one obtains:  
$$
\min\{ X , Z\}\leq X^{{24\over 301}} Z^{{277\over 301}} 
=T^{{1515\over 4816}}\;.
$$ 
These bounds would imply, firstly, that 
$$
S\ll H T^{\varepsilon} X + H T^{\varepsilon} Y 
\ll {M\over T^{\varepsilon}} + H T^{{3785\over 12064}+2\varepsilon} 
<{M\over T^{\varepsilon}} + H T^{c-\varepsilon}\leq {M+H U\over T^{\varepsilon}}
\leq {2M\over T^{\varepsilon}}
$$ 
(with the latter part of this following by virtue of (13.14), (13.1), (13.2) and (13.4), given that  
$3785/12064<0.31375$, whereas $1273/4053>0.31408$) and, secondly, that  one has:
$$
{S\over H}\ll T^{{1515\over 4816}+\varepsilon} + T^{{3785\over 12064}+2\varepsilon} 
\leq T^{c-\varepsilon}\leq {U\over T^{\varepsilon}}\quad\ \hbox{if $\,c>{1515\over 4816}$}
$$ 
(with the inequality in the middle following by (13.1), (13.12) and (13.14)). 
Therefore, in the event 
that (6.7) and (6.8) hold,  we obtain (13.8). Since we have found that (13.8) is obtained whether or not the bound (13.16) holds, this completes our proof in respect of the cases in which the condition (13.15) is satisfied. 
\par 
Suppose now that (13.15) does not hold, so that one has 
$$
H< M^{-9} T^4 (\log T)^{{171\over 140}}\;.\eqno{(13.17)}
$$ 
This implies $M<H^{-1/9} T^{4/9} (\log T)^{19/140}\leq T^{4/9} (\log T)^{19/140}$. 
Therefore, provided that $T$ is sufficiently large in terms of the small positive constant $\delta_0$,  
it will be the case that the inequality (13.5) does not hold; this is, by hypothesis, incompatible with it simultaneously being the case that (13.6) does not hold, and so we may henceforth assume that the inequality 
(13.6) is satisfied. 
\par 
The remainder of this proof rests on the application of Part~(B) of Theorem~2. 
Our first task, therefore, is to verify that the condition (6.11) is satisfied 
(it already having been noted that (13.4) gives (6.10)). We recall that the inequality 
(13.11) was found to hold (and that this was subject only to the hypotheses of the lemma). 
Therefore (13.11) holds in the present case, and by it and (13.17), we may deduce that 
$$\eqalign{H <\left( {B_5 M (\log T)^{{969\over 40960}}\over T^{{643\over 2048}}}\right)^{\!\!{4096\over 4215}}\left( {T^4 (\log T)^{{171\over 140}}\over M^9}\right)^{\!\!{119\over 4215}} 
 &=B_5^{{4096\over 4215}}\left( {M^{{4096-1071\over 5}}  
(\log T)^{{969\over 50}+{2907\over 100}}\over T^{{1286-476\over 5}}}\right)^{\!\!{1\over 843}} \cr 
 &=B_5^{{4096\over 4215}}\left( {M^{605} (\log T)^{{969\over 20}}\over T^{162}}\right)^{\!\!{1\over 843}} \;.}
$$
Now Theorem~2 would remain valid if one substituted the constant 
$B_5'=\min\{ B_5 , B_4^{4215/4096}\}$ in place of the constant $B_5$ in (6.6)  
(indeed, this would either have no effect, or would slightly 
weaken the content of the theorem). We may therefore assume that the 
constants $B_5$ and $B_4$ in Theorem~2 satisfy $B_5^{4096/4215}\leq B_4$. 
By this, and the upper bound on $H$ that was just obtained (above), we find that the case $\nu =7$ of the condition 
(6.11) will be satisfied if it is the case that one has:  
$$
H\leq B_4\left( {M^3\over T}\right)^{\!\!{23\over 55}}\;.
$$ 
To see that this does hold, we recall that (13.6) was shown to hold (subject to our 
assumption (13.17)), and note that, by the combination of the first inequality in (13.4) with (13.6), (13.2) and (13.1), one has 
$$
\left( {1\over H}\right)\left( {M^3\over T}\right)^{\!\!{23\over 55}} 
\geq \left( {1\over H}\right)\left( {U^3 H^3\over T}\right)^{\!\!{23\over 55}} 
={U^{{69\over 55}} H^{{14\over 55}}\over T^{{23\over 55}}} 
>{U^{{69\over 55}} \left( U^{{7\over 2}} T^{-1}\right)^{\!\!{14\over 55}}\over T^{{23\over 55}}} 
={U^{{118\over 55}}\over T^{{37\over 55}}}\geq T^{\left( {118\over 55}\right)\left( c-{37\over 118}\right)} 
>{1\over B_4}
$$
(with the final inequality following from the relations $c>1273/4053=0.3140...\ >37/118=0.3135...\ $, 
provided that $T$ is sufficiently large in terms of the small positive constant $B_4$).
This completes our verification of the condition (6.11). 
\par 
Since both (6.10) and (6.11) are satisfied, 
it follows by (13.9), (13.10) and Part~(B) of Theorem~2 that either 
$$
S\ll T^{{161\over 790}+\varepsilon} M^{{19\over 79}} H^{{417\over 395}} 
+T^{{1031\over 2686}+\varepsilon} M^{-{407\over 2686}} H^{{2469\over 2686}}\;,\eqno{(13.18)}
$$ 
or else one has the bounds stated in (6.12) and (6.13). 
Here we observe that one has:
$$
{T^{{161\over 790}+\varepsilon} M^{{19\over 79}} H^{{417\over 395}}\over 
H (H/M)^{{103\over 3950}} T^{{1273\over 3950}+\varepsilon}} 
=\left( {M^{1053} H^{117}\over T^{468}}\right)^{\!\!{1\over 3950}} 
=\left( {M^9 H\over T^4}\right)^{\!\!{117\over 3950}} 
\ll T^{\varepsilon}
$$
(the last inequality following by virtue of our assumption (13.17)). 
By this and the first calculation appearing below (13.16), we find that, given the hypothesis (13.1), it 
must follow from (13.18) that one has: 
$$
{S\over H}\ll {U\over T^{\varepsilon}} +T^{{1031\over 2686}+\varepsilon} 
M^{-{407\over 2686}} H^{-{217\over 2686}}\;.
$$ 
Recall now that the assumption (13.17) ensured that (13.6) must hold. By (13.4) and (13.6), we have 
$H>U^{7/2} T^{-1}$ and $M\geq UH>U^{9/2} T^{-1}$. By these inequalities, the above bound for 
$S/H$ and the hypothesis (13.2), we find that 
$$
{S\over H} 
\ll {U\over T^{\varepsilon}} 
+ \left( {T\over U^{9/2}}\right)^{\!\!{407\over 2686}} \left( {T\over U^{7/2}}\right)^{\!\!{217\over 2686}} 
T^{{1031\over 2686}+\varepsilon} 
={U\over T^{\varepsilon}} 
+\left( {T^{{1655\over 5277}}\over U}\right)^{\!\!\!{5277\over 2686}} T^{\varepsilon} U
\ll  {U\over T^{\min\left\{ \varepsilon\,,\,\left( {5277\over 2686}\right)\left( c-{1655\over 5277}\right)-\varepsilon\right\} }}\;.
$$ 
Therefore, given that $1655/5277<0.3137$, whereas $1273/4053>0.314$, the hypothesis 
(13.1) and inequality (13.14) are enough to ensure that we obtain the bound 
$S\ll UH T^{-\varepsilon}$. By this and the first inequality of (13.4), we may conclude that 
the bound in (13.8) holds when one has both (13.17) and (13.18). 
\par 
The only cases that remain to be considered are those in which one has (13.17) and (instead of (13.18)) the bounds (6.12) and (6.13). By (6.12), we have $M\gg T^{1/3} H^{13/15}$. This leads to bounds for 
a couple of the terms occurring on the right-hand side of (6.13): 
$$
T^{113\over 221} M^{-{118\over 221}} H^{276\over 221}\ll T^{{1\over 3}} H^{{2606\over 3315}}\quad 
{\text and}\quad 
T^{79\over 204} M^{-{11\over 68}} H^{191\over 204} 
\ll T^{{1\over 3}} H^{{203\over 255}} =T^{{1\over 3}} H^{{2639\over 3315}}\;.\eqno{(13.19)}
$$
Note the greater magnitude of the latter bound: regarding it, we observe that, since the inequality 
(13.6) holds, one  has 
$$
{T^{{1\over 3}} H^{{203\over 255}}\over H} 
=T^{{1\over 3}} H^{-{52\over 255}} 
<T^{{1\over 3}} \left( {T\over U^{{7\over 2}} }\right)^{\!\!{52\over 255}} 
=U\left( {T^{{137\over 437}}\over U}\right)^{\!\!\!{437\over 255}} 
\leq {U\over T^{\left( {437\over 255}\right)\left( c-{137\over 437}\right) }}\leq {U\over T^{\varepsilon}} 
\eqno{(13.20)}
$$
(the last two inequalities following by virtue of (13.2), (13.1) and (13.14), given that one has $137/437<0.3136$, 
whereas $1273/4053>0.314$). 
\par 
By (13.17) and (6.12), we have also $M^9 H<T^4 (\log T)^{171/140}$ and 
$M^{-155} H^{189}\ll T^{-46} (\log T)^{969/140}$, and so:
$$
M^{-{11\over 208}} H^{433\over 1040} 
=\left( M^9 H\right)^{{709\over 24128}} \left( {H^{189}\over M^{155}}\right)^{\!\!{19\over 9280}} 
\ll \left( T^4 (\log T)^{{171\over 140}}\right)^{\!{709\over 24128}} 
\left( {(\log T)^{{969\over 140}}\over T^{46}}\right)^{\!\!\!{19\over 9280}}\;.
$$
By this bound, together with (13.1), (13.2) and (13.14), it follows that we have  
$$
{T^{151\over 520} M^{-{11\over 208}} H^{1473\over 1040}\over H} 
\ll T^{{3785\over 12064}+\varepsilon} <T^{c-\varepsilon}\leq {U\over T^{\varepsilon}}\eqno{(13.21)}
$$
(note what this has in common with the bound at the end of the paragraph containing (13.15)). 
\par 
By (13.19), (13.20) and (13.21), the bound (6.13) for $S$ implies:
$$
{S\over H} 
\ll {U\over T^{\varepsilon}} 
+T^{\varepsilon}\min\left\{ X_1\;,\, Z_1 \right\} 
+T^{\varepsilon}\min\left\{ Y_1\;,\, Z_1 \right\}\;, \eqno{(13.22)}
$$
where 
$$
X_1 = T^{13\over 160} M^{125\over 192} H^{-{141\over 320}} \;,\qquad 
Y_1 = T^{32\over 153} M^{19\over 51} H^{-{329\over 612}} \quad {\text and}\quad 
Z_1 = T^{17\over 80} M^{7\over 32} H^{{11\over 160}}\;.\eqno{(13.23)}
$$
We note, firstly, that it follows from (13.22) and the first inequality in (13.4) that one has 
$$
{S\over M}\ll \left( {1\over T^{\varepsilon}}\right)\left( 1 + 
\left( {H\over M}\right) X_1 T^{2\varepsilon} + \left( {H\over M}\right) Y_1 T^{2\varepsilon} \right) \;.\eqno{(13.24)}
$$
By (13.23), one has 
$$
\left( {H\over M}\right) X_1 = {T^{13\over 160} H^{{179\over 320}}\over M^{{67\over 192}} } 
=\left( {H T^{{78\over 537}}\over M^{{335\over 537}}}\right)^{\!\!\!{537\over 960}}\;.
$$ 
Moreover, the bound (6.12) and assumption (13.17) imply that one has here 
$$
H\ll\left( {M^{155\over 189} (\log T)^{323\over 8820}\over T^{{46\over 189}} }\right)^{\!\!\!{20349\over 20764}} 
\left( {T^4 (\log T)^{{171\over 140}}\over M^9}\right)^{\!\!\!{415\over 20764}} 
= {M^{{335\over 537}} (\log T)^{{12521\over 207640}} \over T^{{4939\over 31146}} }\;,
$$
and so (given (13.14)): 
$$
\left( {H\over M}\right) X_1 T^{2\varepsilon} 
\ll \left( {(\log T)^{{12521\over 207640}} \over T^{{4939\over 31146}-{78\over 537}} } \right)^{\!\!\!{537\over 960}} 
T^{2\varepsilon} \ll T^{3\varepsilon -{83\over 11136}} < 1\;.\eqno{(13.25)}
$$
Regarding the final term in (13.24), we note that, since (6.12) implies 
$M\gg T^{1/3} H^{13/15}$, it follows from (13.23) that one has
$$
\left( {H\over M}\right) Y_1 T^{2\varepsilon} 
={T^{{32\over 153}+2\varepsilon}  H^{{283\over 612}} \over M^{{32\over 51}} } 
\ll {T^{{32\over 153}+2\varepsilon}  H^{{283\over 612}} \over 
\left(T^{{1\over 3}} H^{{13\over 15}}\right) ^{\!\!{32\over 51}} } 
= {T^{2\varepsilon}\over H^{{249\over 3060}}}\;,\eqno{(13.26)}
$$
where, since (13.6) holds, one has $H>U^{7/2} T^{-1}\geq T^{(7/2)c-1}$, with 
$(7/2)c-1>115/1158>153/2490\,$ (by (13.1)). 
By (13.24), (13.25), (13.26) and (13.14), we find that 
$$
S\ll {M\over T^{\varepsilon}}\;.\eqno{(13.27)}
$$ 
\par 
Note that (13.22) also implies the bound 
$$
{S\over H}\ll {U+X_2 T^{2\varepsilon} +Y_2 T^{2\varepsilon} \over T^{\varepsilon}} \;,\eqno{(13.28)}
$$ 
where 
$$
X_2 = X_1^{{24\over 301}} Z_1^{{277\over 301}} \quad {\text and}\quad 
Y_2 = Y_1^{{408\over 5723}} Z_1^{{5315\over 5723}}\;,
$$
so that, by (13.23) and (13.17), one has: 
$$
X_2 =\left( T^{{973\over 16}} M^{{2439\over 32}} H^{{271\over 32}}\right)^{\!{1\over 301}} 
=\left( M^9 H\right)^{{271\over 9632}} T^{{139\over 688}} 
<\left( T^4 (\log T)^{{171\over 140}}\right)^{{271\over 9632}} T^{{139\over 688}} 
\ll T^{{1515\over 4816}+\varepsilon}
$$
and, similarly,  
$$
Y_2 = \left( T^{{58309\over 48}} 
M^{{42069\over 32}} H^{{14023\over 96}} \right)^{\!{1\over 5723}} 
<\left( T^4 (\log T)^{{171\over 140}}\right)^{\!{14023\over 549408}} T^{{58309\over 274704}} 
\ll T^{{28785\over 91568}+\varepsilon}=T^{{1515\over 4816+(64/19)}+\varepsilon}\;.
$$ 
By (13.28), the above bounds for $X_2$ and $Y_2$, and (13.2) and (13.12),  we have 
$$
S\ll \left( {U H\over T^{\varepsilon} }\right)\left( 1 +{X_2 T^{2\varepsilon} + Y_2 T^{2\varepsilon}\over T^c} \right)  
\ll \left( {U H\over T^{\varepsilon} }\right) 
\left( 1+{T^{3\varepsilon}\over T^{c-{1515\over 4816}}}\right)\ll {UH\over T^{\varepsilon}} 
\quad\ {\text if}\ \,c>{1515\over 4816}\;.
$$
This bound on $S$, together with 
that in (13.27), imply what is stated in (13.8), and so complete the proof.$\quad\blacksquare$

\medskip
\noindent
{\bf Lemma~13.2.}
{\sl Let $t,U>0$ satisfy $(t\pi)^{1/2}\leq U^2\leq t/(2\pi)$. Let $I(t,U)$ be as 
defined in Theorem~3, above. Then either it is the case that 
$$
I(t,U)\ll \log t\;, \eqno{(13.29)}
$$ 
or else there exists some 
$$
T\in\left[ {t\over 4\pi} \,,\, {3t\over 4\pi}\right] \;,\eqno{(13.30)}
$$
some 
$$
M\in\left[ 1 \,,\, (2T)^{{1\over 2}}\right]\;,\eqno{(13.31)}
$$ 
and some 
$$
M_1\in \left[ {M\over 2} \,,\, M\right]\;,\eqno{(13.32)}
$$ 
such that the sum 
$$
S_{*}\left( T ; U ; M , M_1\right) 
=\sum_{0<h\leq \left( e^{1/U} -1\right) M/2}\ \,\sum_{M_1 <m\leq M} \left( {m+h\over m-h}\right)^{\!\!2\pi i T}\;, 
\eqno{(13.33)}
$$ 
satisfies 
$$
{\left| S_{*}\left( T ; U ; M , M_1\right)\right| \over M} \gg {I(t,U)\over \log t}\;.\eqno{(13.34)}
$$} 

\noindent
{\bf Proof.}\ 
By the case $r=1$ of the results contained in [W04, Lemma~10.2 and Lemma~10.3] it follows that, 
for some  $\tau$ lying in the interval $[t/2 , 3t/2]$, one has 
$$
I(t,U)\ll\biggl(\,\sum_{0<m\leq K_0} {1\over m}\biggr) + t^{{1\over 2}} U^{-2} + U^2 t^{-{3\over 2}} 
+\sum_{j=0}^{\infty} \left| E\left( U; \tau ; K_j \right)\right| \;,
$$ 
where 
$$
K_j = 2^{-j} \left( {t\over 2\pi}\right)^{\!\!{1\over 2}}\quad\,{\text and}\quad\, 
E(U ; \tau ; K) =\sum_{{K\over 2} < k\leq K} {1\over k} \sum_{1\leq d\leq \left( e^{1/U} -1\right) K/2} 
\left( {k+d\over k-d}\right)^{\!\!-i\tau}\;.
$$ 
Given the hypotheses of the lemma concerning $t$ and $U$, it follows that we have the bound 
$$
{I(t,U)\over\log t}\ll 1 + {1\over |{\cal J}|}\sum_{j\in{\cal J}} \left| E\left( U; \tau ; K_j \right)\right|\;,
$$ 
where ${\cal J}=\{ j\in{\Bbb Z} \;:\; 1\leq 2^j\leq K_0\}$.  
Therefore, either it is the case that the relation (13.29) holds, or else we must have 
$$
{I(t,U)\over\log t}\ll \left| E\left( U; \tau ; K_{j_{*}} \right)\right|\quad\,{\text for\ some}\ \, 
j_{*}\in{\cal J}\;.
$$ 
With regard to the latter of these two cases, we observe that, by partial summation and the 
invariance of the absolute value under complex conjugation, it follows that if $K>0$ then one has 
$$
E\left( U; \tau ; K \right) \ll {1\over K}\left|\sum_{K'< k\leq K} 
\ \sum_{1\leq d\leq \left( e^{1/U} -1\right) K/2} \left( {k+d\over k-d}\right)^{\!i\tau}\right| \quad\,{\text for\ some}\ \, 
K'\in\left[ {K\over 2} , K\right]\;.
$$
By applying this with $K=K_{j_{*}}$ and then substituting 
$2\pi T$, $M$, $M_1$, $m$ and $h$ for $\tau$, $K_{j_{*}}$, $K_{j_{*}}'$, $k$ and $d\,$ (respectively),  
we obtain what is described in (13.30)-(13.34).$\quad\blacksquare$ 
\medskip

\noindent
{\bf Lemma~13.3.}
{\sl  Let $U\geq 1$. Suppose that $T,M,M_1>0$ satisfy the conditions (13.31) and (13.32). 
Put  
$$
D(T , U , M)=\min\left\{ \left( {e^{1/U} -1\over 2}\right) M \;,\, {U^{{7\over 2}}\over T}\right\}\;.\eqno{(13.35)}
$$ 
Then one has 
$$
\sum_{0<h\leq D(T , U , M)}\ \sum_{M_1 <m\leq M} \left( {m+h\over m-h}\right)^{\!\!2\pi i T} 
\ll M\;.\eqno{(13.36)}
$$} 

\noindent
{\bf Proof.}\ 
Note firstly that, for $h\in{\Bbb N}$, one has 
$$
\sum_{M_1 <m\leq M} \left( {m+h\over m-h}\right)^{\!\!2\pi i T} 
=\sum_{M_1 <m\leq M} {\text e}\left( f_h (m)\right) = W_h \quad\hbox{(say),}
$$ 
where $f_h(x)=T\log (x+h) - T\log (x-h)$. One can show moreover that, when 
$h/M$ is sufficiently small (in absolute terms), the exponential sum $W_h$ may be estimated through 
an application of the theory of exponent pairs: one will then obtain, in particular, the 
bounds 
$$
W_h\ll \left( {hT\over M^2}\right)^{{2\over 7}} M^{{4\over 7}} +{M^2\over hT} 
=(hT)^{{2\over 7}} +O\left( {1\over h}\right)\ll (hT)^{{2\over 7}}\;,\eqno{(13.37)}
$$ 
which derive from the exponent pair $BAAB(0,1) =(2/7 , 4/7)\,$  
(i.e. we are here applying the case $k=2/7$, $l =4/7$, $s=2$, $y=2hT$, $N=M/2$ of 
the result stated in [G\&K91, Equation~(3.3.4)]). 
Since (13.35) implies $D(T,U,M)\leq (e^{1/U} -1)M/2$, where we have 
$0< e^{1/U}-1=e^{1/U}-e^0<((1/U)-0)e^{1/U}\leq e/U$, it follows that if  
$U$ is sufficiently large (in absolute terms) then, by (13.37) and (13.35) (again), one will have
$$
\sum_{0<h\leq D(T , U , M)} \left| W_h\right| 
\ll \sum_{0<h\leq D(T , U , M)} (hT)^{{2\over 7}} 
\ll \left( D(T , U , M)\right)^{{9\over 7}} T^{{2\over 7}} 
\leq \left( {e^{1/U} -1\over 2}\right) \!M \left( {U^{{7\over 2}}\over T}\right)^{\!\!{2\over 7}} \!T^{{2\over 7}} 
\leq {e M\over 2}\;,
$$
and so will obtain the result stated in (13.36). 
\par 
The only cases of the lemma requiring further 
proof are those in which one has $U\leq U_0$, with $U_0$ equal to a certain positive absolute constant. 
We note that, by (13.35) and the hypotheses concerning $U$, $T$ and $M$,  one has 
$D(T , U , M)\leq U^{7/2}/T \leq 2 U^{7/2}$. 
Since this trivially implies the bound   
$$
\sum_{0<h\leq D(T , U , M)} \left| W_h\right|\leq D(T , U , M) M\ll U^{{7\over 2}} M\;,
$$ 
we therefore find that (13.36) holds when $U$ is less than or equal to the absolute constant $U_0.\quad\blacksquare$ 

\medskip
\noindent
{\bf Lemma~13.4.}
{\sl  Let $c$ be a constant satisfying 
$$
{1515\over 4816} < c < {1\over 3}\;.\eqno{(13.38)}
$$
The, for all $t , \Delta\geq 1$ such that 
$$
t^c \log t \leq \Delta \leq (2t)^c \log(2t)\;,\eqno{(13.39)}
$$ 
the sum 
$$
G^{+}(t , \Delta) 
=\sum_{\scriptstyle mn\leq t/2\pi\atop\scriptstyle 0<\Delta\log(m)-\Delta\log(n)\leq\log t}
{(m/n)^{it}\over\sqrt{mn}\log(m/n)}
\,\exp\left( -\left( \textstyle{1\over 2}\,\Delta\log(m/n)\right)^{2}\right) \eqno{(13.40)}
$$
satisfies 
$$
G^{+}(t , \Delta) \ll (\log t)^2 \Delta\;.\eqno{(13.41)}
$$} 

\noindent
{\bf Proof.}\ 
It may be assumed that $t$ satisfies a condition of the form $T\geq T_0$, where $T_0$ denotes an arbitrarily large positive constant: for the bound (13.41) is trivial when one has $1\leq t\ll 1$. 
In particular we may assume that  $t\geq T_0 >(8\pi /3)^7$ and, given (13.39) and (13.38) (in which 
$1515/4816>2/7$), may also assume that 
$$
(t/4)^{{1\over 3}}\geq \Delta\log t\quad\,{\text and}\quad\,{\Delta\over\log t}\geq t^{{2\over 7}}\geq\log t 
\geq\log T_0 > 7 \log {8\pi \over 3}\;.\eqno{(13.42)}
$$ 
This justifies the application of the bound for $G^{+}(t , \Delta)$ that is noted in 
[W10, Equation~(6.7)]. From that bound it follows that either (13.41) holds, or else one has 
a bound of the form 
$$
G^{+}(t , \Delta) \ll (\log t)^2 H^{-1} \left| S_F \left( T ; H , H_1 ; M , M_1 \right)\right| 
=(\log t)^2 H^{-1} |S|\quad\ \hbox{(say),}\eqno{(13.43)}
$$ 
where $S=S_F ( T ; H , H_1 ; M , M_1 )$ is as defined in (7.2), while $F : [1/3 , 3]\rightarrow{\Bbb R}$ 
is the function given by (13.7), with $b$ some constant that is equal to either $0$ or $1$,  and 
$(T , H , H_1 , M , M_1 )$ is some point of ${\Bbb R}^5$ such that $T$, $H$, $H_1$, $M$ and $M_1$ 
satisfy: 
$$
\left| T-{t\over 2\pi}\right|\leq {\Delta\log t\over 2\pi}\leq {t\over 8\pi}\eqno{(13.44)}
$$
and 
$$
{\Delta H\over \log T}\leq M\leq 4T^{{1\over 2}}\;,\eqno{(13.45)}
$$ 
as well as the conditions in (13.3) and either  the inequality 
$$
H>{(\Delta \log T )^{{7\over 2}} \over T}\;,\eqno{(13.46)}
$$ 
or else a condition of the form (13.5) in which $\delta_0$ is   
equal to a certain positive absolute constant.
In cases where (13.41) holds there is nothing further to prove. Therefore 
it may henceforth be assumed that the function $F$ and real parameters $T$, $H$, $H_1$, $M$ and $M_1$ 
fit the description just given, and are, moreover, such that the relation $(13.43)$ holds. 
\par 
With the application of Lemma~13.1 in mind, we put 
$$
U={\Delta \over\log T}\;,
$$ 
so that, by (13.45), the condition (13.4) is satisfied. 
The condition (13.1) is implied by (13.38), and since   
$t>(8\pi /3)^7$, it follows by (13.38), (13.39) and (13.44) that the condition 
(13.2) is satisfied also. By the assumptions we have made, 
the inequalities in (13.3) are satisfied, the function $F(x)$ is as stated in (13.7), and 
either (13.5) holds, or else we have (13.46). 
Moreover, if (13.46) holds, then (given we have (13.44) and $t>(8\pi /3)^7$) 
it implies 
$$
H>\left( {U^{{7\over 2}}\over T}\right) (\log T)^7 
>\left( {U^{{7\over 2}}\over T}\right) \left( {6 \log t\over 7}\right)^{\!\!7} 
>\left( {U^{{7\over 2}}\over T}\right) \left( 6\log{8\pi\over 3}\right)^{\!\!7} \;,
$$ 
and so gives the inequality (13.6). 
Therefore we may conclude that Lemma~13.1 applies, so that the bound (13.8) is obtained. 
By (13.8), (13.38) and (13.44), we have 
$$
S_F \left( T ; H , H_1 ; M , M_1 \right)\ll {UH\over \log T}={\Delta H\over (\log T)^2} 
\ll {\Delta H\over (\log t)^2}\;.
$$ 
By this, (13.43) and (13.42), we find that one has $G^{+}(t,\Delta)\ll\Delta\ll (\log t)^2 \Delta.\quad\blacksquare$

\bigskip 

\noindent
{\bf The proof of Theorem~3.}\ 
We put $c={1273\over 4053}+\varepsilon$ and $U=t^c$. 
In view of Remark~(iv) following the statements of Theorem~3 and Theorem~4 (at the beginning of this section), 
it will be enough to consider only those cases in 
which one has $0<\varepsilon\leq 2^{-11}$, and so (given that $1273/4053<3^{-1}-2^{-6}$) we may certainly assume 
that $c$ satisfies the inequalities in (13.1). Since we have only to bound 
$I(t,t^c)$ for all sufficiently large positive values of $t$, we may certainly assume also that 
$t\geq (2\pi)^3$, that $U\geq (2\pi)^{3c}>1$, and that 
any $T$ satisfying (13.30) will (by virtue of the implied inequality $T>t/4\pi$) 
certainly be large enough to permit the application of Lemma~13.1 (should all the other hypotheses of that lemma happen to be satisfied). Then, given that $1273/4053=0.3140...\ >5/16=0.3125$, 
it follows by (13.1), Lemma~13.2 and Lemma~13.3 that either it is the case that 
$$
I\left( t , t^c\right)=I(t,U)\ll \log t\;,\eqno{(13.47)}
$$ 
or else, for some $(T,M,M_1)\in{\Bbb R}^3$ satisfying (13.30), (13.31) and (13.32), one has 
$$
{I\left( t , t^c \right)\over\log t}={I(t,U)\over\log t}\ll {1\over M} 
\left| \sum_{U^{7/2} T^{-1} <h\leq\left( e^{1/U}-1\right) M/2}\ \,\sum_{M_1 <m\leq M} 
\left( {m+h\over m-h}\right)^{2\pi iT}\right| \eqno{(13.48)}
$$ 
(which would imply also that $M>2(e^{1/U }-1)^{-1} U^{7/2} T^{-1}$). 
Only the latter of these two cases requires further consideration: for the validity of the 
bound in (13.47) is what we are seeking to establish in this proof. 
Accordingly, we note that by splitting the sum in (13.48) at points where 
$h\in\{ [2^{-j}(e^{1/U} -1)M] \;:\; j\in{\Bbb N}\}$, and then applying the triangle inequality and the principle that the arithmetic mean of $N$ real numbers will not exceed the greatest of those numbers, 
it may be deduced that, for some $(H,H_1)\in{\Bbb R}^2$ satisfying both 
$$
{U^{{7\over 2}}\over T} < H \leq\left( {e^{1/U} -1\over 2}\right) M\;,\eqno{(13.49)}
$$
$$
{H\over 2}\leq H_1\leq H\;,\eqno{(13.50)}
$$ 
one has: 
$$
{\left| S_L \right|\over M}  
\gg {I(t,U)\over (\log t)^2} = {I\left( t , t^c \right)\over (\log t)^2}\;,\eqno{(13.51)}
$$ 
where $L$ denotes the function $L(x)=\log(x)\,$ ($1/3\leq x\leq 3$) and $S_L$ is 
the sum $S_L ( T ; H , H_1 ; M , M_1 )$ 
that is given by the case $F(x)=L(x)$ of (7.2). 
\par 
Based on observations made earlier, we have $U\geq (2\pi)^{3c}>(2\pi)^{15/16}>2>1/\log 2$, and so 
the rightmost inequality if (13.49) will imply $H\leq (1/U)\exp( (1/U) -\log 2) M < M/U$. 
This, together with the fact that $T$, $M$, $M_1$, $H$ and $H_1$ satisfy (13.31), (13.32), (13.49) and (13.50), is enough to ensure that (13.3), (13.4) and (13.6) hold. Since $c$ satisfies (13.1) and $T$ satisfies (13.30), we have also $T^{-c}U=(t/T)^c\in[(4\pi /3)^c , (4\pi)^c]$, in which  
$(4\pi /3)^c >1$ and $(4\pi)^c < (4\pi)^{1/3} <3$. Consequently we find that (13.2) holds as well (as does (13.7), when $F(x)=L(x)=\log(x)$ and $b=0$). Since we are assured of having $T$ be sufficiently large for 
Lemma~13.1 to apply, it therefore follows by that lemma that the bound (13.8) is obtained when $S=S_L$, so that we must have $S_L\ll (\log T)^{-1}\max\{ UH , M\}\ll M/(\log T)\,$ (with the final inequality holding by virtue of it being the case that the condition (13.4) is satisfied). Since (13.30) holds, and since we assume that $t\geq (2\pi)^3$, we have also $\log T>\log(t) - 2\log(2\pi)\geq (1/3)\log(t)$, and so 
may deduce from the bound just obtained for $S_L$ that one has $|S_L|/M\ll 1/\log t$. 
By this and (13.51), we obtain the desired estimate (13.47).$\quad\blacksquare$ 

\bigskip 

\noindent{\bf The proof of Theorem~4.}\
We put $c={1515\over 4816}+{\varepsilon \over 2}$. It will suffice to consider only cases in which 
$\varepsilon$ lies in the interval $(0 , 1/28]\,$ (say), and so we may certainly assume that 
the constant $c$ satisfies the inequalities in (13.38). 
\par 
Suppose now that $T$ satisfies $T>(8/(1-3c))^{12/(1-3c)}\,$ (for example). We then put 
$\Delta = T^c \log T$. By (13.38) and our supposition concerning $T$ it follows that we have 
$1\leq\Delta\leq (T/4)^{1/3}$, and so, as an immediate corollary of the estimates 
contained in [H\&H90, Lemma~8.1], we find that one has 
$$
E(T)-E(T/2)=O\left( \Delta (\log T)^2\right) +G(T) -G(T/2)\;,\eqno{(13.52)}
$$ 
where $G(t)$ denotes a certain sum that is defined in  [H\&H90, Lemma~8.1]: 
one can show, in particular, that 
$$
G(t)= 4 {\text Im}\left( G^{+}(t , \Delta)\right) +O(1)\qquad\quad\hbox{($t\in\{ T , T/2\}$),}\eqno{(13.53)}
$$ 
where $G^{+}(t , \Delta)$ is the sum defined in (13.40) 
(to show this requires essentially nothing more than the properties of complex conjugation and elementary bounds for certain of the terms occurring in the sum $G^{+}(t , \Delta)$). 
By our choice of $\Delta$, the condition (13.39) is satisfied for $t=T$, and also for $t=T/2$. 
Therefore (given that the condition (13.38) is also satisfied) we obtain from Lemma~13.4 the 
upper bound (13.41) for $t=T$, and also for $t=T/2$, and so it follows, by (13.52) and (13.53), that we have: 
$$
E(T)-E(T/2)\ll \Delta (\log T)^2 =T^c (\log T)^3 =T^{{1515\over 4816}+{\varepsilon\over 2}} (\log T)^3 
\ll T^{{1515\over 4816}+\varepsilon}\;.
$$
>From this estimate, and the formula for the sum of a geometric series, we may infer that one has 
$$
E(T)-E\left( T/2^j\right)\ll T^{{1515\over 4816}+\varepsilon}\;,
$$ 
where $j>0$ is the least integer such that $T/2^j\leq (8/(1-3c))^{12/(1-3c)}$. 
To complete the proof we have only to observe that one has $E(U)\ll 1\ll U^{{1515\over 4816}+\varepsilon}$ 
for all $U$ satisfying $1\leq U\leq (8/(1-3c))^{12/(1-3c)}\,$ 
(this being a trivial corollary of the elementary fact that, since $|\zeta(1/2+it)|^2\geq 0$ for all real $t$,   one must have  
$|E(U)|\leq |E(V)|+2(\log(2\pi V)+2\gamma -1)V$ whenever $1\leq U\leq V$).$\quad\blacksquare$ 

\bigskip

\beginsection References.
\par 
\frenchspacing 

\item{[B]} J.~Bourgain, {Decoupling, exponential sums and the Riemann 
Zeta Function}, arXiv:1408.5794.

\item{[BCT]} J.~Bennett, A.~Carbery, T.~Tao, 
{ On the multilinear restriction and Kakeya conjectures}, {\it Acta
Math.} 19 (2006), no 2, 261--302.

\item
{[B-D]} J.~Bourgain, C.~Demeter, {The proof of the $\ell^2$-decoupling
conjecture}, arXiv:1403.5335, (to appear in Annals of Math).

\item
{[B-G]} J.~Bourgain, L.~Guth, {Bounds on oscillatory integral operators
based on multilinear estimates}, {\it GAFA 21}
(1211), no 6, 1239-1295

\item{[G\&K91]} S.~W. Graham and G. Kolesnik, { Van der Corput's Method of Exponential Sums}, 
{\it London Mathematical Society Lecture Notes}  126 (Cambridge University Press, Cambridge, 1991). 

\item{[H\&H90]} D.~R. Heath-Brown and M.~N. Huxley, Exponential sums with 
a difference, {\it Proc. London Math. Soc.} (3) 61 (1990), 227-250.

\item{[H96]} M.~N. Huxley, {Area, Lattice Points and Exponential Sums}, 
{\it London Mathematical Society Monographs New Series} 13 (Oxford
University Press, Oxford, 1996).

\item{[H03]} M.~N. Huxley, Exponential sums and lattice points III, {\it Proc. London Math. Soc.} (3) 87 (2003), 591-609.

\item{[H04]} M.~N Huxley, Resonance curves in the Bombieri-Iwaniec method, {\it Funct. Approx. Comment. Math.} 32 (2004), 7-49.

\item{[H05]} M.~N. Huxley, Exponential sums and the Riemann zeta function V', {\it Proc. London Math. Soc.} (3) 90 (2005), 1-41.

\item{[W04]} N. Watt, On the mean squared modulus of a Dirichlet $L$-function over a short segment of the critical line, {\it Acta Arith.} 111 (2004), 307-403. 
\item{[W10]} N.~Watt, A note on the mean square of $|\zeta({1\over 2}+it)|$, {\it J. London Math. Soc.} 82 (2) (2010), 279-294. doi: 10.1112/jlms/jdq024

\enddocument